\documentclass[10pt,reqno]{amsart}

\usepackage{amsthm,amsmath,amstext,amssymb,amscd,euscript, mathrsfs, dsfont,multicol,times,enumerate,subfig,sidecap}
\usepackage{enumerate}
\usepackage{amsmath, amssymb, amsthm}
\numberwithin{equation}{section}
\usepackage{mathrsfs}
\usepackage{esint}
\usepackage{xcolor}
\usepackage{mathtools}
\usepackage{bm}
\usepackage{bbm}

\usepackage[colorlinks=true, pdfstartview=FitV, linkcolor=blue, citecolor=red, urlcolor=blue]{hyperref} 

\newtheorem{thm}{Theorem}[section]
\newtheorem{corollary}[thm]{Corollary}
\newtheorem{lem}[thm]{Lemma}
\newtheorem{prop}[thm]{Proposition}
\theoremstyle{definition}
\newtheorem{defn}[thm]{Definition}
\newtheorem{example}[thm]{Example}
\newtheorem{assumption}[thm]{Assumption}
\newtheorem{rem}[thm]{Remark}

\newcommand\bN{\mathbb{N}}

\newcommand\bR{\mathbb{R}}

\newcommand\cF{\mathcal{F}}

\newcommand\cpar{\text{$[$\kern-.38em$[$}}
\newcommand\cbrk{\text{$]$\kern-.15em$]$}}
\newcommand\opar{\text{\,\raise.2ex\hbox{${\scriptstyle
|}$}\kern-.34em$($}}
\newcommand\obrk{\text{$)$\kern-.34em\raise.2ex\hbox{${\scriptstyle |}$}}\,}

\def\XXint#1#2#3{{\setbox0=\hbox{$#1{#2#3}{\int}$ }
\vcenter{\hbox{$#2#3$ }}\kern-.58\wd0}}

\makeatletter
\def\@tocline#1#2#3#4#5#6#7{\relax
  \ifnum #1>\c@tocdepth 
  \else
    \par \addpenalty\@secpenalty\addvspace{#2}%
    \begingroup \hyphenpenalty\@M
    \@ifempty{#4}{%
      \@tempdima\csname r@tocindent\number#1\endcsname\relax
    }{%
      \@tempdima#4\relax
    }%
    \parindent\z@ \leftskip#3\relax \advance\leftskip\@tempdima\relax
    \rightskip\@pnumwidth plus4em \parfillskip-\@pnumwidth
    #5\leavevmode\hskip-\@tempdima
      \ifcase #1
       \or\or \hskip 1em \or \hskip 2em \else \hskip 3em \fi%
      #6\nobreak\relax
    \dotfill\hbox to\@pnumwidth{\@tocpagenum{#7}}\par
    \nobreak
    \endgroup
  \fi}
\makeatother

\begin{document}

\title[Stable-like Operators related to Bernstein Functions in weighted Lebesgue spaces]{Stable-like Operators related to Bernstein Functions in weighted Lebesgue spaces: Boundedness, Regularity, Embedding, and Martingale Problem}

\author[J.-H. Choi]{Jae-Hwan Choi}
\address{Department of Mathematics, The Catholic University of Korea, 43 Jibong-ro, Bucheon-si, Gyeonggi-do, 14662, Republic of Korea}
\email{mathjhchoi@catholic.ac.kr}

\author[J. Kang]{Jaehoon Kang}
\address{Department of Applied Mathematics and Institute for Integrated Mathematical Sciences, Hankyong National University, 327 Jungang-ro, Anseong-si, Gyeonggi-do 17579, Republic of Korea} 
\email{jaehoon.kang@hknu.ac.kr}

\author[D. Park]{Daehan Park}
\address{Department of Mathematics, Kangwon National University, 1 Kangwondaehakgil, Chucheon-si, Gangwon-do, 24341, Republic of Korea} 
\email{daehanpark@kangwon.ac.kr}

\thanks{}



\begin{abstract}
We study variable-coefficient nonlocal operators associated with a Bernstein function $\phi$ satisfying global weak scaling conditions.
Let $j_d^\phi$ denote the radial L\'evy density of the corresponding subordinate Brownian motion.
The kernels may be nonsymmetric, and the coefficients are measurable in the jump variable and satisfy a prescribed uniform modulus of continuity in the spatial variable. 
In the critical scaling regime, we impose an annular cancellation condition.

For jump profiles satisfying $0\le V(y)\le Cj_d^\phi(|y|)$, we prove that the operator $\mathcal L_V^a$ is bounded from $H_p^{\phi,2}(w)$ to $L_p(w)$ for every $1<p<\infty$ and $w\in A_p(\mathbb R^d)$.
For $V(y)=j_d^\phi(|y|)$ and uniformly positive real coefficients, we establish maximal regularity for parabolic equations in mixed-norm spaces with Muckenhoupt weights in space and time. 
The parabolic coefficients may be merely measurable in time. 
We also obtain elliptic resolvent estimates in weighted $L_p$ spaces for sufficiently large resolvent parameters.

Under suitable integrability conditions, we prove a weighted Sobolev--H\"older embedding with an explicit local modulus of continuity in space and time. 
The modulus depends on $\phi$, the integrability exponents, and the Muckenhoupt indices of the weights.

Finally, we prove global well-posedness of the time-dependent martingale problem. 
For the resulting time-inhomogeneous strong Markov family, we derive a Krylov-type estimate in weighted mixed norms. 
Its local normalization involves the weighted measures of an intrinsic time interval and a spatial ball.
\end{abstract}

\maketitle
\tableofcontents

\section{Introduction and Main Results}

\subsection{Motivation and goals}

A basic motivation for the class of nonlocal operators studied in this paper comes from jump-type stochastic differential equations with state-dependent coefficients. 
Let $L=(L_t)_{t\ge0}$ be a pure-jump L\'evy process on $\mathbb{R}^d$ with L\'evy measure $\nu(\mathrm{d}y)=V_L(y)\,\mathrm{d}y$.
The corresponding translation-invariant nonlocal operator is
$$
\mathcal L_L f(x)
=
\int_{\mathbb{R}^d}
\left(
f(x+y)-f(x)-\nabla f(x)\cdot y\,\mathbf 1_{\{|y|\le1\}}
\right)
V_L(y)\,\mathrm{d}y .
$$
Now consider the stochastic differential equation
\begin{equation}\label{eq:motivation_sde}
    \mathrm{d}X_t=\sigma(X_{t-})\,\mathrm{d}L_t .
\end{equation}
The coefficient $\sigma(x)$ changes a jump $y$ of the driving process into the state-dependent jump $\sigma(x)y$. 
Formally, the generator associated with \eqref{eq:motivation_sde} is therefore
$$
\mathcal L_X f(x)
=
\int_{\mathbb{R}^d}
\left(
f(x+\sigma(x)y)-f(x)
-\nabla f(x)\cdot \sigma(x)y\,\mathbf 1_{\{|y|\le1\}}
\right)
V_L(y)\,\mathrm{d}y .
$$
If $\sigma(x)$ is invertible, then the change of variables $z=\sigma(x)y$ gives
$$
\mathcal L_X f(x)
=
\int_{\mathbb R^d}
\left(
f(x+z)-f(x)
-\nabla f(x)\cdot z\,
\mathbf 1_{\{|\sigma(x)^{-1}z|\le1\}}
\right)
V_L(\sigma(x)^{-1}z)
|\det \sigma(x)^{-1}|\,\mathrm{d}z .
$$
When $V_L$ is symmetric, the difference between the above truncation and the standard truncation $\mathbf 1_{\{|z|\le1\}}$ contributes only an odd drift term, which vanishes. 
Hence, the operator can be written as
$$
\mathcal L_X f(x)
=
\int_{\mathbb R^d}
\left(
f(x+z)-f(x)
-\nabla f(x)\cdot z\,\mathbf 1_{\{|z|\le1\}}
\right)
V_L(\sigma(x)^{-1}z)
|\det \sigma(x)^{-1}|\,\mathrm{d}z .
$$
In particular, if $L$ is a subordinate Brownian motion associated with a Bernstein function $\phi$, then by denoting $V_{L}(y) = j^{\phi}_{d}(|y|)$, the preceding operator takes the form
\begin{equation}\label{eq:motivation_operator}
\mathcal L_X f(x)
=
\int_{\mathbb R^d}
\left(
f(x+z)-f(x)
-\nabla f(x)\cdot z\,\mathbf 1_{\{|z|\le1\}}
\right)
a_\sigma(x,z)j_d^\phi(|z|)\,\mathrm{d}z ,
\end{equation}
where
\begin{equation}\label{eq:motivation_coefficient}
a_\sigma(x,z)
:=
\frac{
j_d^\phi(|\sigma(x)^{-1}z|)
|\det \sigma(x)^{-1}|
}{
j_d^\phi(|z|)
}.
\end{equation}
Thus, even when the driving noise is spatially homogeneous, a state-dependent coefficient naturally produces a nonlocal operator with a spatially inhomogeneous jump kernel. 
This motivates the study of L\'evy-type operators whose jump intensity depends on the current state and whose jump kernel is comparable to that of a subordinate Brownian motion.

For a general coefficient $a(x,y)$, which need not be of the special form \eqref{eq:motivation_coefficient}, the construction or identification of an underlying Markov process associated with \eqref{eq:motivation_operator} is a delicate problem. 
In probabilistic terms, this is closely related to the well-posedness of the martingale problem for the corresponding integro-differential operator. 
From the analytic point of view, it is equally natural to ask whether such operators admit a weighted regularity theory.

Motivated by these considerations, we study nonlocal operators of the form
$$
\mathcal L_\phi^a u(x)
:=
\int_{\mathbb R^d}
\left(
u(x+y)-u(x)-\nabla u(x)\cdot y\,\chi_{\delta_0,\delta_1}(y)
\right)
a(x,y)j_d^\phi(|y|)\,\mathrm{d}y ,
$$
where $j_d^\phi$ is the jumping kernel associated with the Bernstein function $\phi$, and $\chi_{\delta_0,\delta_1}$ is a indicator function determined by the weak scaling exponents $\delta_{0},\delta_{1}$ of $\phi$. 
The coefficient $a=a(x,y)$ is assumed to be bounded and measurable in the jump variable $y$, and to have a prescribed modulus of continuity in the spatial variable $x$. 
This class includes the operators arising from \eqref{eq:motivation_sde} under natural nondegeneracy and regularity assumptions on $\sigma$, but also allows substantially more general spatially inhomogeneous and nonsymmetric perturbations.

The purpose of this paper is to develop a weighted regularity theory for $\mathcal L_\phi^a$ and to connect this analytic theory with probabilistic properties of the associated jump processes. 
The main results can be summarized as follows.

\begin{enumerate}[1.]
    \item We prove \emph{the weighted boundedness of the nonlocal operator}.
    More precisely, for $w\in A_p(\mathbb R^d)$ and $1<p<\infty$, we show that $\mathcal L_\phi^a$ is bounded from the weighted Sobolev space $H_p^{\phi,2}(w)$ to $L_p(w)$. 

    \item Using this boundedness, we establish \emph{weighted maximal regularity estimates}  of elliptic and parabolic equations involving nonlocal operators.
    For parabolic equations, we use mixed-norm spaces with Muckenhoupt weights in both the space and time variables.

    \item Besides the weighted regularity estimates, we establish a \emph{weighted Sobolev--Morrey embedding} adapted to the nonlocal operator generated by $\phi$. 
    We prove that functions in the corresponding weighted Sobolev space admit locally H\"older-type continuous representatives, with a modulus determined by $\phi$. 
    This embedding is then used to convert weighted Sobolev regularity into pointwise resolvent estimates, which is crucial for the martingale problem argument and the Krylov-type bounds.

    \item Finally, we derive probabilistic consequences from the weighted regularity theory. 
    In the real-valued, uniformly nondegenerate setting, we prove \emph{the well-posedness of the time-dependent martingale problem associated with the family of operators $\{\mathcal L_{\phi,t}^a\}_{t\ge0}$}, or equivalently with the space-time operator $\partial_t+\mathcal L_{\phi,t}^a$. 
    We provide \emph{Krylov-type occupation bounds} for the corresponding time-inhomogeneous Markov family.
\end{enumerate}

\subsection{Setting}

Before stating the main results, we fix the notation and assumptions used throughout the paper.

\subsubsection*{Bernstein functions and the associated nonlocal operators}

We first provide the definition and properties of Bernstein functions and the nonlocal operator $\phi(\Delta)$.

\begin{defn}\label{def 02.23.18:54}
A function $\phi:(0,\infty)\to[0,\infty)$ is called a \emph{Bernstein function} if $\phi$ is infinitely differentiable and
$$
    (-1)^n \phi^{(n)}(r)\le 0,
    \qquad r>0,\quad n\in\mathbb{N}.
$$
\end{defn}

\begin{rem}\label{rmk 12.30.16:40}
We recall some facts on Bernstein functions, which will be used in this paper.

\begin{enumerate}[(i)]
\item
If $\phi$ is a Bernstein function with $\phi(0+)=0$, then there exist a constant $b_\phi\ge0$ and a measure $\mu_\phi$ on $(0,\infty)$ satisfying
$$
    \int_0^\infty (1\wedge t)\,\mu_\phi(\mathrm{d}t)<\infty
$$
such that
\begin{equation}\label{eqn 12.30.16:03}
    \phi(r)
    =
    b_\phi r
    +
    \int_0^\infty
    \bigl(1-\mathrm{e}^{-rt}\bigr)\,\mu_\phi(\mathrm{d}t),
    \qquad r>0 .
\end{equation}
Conversely, any function of the form \eqref{eqn 12.30.16:03} is a Bernstein function. 
We call $b_\phi$ the drift of $\phi$ and $\mu_\phi$ the L\'evy measure of $\phi$.
We refer to \cite{schbern} for further properties of Bernstein functions.

\item
For every $n\in\mathbb{N}$, there exists a constant $C(n)>0$ such that
\begin{equation}\label{eqn 01.28.13:23}
    |r^n\phi^{(n)}(r)|
    \le C(n)\phi(r),
    \qquad r>0 .
\end{equation}

\item
For every $s\in(0,1)$, the function $\phi^s(r):=(\phi(r))^s$ is again a Bernstein function; see \cite[Corollary 3.8(iii)]{schbern}.
\end{enumerate}
\end{rem}

For $u\in\mathcal S(\mathbb R^d)$, we define $\phi(\Delta)$ by the Fourier multiplier
$$
    \mathcal F_d[\phi(\Delta)u](\xi)
    =
    -\phi(|\xi|^2)\mathcal F_d[u](\xi).
$$
Thus $-\phi(\Delta)$ is the positive operator with symbol $\phi(|\xi|^2)$. 
More generally, for $\gamma\in\mathbb R$, we define
\begin{equation}\label{26.06.05.11.35}
    \mathcal F_d[(1-\phi(\Delta))^{\gamma/2}u](\xi)
    =
    (1+\phi(|\xi|^2))^{\gamma/2}\mathcal F_d[u](\xi).
\end{equation}
For $\gamma>0$, we also write $\phi^\gamma(\Delta)$ for the operator with Fourier multiplier $-(\phi(|\xi|^2))^\gamma$, that is,
$$
    \mathcal F_d[\phi^\gamma(\Delta)u](\xi)
    =
    -(\phi(|\xi|^2))^\gamma \mathcal F_d[u](\xi).
$$
When only norms are involved, the sign convention is immaterial.

The operator $\phi(\Delta)$ is the infinitesimal generator of a subordinate Brownian motion. 
Indeed, let $S=(S_t)_{t\ge0}$ be a subordinator with Laplace exponent $\phi$,
$$
    \mathbb{E}[\mathrm{e}^{-\lambda S_t}]
    =
    \mathrm{e}^{-t\phi(\lambda)},
    \qquad \lambda>0,
$$
and let $B=(B_t)_{t\ge0}$ be a $d$-dimensional Brownian motion independent of $S$. 
Then $X_t=B_{S_t}$ is a subordinate Brownian motion with characteristic exponent $\phi(|\xi|^2)$.

We denote by $j_d^\phi$ the radial jump density of this process. 
With the above convention for Bernstein function $\phi$,
$$
    j_d^\phi(r)
    =
    \int_0^\infty
    (4\pi t)^{-d/2}
    \exp\left(-\frac{r^2}{4t}\right)
    \mu_\phi(\mathrm{d}t),
    \qquad r>0 .
$$
Also, the L\'evy measure of the subordinate Brownian motion $X$ is $j_d^\phi(|y|)\,\mathrm{d}y$, and for smooth functions $u$,
$$
\phi(\Delta)u(x)
=
b_\phi\Delta u(x)
+
\int_{\mathbb R^d}
\left(
u(x+y)-u(x)-\nabla u(x)\cdot y\,\mathbf 1_{\{|y|\le1\}}
\right)
j_d^\phi(|y|)\,\mathrm{d}y .
$$
In the model case $\phi(r)=r^{\alpha/2}$ with $\alpha\in(0,2)$, one has
$$
    \phi(\Delta)=\Delta^{\alpha/2},
    \quad
    j_d^\phi(r)=c_{d,\alpha}r^{-d-\alpha}.
$$

We impose the following weak scaling condition on $\phi$.

\begin{assumption}[$\delta_0,\delta_1$]\label{asm 07.21.13:06}
Let $\phi:(0,\infty)\to(0,\infty)$ be a Bernstein function. 
We assume that there exist constants $0<\delta_0\le\delta_1<1$ and $c_0,c_1>0$ such that
$$
    c_0\left(\frac Rr\right)^{\delta_0}
    \le
    \frac{\phi(R)}{\phi(r)}
    \le
    c_1\left(\frac Rr\right)^{\delta_1},
    \qquad 0<r\le R<\infty .
$$
\end{assumption}

Under Assumption \ref{asm 07.21.13:06}, the drift part vanishes, that is, $b_\phi=0$. 
Moreover, the jumping kernel $j^{\phi}_{d}(|y|)$ of $\phi(\Delta)$ satisfies
\begin{equation}\label{eqn 02.25.17:23}
    C^{-1}\frac{\phi(|y|^{-2})}{|y|^d}
    \le
    j_d^\phi(|y|)
    \le
    C\frac{\phi(|y|^{-2})}{|y|^d},
    \qquad y\neq0,
\end{equation}
where $C>0$ depends only on $d,c_0,c_1,\delta_0,\delta_1$.

\subsubsection*{Muckenhoupt weights}

We next introduce the class of weights used throughout the paper.

\begin{defn}
Let $1<p<\infty$. A nonnegative locally integrable function $w$ on
$\mathbb R^d$ is called a Muckenhoupt $A_p$ weight if
$$
    [w]_{A_p(\mathbb R^d)}
    :=
    \sup_B
    \left(\frac1{|B|}\int_B w(x)\,\mathrm{d}x\right)
    \left(
    \frac1{|B|}\int_B w(x)^{-1/(p-1)}\,\mathrm{d}x
    \right)^{p-1}
    <\infty,
$$
where the supremum is taken over all open balls $B\subset\mathbb R^d$. 
For a Muckenhoupt $A_{p}$ weight $w$, we write $w\in A_p(\mathbb R^d)$.
When $d=1$, we also write $w\in A_p(\mathbb R)$.
\end{defn}

For a measurable set $E\subset\mathbb R^d$, we write
$$
    w(E):=\int_E w(x)\,\mathrm{d}x .
$$

\begin{rem}\label{rmk 02.02.23:45}

\begin{enumerate}[(i)]
\item
For a locally integrable function $f$ on $\mathbb R^d$, the Hardy--Littlewood maximal function (or Hardy--Littlewood maximal operator) is defined by
$$
    \mathbb M f(x)
    :=
    \sup_{B\ni x}\frac1{|B|}\int_B |f(y)|\,dy .
$$
It is well-known that $\mathbb M$ is bounded on $L_p(\mathbb R^d,w(x)\,dx)$ if and only if $w\in A_p(\mathbb R^d)$; see \cite{M72}. 
For functions of one real variable, usually the time variable, we write $\mathcal M$ for the one-dimensional maximal operator.

\item
The class $A_p(\mathbb R^d)$ is invariant under translations, dilations, and orthogonal transformations. 
More precisely, if $w\in A_p(\mathbb R^d)$, $c>0$, $Q$ is an orthogonal matrix, and $y\in\mathbb R^d$, then
$$
    w_c(x):=w(cx),
    \qquad
    w_Q(x):=w(Qx),
    \qquad
    w_y(x):=w(x+y)
$$
belong to $A_p(\mathbb R^d)$ and
$$
    [w_c]_{A_p(\mathbb R^d)}
    =
    [w_Q]_{A_p(\mathbb R^d)}
    =
    [w_y]_{A_p(\mathbb R^d)}
    =
    [w]_{A_p(\mathbb R^d)}.
$$

\item
Let $p'=p/(p-1)$ and $w\in A_p(\mathbb R^d)$. 
Set $\bar w(x):=w(x)^{-1/(p-1)}$.
Then $\bar w\in A_{p'}(\mathbb R^d)$ and
$[\bar w]_{A_{p'}(\mathbb R^d)}=[w]_{A_p(\mathbb R^d)}^{1/(p-1)}$.
Moreover, the dual space of $L_p(w):=L_p(\mathbb R^d,w(x)\,dx)$ is identified with $L_{p'}(\bar w)$ through the pairing
$$
    \langle f,g\rangle
    =
    \int_{\mathbb R^d} f(x)g(x)\,\mathrm{d}x .
$$

\item
Let $w_0\in A_{p_0}(\mathbb R^d)$ and $w_1\in A_{p_1}(\mathbb R^d)$. 
If
$$
    \frac1p
    =
    \frac{1-\theta}{p_0}
    +
    \frac{\theta}{p_1},
    \qquad \theta\in(0,1),
$$
then $w:=w_0^{p(1-\theta)/p_0}w_1^{p\theta/p_1}
    \in A_p(\mathbb R^d)$.

\item \label{itm 11.03.12:46}
The classes $A_p$ are increasing in $p$:
$$
    A_p(\mathbb R^d)\subset A_q(\mathbb R^d),
    \qquad 1\le p<q<\infty,
$$
and for $w \in A_{p}(\mathbb{R}^{d})$, there exists $r<p$ such that $w \in A_{r}(\mathbb{R}^{d})$. We understand $A_1$ in the usual sence.
We set
\begin{equation}\label{eqn 06.29.12:14}
    A_\infty(\mathbb R^d)
    :=
    \bigcup_{1<p<\infty} A_p(\mathbb R^d).
\end{equation}
If $w\in A_p(\mathbb R^d)$ and $E\subset B$ is a measurable subset of a ball $B$, then
\begin{equation*}
    C^{-1}
    \left(\frac{|B|}{|E|}\right)^{(s-1)/s}
    \le
    \frac{w(B)}{w(E)}
    \le
    [w]_{A_p(\mathbb R^d)}
    \left(\frac{|B|}{|E|}\right)^p ,
\end{equation*}
where $s=s(d,p,[w]_{A_p(\mathbb R^d)})>1$ and $C=C(d,p,[w]_{A_p(\mathbb R^d)})>0$. 
In particular, $w$ is doubling: for every $\lambda\ge1$, $z\in\mathbb R^d$, and $r>0$,
\begin{equation}\label{eqn 01.30.01:21}
    w(B_{\lambda r}(z))
    \le
    C\lambda^{dp}w(B_r(z)),
\end{equation}
where $C$ depends only on $d,p$, and $[w]_{A_p(\mathbb R^d)}$.

\item\label{item 06.29.12:15}
If $w\in A_\infty(\mathbb R^d)$, we define $R_w:=\inf\{q\in[1,\infty): w\in A_q(\mathbb R^d)\}$.
Then $R_w<\infty$. 
If $w\in A_p(\mathbb R^d)$, then $R_w\le p$.
Moreover, for every $r\in(R_w,p]$, we have $w\in A_r(\mathbb R^d)$, and therefore
\begin{equation}\label{eqn 02.03.21:00}
    C^{-1}
    \left(\frac{|B|}{|E|}\right)^{(s_r-1)/s_r}
    \le
    \frac{w(B)}{w(E)}
    \le
    [w]_{A_r(\mathbb R^d)}
    \left(\frac{|B|}{|E|}\right)^r
\end{equation}
for every measurable set $E\subset B$, where $s_r=s_r(d,r,[w]_{A_r(\mathbb R^d)})>1$ and $C=C(d,r,[w]_{A_r(\mathbb R^d)})>0$.
\end{enumerate}

For more information on $A_p$ weights, we refer to \cite{Grafakos}.
\end{rem}

\begin{rem}\label{rmk 02.03.00:25}
Let $w\in A_\infty(\mathbb R^d)$. Then $R_w<\infty$ (recall \eqref{eqn 06.29.12:14}), and for every $r>R_w$ one has $w\in A_r(\mathbb R^d)$.
Consequently, for any $\delta_0,\gamma>0$, one can choose $p>R_w$ sufficiently large and then choose $p_{0}\in(R_w,p)$ so that
$$
    2\delta_0\gamma-\frac{dp_{0}}{p}>0.
$$
\end{rem}

\subsubsection*{Weighted Sobolev spaces associated with $\phi(\Delta)$}

\begin{defn}\label{def 02.23.18:55}
Let $1<p<\infty$, $\gamma\in\mathbb R$, $w\in A_p(\mathbb R^d)$, and let
$\phi$ be a Bernstein function. We define
$$
    H_p^{\phi,\gamma}(w)
    =
    H_p^{\phi,\gamma}(\mathbb R^d,w(x)\,\mathrm{d}x)
$$
as the completion of $\mathcal S(\mathbb R^d)$ with respect to the norm
$$
    \|u\|_{H_p^{\phi,\gamma}(w)}
    :=
    \|(1-\phi(\Delta))^{\gamma/2}u\|_{L_p(\mathbb R^d,w(x)\,\mathrm{d}x)} .
$$
When $w\equiv1$, we simply write $H_p^{\phi,\gamma}$. 
We also write $H_p^{\phi,0}(w)=L_p(w)$.
\end{defn}

If $\phi(r)=r$, then $\phi(\Delta)=\Delta$ and $H_p^{\phi,\gamma}(w)$ coincides with the usual weighted Bessel potential space. 
Thus,  $H_p^{\phi,\gamma}(w)$ can be viewed as the weighted Bessel potential spaces generated by the nonlocal operator $\phi(\Delta)$.
Moreover, by \eqref{eqn 01.28.13:23}, the multipliers appearing in \eqref{26.06.05.11.35} satisfy the usual Mikhlin-type estimates in weighted $L_p$ spaces. 
Hence $H_p^{\phi,\gamma}(w)$ shares many basic properties with classical Bessel potential spaces; see Remark \ref{rmk 01.28.14:36} below.

\subsubsection*{The nonlocal operator}

We now introduce the class of nonlocal operators that we study in this paper.
Let $V:\mathbb R^d\to[0,\infty)$ be a measurable function satisfying
\begin{equation}\label{eqn 02.11.14:19}
    V(y)
    \le
    C_V j_d^\phi(|y|),
    \qquad y\neq0 .
\end{equation}
For a bounded measurable coefficient $a:\mathbb R^d\times\mathbb R^d\to\mathbb C$, we define, initially for smooth functions $u$,
$$
\mathcal L_V^a u(x)
:=
\int_{\mathbb R^d}
\left(
u(x+y)-u(x)
-
\nabla u(x)\cdot y\,\chi_{\delta_0,\delta_1}(y)
\right)
a(x,y)V(y)\,\mathrm{d}y,
$$
where
$$
\chi_{\delta_0,\delta_1}(y)
=
\begin{cases}
0, & \delta_1<1/2\\[1mm]
\mathbf 1_{\{|y|\le1\}} & \delta_0\le1/2\le\delta_1,\\[1mm]
1 & \delta_0>1/2.
\end{cases}
$$
The model operator is obtained by taking $V(y)=j_d^\phi(|y|)$.
In this case, we also write $\mathcal L_\phi^a:=\mathcal L_{j_d^\phi}^a $.

In the critical regime $\delta_0\le1/2\le\delta_1$, we impose the following
cancellation condition on the kernel.

\begin{assumption}[Cancellation condition]
\label{24.07.02.14.41}
Assume that $\delta_0\le1/2\le\delta_1$. The coefficient $a:\mathbb R^d\times\mathbb R^d\to\mathbb C$ is bounded and measurable, and for every $x\in\mathbb R^d$ and every $0<r<R<\infty$,
$$
    \int_{r<|y|<R}
    y\,a(x,y)V(y)\,\mathrm{d}y
    =
    0 .
$$
In the model case $V(y)=j_d^\phi(|y|)$, this condition becomes
$$
    \int_{r<|y|<R}
    y\,a(x,y)j_d^\phi(|y|)\,\mathrm{d}y
    =
    0 .
$$
\end{assumption}

Additional assumptions on the spatial regularity and nondegeneracy of $a$ will be specified in the statements of the corresponding theorems.
For probabilistic results, the coefficient will be assumed to be real-valued and uniformly nondegenerate.

\subsection{Weighted boundedness of nonlocal operators}

We first state the weighted $L_p$-boundedness of the nonlocal operator $\mathcal L_V^a$. 
The result shows that $\mathcal L_V^a$ has the same differentiability order as the operator $\phi(\Delta)$.
We impose the following continuity assumption on the coefficient in the spatial variable.

\begin{assumption}
\label{24.08.18.13.17}
Let $a:\mathbb R^d\times\mathbb R^d\to\mathbb C$ be a bounded measurable function. 
We assume that there exists an increasing function
$\nu:[0,\infty)\to[0,\infty)$ such that
$$
\sup_{y\in\mathbb R^d}
|a(x_1,y)-a(x_2,y)|
\leq
\nu(|x_1-x_2|),
\quad x_1,x_2\in\mathbb R^d,
$$
and
$$
\int_{|z|\leq1}
\frac{\nu(|z|)}{|z|^{d+\beta}}\,\mathrm dz
<\infty
$$
for some $\beta\in(0,1)$.
\end{assumption}

\begin{rem}\label{rmk 07.16.16:08}
 Note that since $\nu$ is increasing, we get
 \begin{align*}
 \int_{r<|y|<2r} \frac{\nu(|y|)}{|y|^{d+\beta}} \mathrm{d}y \geq \nu(r)\int_{r<|y|<2r} \frac{1}{|y|^{d+\beta}} \mathrm{d}y \geq Cr^{-\beta}\nu(r) \quad \forall \,  0<r<1/2.
 \end{align*}
 This and Assumption \ref{24.08.18.13.17} imply that $a(x,y)$ is $\beta$-H\"older continuous in $x$ uniformly in $y$.
\end{rem}

\begin{thm}
\label{thm 11.10.18:01}
Let $\phi$ be a Bernstein function satisfying Assumption \ref{asm 07.21.13:06}$(\delta_0,\delta_1)$ with $0<\delta_0\leq\delta_1<1$. 
Let $V$ be a nonnegative function satisfying \eqref{eqn 02.11.14:19}, and let
$a:\mathbb R^d\times\mathbb R^d\to\mathbb C$ satisfy
Assumption \ref{24.08.18.13.17}. 
If $\delta_0\leq 1/2\leq\delta_1$, we additionally assume Assumption \ref{24.07.02.14.41}.
Then, for every $p\in(1,\infty)$ and $w\in A_p(\mathbb R^d)$,
$$
\|\mathcal L_V^a u\|_{L_p(w)}
\leq
C
\|u\|_{H_p^{\phi,2}(w)}.
$$
Here, the constant $C$ depends only on
\begin{align*}
  \|a\|_{L_\infty(\mathbb R^d\times\mathbb R^d)},
c_0,c_1,C_V,d,\delta_0,\delta_1,p,
[w]_{A_p(\mathbb R^d)},
\beta,
\nu.
\end{align*}

\end{thm}
Theorem \ref{thm 11.10.18:01} identifies the natural order of $\mathcal L_V^a$ in weighted $L_p$ spaces. 
In the model case $\phi(r)=r^{\alpha/2}$, the operator $\phi(\Delta)$ is the fractional Laplacian $\Delta^{\alpha/2}$, and the estimate above reduces to the natural fractional-order bound. 
A naive computation would only give a second-order estimate of the form
$$
\|\mathcal L_V^a u\|_{L_p}
\le
C\bigl(\|u\|_{L_p}+\|D^2u\|_{L_p}\bigr),
$$
which is far from optimal when the order of the jump operator is strictly less than two. 
Thus, the point of Theorem \ref{thm 11.10.18:01} is that regularity is measured in the scale generated by $\phi(\Delta)$, rather than by the local second-order operator $D^2$.

Boundedness estimates for nonlocal operators have been studied from several related viewpoints. 
When the coefficient is independent of the spatial variable, that is, when $a=a(y)$, Ba\~nuelos and Bogdan \cite{BB07} obtained $L_p$ bounds for Fourier multipliers arising from modulated jumps of L\'evy processes.
This direction was further developed for nonsymmetric L\'evy processes by Ba\~nuelos, Bielaszewski, and Bogdan
\cite{BBB2011}.

For fractional-type kernels, Dong and Kim \cite{DK12} proved $L_p$-boundedness for nonlocal operators with kernels of the form
$$
K(y)=\frac{a(y)}{|y|^{d+\alpha}},
\qquad 0<\alpha<2,
$$
under rough assumptions on the coefficient $a(y)$; in particular, $a$ is not required to be homogeneous, smooth, or symmetric.
Mikulevi\v{c}ius and Pragarauskas \cite{MP12singular,MP14} obtained $L_p$ estimates for singular integral operators and for fractional integro-differential operators related to jump processes. 
In particular, for bounded measurable $m=m(y)$, with the usual cancellation condition in the critical case $\alpha=1$, they proved estimates of the form
$$
\|L_m^\alpha u\|_{L_p}
\le
C\|(-\Delta)^{\alpha/2}u\|_{L_p},
\qquad 1<p<\infty,
$$
where
$$
L_m^\alpha u(x)
=
\int_{\mathbb R^d}
\left(
u(x+y)-u(x)-\chi_\alpha(y)\nabla u(x)\cdot y
\right)
m(y)\frac{\mathrm{d}y}{|y|^{d+\alpha}} .
$$
They also treated spatially dependent kernels $m(x,y)|y|^{-d-\alpha}$ under additional regularity assumptions on $m$ in the spatial variable and, in that setting, $p$ should be sufficiently large.

Related $L_p$ estimates for nonlocal operators with nonsymmetric measurable kernels were obtained by Kim and Kim \cite{KK16nonlocal}.
Their results cover kernels of the form
$$
K(y)=a(y)J(y),
$$
where $J(y)\,dy$ is a rotationally invariant L\'evy measure and $a(y)$ is merely measurable.

The weighted boundedness problem for spatially dependent fractional kernels was treated by Dong, Jung, and Kim \cite{DJK23nonlocal}. 
They considered operators of the form
$$
\mathcal L^a u(x)
=
\int_{\mathbb R^d}
\left(
u(x+y)-u(x)-\chi_\alpha(y)\nabla u(x)\cdot y
\right)
a(x,y)\frac{\mathrm{d}y}{|y|^{d+\alpha}},
$$
where $a(x,y)$ is bounded, merely measurable in $y$, and H\"older continuous in $x$. 
They proved
$$
\mathcal L^a:
H_p^\alpha(\mathbb R^d,w(x)\,\mathrm{d}x)
\to
L_p(\mathbb R^d,w(x)\,\mathrm{d}x)
$$
for all $p\in(1,\infty)$ and all Muckenhoupt weights $w\in A_p(\mathbb R^d)$. Note that they removed the condition that $p$ is sufficiently large by using the theory of weighted measure spaces.

For nonlocal operators beyond the exact fractional order, Kim, Kim, and Kim \cite{KKK13,KKK19} developed $L_p$ theories for $\phi(\Delta)$-type operators and for equations associated with processes having nonstationary independent increments. 
Mikulevi\v{c}ius and Phonsom \cite{mikulevivcius2017p,mikulevivcius2019cauchy} studied equations with scalable operators and in spaces of generalized smoothness. 
Recently, Janreung, Siripraparat, and Saksurakan \cite{JSS25} considered parabolic integro-differential equations with spatially dependent coefficients in generalized Bessel potential spaces, where the smoothness is defined by L\'evy measures with O-regularly varying profiles.

The main difficulty is the lack of exact homogeneity. 
In the fractional case, many estimates can be reduced to the scaling of the fractional Laplacian.
For a general Bernstein function satisfying only weak scaling, such homogeneous arguments are no longer available. 
To overcome this, we prove a Haj\l asz-type pointwise estimate adapted to the $\phi(\Delta)$: for suitable $\gamma\in(0,1]$,
$$
|u(x)-u(y)|
\le
C
\bigl(\phi(|x-y|^{-2})\bigr)^{-\gamma}
\left(
\mathbb M[\phi^\gamma(\Delta)u](x)
+
\mathbb M[\phi^\gamma(\Delta)u](y)
\right).
$$
This estimate is inspired by the classical pointwise approach, but the modulus is governed by the Bernstein function $\phi$. 

We also obtain sharper bounds in the spatially homogeneous case $a=a(y)$.
In this case, the lower-order term $\|u\|_{L_p(w)}$ can be removed from the right-hand side; see Theorem \ref{thm 01.15.15:07} (iii). 
Moreover, away from the critical scaling range, the cancellation assumption can be relaxed; see Theorem \ref{thm 01.15.15:07} (i)--(ii).

\subsection{Weighted Sobolev regularity theory for nonlocal PDEs}

We next apply the weighted boundedness theorem to develop a weighted Sobolev regularity theory for parabolic and elliptic equations driven by $\mathcal L_V^a$. 
In this subsection, the coefficient is allowed to depend on time.

\begin{assumption}
\label{asm 11.10.16:43}
Let $a=a(t,x,y):\mathbb R\times\mathbb R^d\times\mathbb R^d\to\mathbb R$ be a bounded measurable function. 
We assume the following conditions.

\begin{enumerate}[(i)]
\item
For each $t\in\mathbb R$, the function $a(t,\cdot,\cdot)$ satisfies Assumption \ref{24.08.18.13.17}. 
If $\delta_0\leq 1/2\leq\delta_1$, then $a(t,\cdot,\cdot)$ also satisfies Assumption \ref{24.07.02.14.41}. 
The constants and the modulus of continuity appearing in these assumptions are independent of $t$.

\item
There exist constants $0<a_0\leq a_1<\infty$ such that
$$
a_{0} \leq a(t,x,y) \leq a_{1}
$$
for all $(t,x,y)\in\mathbb R\times\mathbb R^d\times\mathbb R^d$.
\end{enumerate}
\end{assumption}

\begin{defn}
Let $1<p,q<\infty$, $\gamma \in \mathbb{R}$, $w_1\in A_p(\mathbb R^d)$,
$w_2\in A_q(\mathbb R)$, and $0<T<\infty$.

\begin{enumerate}[(i)]
\item
We denote by $H^{\phi,\gamma}_{q,p}(w_2,w_1,T)$ the space of all measurable functions $u:(0,T)\to H_p^{\phi,\gamma}(w_1)$ such that
$$
\|u\|_{H^{\phi,\gamma}_{q,p}(w_2,w_1,T)}^q
:=
\int_0^T
\|u(t,\cdot)\|_{H_p^{\phi,2}(w_1)}^q
w_2(t)\,\mathrm dt
<\infty.
$$
Similarly, $L_{q,p}(w_2,w_1,T)$ denotes the space of all measurable functions $u:(0,T)\to L_p(w_{1})$ satisfying
$$
\|u\|_{L_{q,p}(w_2,w_1,T)}^q
:=
\int_0^T
\|u(t,\cdot)\|_{L_p(w_1)}^q
w_2(t)\,\mathrm dt
<\infty.
$$

\item
We write $u\in C_p^{1,\infty}([0,T]\times\mathbb R^d)$ if, for every multi-index $\alpha$, both $D_x^\alpha u$ and $\partial_tD_x^\alpha u$ belong to $C([0,T];L_p)$.

\item
We say that $u\in\mathbb H^{\phi,2}_{q,p,0}(w_2,w_1,T)$ if there exists a sequence $u_n\in C_p^{1,\infty}([0,T]\times\mathbb R^d)$ with $u_n(0,\cdot)=0$ such that $u_n$ converges to $u$ in $H^{\phi,2}_{q,p}(w_2,w_1,T)$ and $\partial_tu_n$ is a Cauchy sequence in $L_{q,p}(w_2,w_1,T)$. 
In this case, $\partial_tu$ is defined as the limit of $\partial_tu_n$ in $L_{q,p}(w_2,w_1,T)$.
The norm on $\mathbb H^{\phi,2}_{q,p,0}(w_2,w_1,T)$ is defined by
$$
\|u\|_{\mathbb H^{\phi,2}_{q,p,0}(w_2,w_1,T)}
:=
\|\partial_tu\|_{L_{q,p}(w_2,w_1,T)}
+
\|u\|_{H^{\phi,2}_{q,p}(w_2,w_1,T)}.
$$
When $w_1\equiv1$ or $w_2\equiv1$, we omit the corresponding weight from the notation.
\end{enumerate}
\end{defn}

\begin{rem}
Let $u_{n}$ be a defining sequence of $u \in \mathbb{H}^{\phi,2}_{q,p,0}(w_{2},w_{1},T)$. Then for any for any $\eta \in C^{\infty}_{c}([0,T)\times \mathbb{R}^{d}$), one can simply check that
\begin{align}\label{eqn 12.24.15:39}
\int_{0}^{T} \int_{\mathbb{R}^{d}} \partial_{t} u_{n}(t,x) \eta(t,x) \mathrm{d}x \mathrm{d}t = -\int_{0}^{T} \int_{\mathbb{R}^{d}} u_{n}(t,x) \partial_{t}\eta(t,x) \mathrm{d}x \mathrm{d}t.
\end{align}
By H\"older's inequality, we obtain the following estimates
\begin{align}\label{eqn 12.24.15:37}
&\int_{0}^{T} \int_{\mathbb{R}^{d}} \left| \partial_{t}u_{n}(t,x) - \partial_{t}u(t,x) \right| \left| \eta(t,x) \right| \mathrm{d}x \mathrm{d}t   \nonumber
\\
& \leq \int_{0}^{T}  \|\partial_{t}u_{n}(t,\cdot) - \partial_{t}u(t,\cdot)\|_{L_{p}(w_{1})} \, \|\eta(t,\cdot)\|_{L_{p'}(\bar{w}_{1})} \mathrm{d}t   \nonumber
\\
&\leq \|\partial_{t}u_{n} - \partial_{t}u\|_{L_{q,p}(w_{2},w_{1},T)} \, \|\eta\|_{L_{q',p'}(\bar{w}_{2},\bar{w}_{1},T)},
\end{align}
where $p' = p/(p-1)$, $q'=q/(q-1)$, $\bar{w}_{1}=w_{1}^{-1/(p-1)} \in A_{p'}(\mathbb{R}^{d})$, and $\bar{w}_{2} = w_{2}^{-1/(q-1)} \in A_{q'}(\mathbb{R})$. Similarly, we also get
\begin{equation}
\label{eqn 12.24.15:37-2}
\begin{aligned}
&\int_{0}^{T} \int_{\mathbb{R}^{d}} \left|  u_{n}(t,x) - u(t,x) \right|  \left| \partial_{t}\eta(t,x)  \right| \mathrm{d}x \mathrm{d}t
 \\
 \leq& \|u_{n} - u\|_{L_{q,p}(w_{2},w_{1},T)} \, \|\partial_{t}\eta\|_{L_{q',p'}(\bar{w}_{2},\bar{w}_{1},T)}.
\end{aligned}
\end{equation}
The estimates \eqref{eqn 12.24.15:37}, and \eqref{eqn 12.24.15:37-2} allow us to take $n \to \infty$ in \eqref{eqn 12.24.15:39} to deduce
\begin{align*}
\int_{0}^{T} \int_{\mathbb{R}^{d}} \partial_{t} u(t,x) \eta(t,x) \mathrm{d}x \mathrm{d}t = -\int_{0}^{T} \int_{\mathbb{R}^{d}} u(t,x) \partial_{t}\eta(t,x) \mathrm{d}x \mathrm{d}t.
\end{align*}
This implies that the (time) trace of $u \in \mathbb{H}^{\phi,2}_{q,p,0}(w_{2},w_{1},T)$ is zero. Since our main object is an equation with zero initial condition, it is enough to consider only the space $\mathbb{H}^{\phi,2}_{q,p,0}(w_{2},w_{1},T)$.


\end{rem}

\begin{rem}
(i) One can easily check that $\mathbb{H}^{\phi,2}_{q,p,0}(w_{2},w_{1},T)$ is a Banach space due to its definition. 

(ii) Let $f \in H^{\phi,2}_{p}(w_{1})$, and let $\eta \in C^{\infty}_{c}(B_{1})$. For each $\varepsilon >0$, define
\begin{align*}
f^{\varepsilon}(x) = \int_{\mathbb{R}^{d}} \varepsilon^{-d} \eta(y/\varepsilon) f(x-y) \mathrm{d}y.
\end{align*}
Then one can check that (see \textit{e.g.} \cite[Remark 2.4 (ii)]{HKP20fractional}) as $\varepsilon \downarrow 0$, $f^{\varepsilon}$ converges to $f$ in $H^{\phi,2}_{p}(w_{1})$. Therefore, by considering a multiplication with a smooth cut-off function, one can also check that $C^{\infty}_{c}(\mathbb{R}^{d})$ is dense in $H^{\phi,2}_{p}(w_{1})$.

(iii) Let $u \in \mathbb{H}^{\phi,2}_{q,p,0}(w_{2},w_{1},T)$ and take $\eta_{1} \in C^{\infty}_{c}((1,2))$, $\eta_{2} \in C^{\infty}_{c}(B_{1})$, and $\eta_{3} \in C^{\infty}_{c}(B_{2})$ so that $\eta_{3}(x)=1$ for $|x|\leq 1$. For $\varepsilon_{1},\varepsilon_{2},\varepsilon_{3} >0$, define
\begin{align*}
u^{\varepsilon_{1},\varepsilon_{2},\varepsilon_{3}}(t,x) = \eta_{3}(\varepsilon_{3}x) \int_{0}^{\infty} \int_{\mathbb{R}^{d}} \varepsilon_{1}^{-d-1}\eta_{1}(s/\varepsilon_{1})\eta_{2}(y/\varepsilon_{1}) \bar{u}(t-s,x-y) \mathrm{d}y \mathrm{d}s,
\end{align*} 
where $\bar{u}(t,x) = u(t,x) \mathbf{1}_{t\in [0,T]}$ is an extension of $u$. 
Then one can check that $u^{\varepsilon_{1},\varepsilon_{2},\varepsilon_{3}} \in C^{\infty}_{c}(\mathbb{R}^{d+1}_{+})$ and converges to $u$ in $\mathbb{H}^{\phi,2}_{q,p,0}(w_{2},w_{1},T)$ as $\varepsilon_{1},\varepsilon_{2},\varepsilon_{3} \downarrow 0$. 
This certainly implies that $C^{\infty}_{c}(\mathbb{R}^{d+1}_{+})$ is dense in $\mathbb{H}^{\phi,2}_{q,p,0}(w_{2},w_{1},T)$.
\end{rem}

\begin{thm}
\label{thm 11.24.15:58}
Let $\phi$ be a Bernstein function satisfying Assumption \ref{asm 07.21.13:06}$(\delta_0,\delta_1)$ with $0<\delta_0\leq\delta_1<1$. Let $1<p,q<\infty$, $w_1\in A_p(\mathbb R^d)$, $w_2\in A_q(\mathbb R)$, and $0<T<\infty$.
Assume that
$$
V(y)=j_d^\phi(|y|),
$$
and that $a=a(t,x,y)$ satisfies Assumption \ref{asm 11.10.16:43}.
Then, for every $f\in L_{q,p}(w_2,w_1,T)$, there exists a unique solution $u\in\mathbb H^{\phi,2}_{q,p,0}(w_2,w_1,T)$ to
\begin{align}\label{eqn 11.24.16:10}
\partial_tu
=
\mathcal L_V^a u+f,
\qquad
u(0,\cdot)=0.
\end{align}
Moreover, the solution satisfies
\begin{align}\label{eqn 12.02.15:54}
\|u\|_{\mathbb H^{\phi,2}_{q,p,0}(w_2,w_1,T)}
\leq
C
\|f\|_{L_{q,p}(w_2,w_1,T)}.
\end{align}
Here the constant $C$ is independent of $f$ and depends only on 
\begin{align*}
a_{0},a_{1},
c_0,c_1,C_V,d,\delta_0,\delta_1,p,q,
[w_{1}]_{A_p(\mathbb R^d)},[w_{2}]_{A_{q}(\mathbb{R})}
\beta,
\nu,T.
\end{align*}
\end{thm}
We also prove the elliptic counterpart in
Theorem \ref{thm 01.26.16:45}.

\begin{rem}
\label{26.08.06.14.07}
Suppose $f\in L_{q,p}(w_2,w_1,T)\cap L_{q^\ast,p^\ast}(w_2^\ast,w_1^\ast,T)$, and let  $u\in\mathbb H^{\phi,2}_{q,p,0}(w_2,w_1,T)$ and $u^\ast\in\mathbb{H}^{\phi,2}_{q^\ast,p^\ast,0}(w_2^\ast,w_1^\ast,T)$ be corresponding solutions. In this remark, we show that $u=u^{\ast}$ in the sense of distributions.

By exploiting Remark \ref{rmk 02.02.23:45}, it is possible to choose $1<\rho<r$ close enough to $1$ so that
\begin{gather*}
w_{1} \in A_{p/r}(\mathbb{R}^{d}), \quad w_{2} \in A_{q/r}(\mathbb{R}), \quad w^{\ast}_{1} \in A_{p^{\ast}/r}(\mathbb{R}^{d}), \quad w^{\ast}_{2} \in A_{q^{\ast}/r}(\mathbb{R}),
\\
w_{1}^{-\frac{r}{p-r}} \in A_{\frac{p}{p-\rho}}(\mathbb{R}^{d}), \quad  w_{2}^{-\frac{r}{q-r}} \in A_{\frac{q}{q-\rho}}(\mathbb{R}^{d}), \quad {w^{\ast}_{1}}^{-\frac{r}{p^{\ast}-r}} \in A_{\frac{r}{p^{\ast}-\rho}}(\mathbb{R}^{d}), \quad {w^{\ast}_{2}}^{-\frac{r}{q^{\ast}-r}} \in A_{\frac{q^{\ast}}{q^{\ast}-s}}(\mathbb{R}^{d}).
\end{gather*}
Therefore, using H\"older's inequality and \eqref{eqn 01.30.01:21}, for any $h \in L_{p}(w_{1})$, it follows that
\begin{equation}\label{eqn 11.03.15:25}
\begin{aligned}
\int_{B_R(0)}|h(x)|^r\,\mathrm{d}x 
&\leq \left( \int_{B_{R}(0)} |h(x)|^{p}w_{1}(x)\mathrm{d}x \right)^{r/p} \left( \int_{B_{R}(0)} \left( w_{1}(x) \right)^{-\frac{r}{p-r}} \mathrm{d}x \right)^{(p-r)/p}
\\
& \leq \|h\|^{r}_{L_{p}(w_{1})} C R^{d\left( \frac{p}{p-\rho}\frac{p-r}{p} \right)} = C R^{d-d\frac{r-\rho}{p-\rho}} \|h\|^{r}_{L_{p}(w_{1})} \quad \forall\, R\geq1.
\end{aligned}
\end{equation}
Choose $0<\delta<d\frac{r-\rho}{p-\rho}$ and set $W(x)=(1+|x|)^{-d+\delta} \in A_{1}(\mathbb{R}^{d}) \subset A_{r}(\mathbb{R}^{d})$.
Since
\begin{align*}
W(x)\le 2^{-k(d-\delta)} \quad \forall\, 2^k\le |x|<2^{k+1} \quad (k=1,2,\cdots),
\end{align*}
for any $h \in L_{p}(w_{1})$, we have
\begin{align*}
\|h\|_{L_r(W)}^r
&\le \int_{B_1(0)} |h(x)|^r\,\mathrm{d}x
  +\sum_{k=0}^\infty 2^{-k(d-\delta)}
    \int_{B_{2^{k+1}}(0)} |h(x)|^r\,\mathrm{d}x\\
&\le C\left(1+\sum_{k=0}^\infty 2^{-k(d\frac{r-\rho}{p-\rho}-\delta)}\right)
       \|h\|_{L_p(w_1)}^r
 \le C\|h\|_{L_p(w_1)}^r.
\end{align*}
Applying the same argument as in \eqref{eqn 11.03.15:25} for $w_{2}$, for any $F \in L_{q,p}(w_{2},w_{1},T)$,
\begin{align*}
\|F\|_{L_{r,r}(1,W,T)}
&\le C\left(\int_0^T \|F(s,\cdot)\|_{L_p(w_1)}^r\,\mathrm{d}s\right)^{1/r}\\
&\le C\left(\int_0^T w_2(s)^{-r/(q-r)}\,\mathrm{d}s\right)^{1/r-1/q}
       \|F\|_{L_{q,p}(w_2,w_1,T)}.
\end{align*} 
The same argument yields that $L_{q^{\ast},p^{\ast}}(w^{\ast}_{2},w^{\ast}_{1},T) \subset L_{r,r}(1,W,T)$.
Therefore, under the assumption for $f$, $u$, and $u^{\ast}$, it follows that $u-u^{\ast} \in \mathbb{H}^{\phi,2}_{r,r,0}(1,W,T)$ is a solution to $\partial_{t}w = \mathcal{L}^{a}_{V}w$ with zero initial condition.
Hence, \eqref{eqn 12.02.15:54} gives $u=u^\ast$ in $\mathbb{H}^{\phi,2}_{r,r,0}(1,W,T)$, and hence in distribution.
\end{rem}

Theorem \ref{thm 11.24.15:58} gives a regularity estimate in weighted mixed-norm Sobolev space for parabolic equations driven by nonlocal operators. The Sobolev and mixed-norm regularity theory for nonlocal equations has a long history. 
Since the literature is extensive, we discuss only the results most closely related to our setting. 
For fractional kernels $|y|^{-d-\alpha}$, Mikulevi\v{c}ius and Pragarauskas \cite{MP92} developed an early $L_p$ theory for Cauchy problems in Sobolev and H\"older spaces. 
Later, Mikulevi\v{c}ius and Pragarauskas \cite{MP14} studied parabolic integro-differential equations whose principal kernels have the form $m(t,x,y)|y|^{-d-\alpha}$, with $m$ measurable in the jump variable and H\"older continuous in the spatial variable; their results are also connected with the martingale problem.
Zhang \cite{Z13-2} obtained $L_p$-maximal regularity for nonlocal parabolic equations with singular and nonsymmetric L\'evy operators, and also proved solvability results for equations with spatially dependent nonsmooth kernels \cite{Z13}.

Several subsequent works developed the theory in different directions.
Dong and Kim \cite{DK12} proved $L_p$ estimates for elliptic nonlocal equations with rough fractional kernels.
Dong, Jung, and Kim \cite{DJK23nonlocal} established the boundedness of spatially dependent fractional operators and the corresponding elliptic and parabolic $L_p$ estimates in weighted spaces. 
Dong and Liu \cite{DL23} studied fractional parabolic equations with space-time nonlocal operators.
More recently, Dong and Ryu \cite{DR26weight} established weighted mixed-norm estimates for equations with a Caputo time derivative and spatially translation-invariant nonlocal operators of fixed order $\sigma\in(0,2)$, whose kernels are merely measurable in time.
Their estimates allow Muckenhoupt weights in both space and time.
We also refer to Choi, Kim, and Ryu \cite{CKR23} and Dong and Ryu \cite{DR24} for nonlocal equations on domains, to Dong and Ryu \cite{DR26} for equations with general L\'evy measures.

For kernels beyond the exact fractional scaling, Kim, Kim, and Kim \cite{KKK13,KKK19} developed $L_p$ theories for equations associated with $\phi(-\Delta)$-type operators and with processes having non-stationary independent increments. 
Mikulevi\v{c}ius and Phonsom \cite{mikulevivcius2017p,mikulevivcius2019cauchy} studied parabolic and elliptic equations with scalable operators and equations in spaces of generalized smoothness. 
Related developments include the $L_q(L_p)$ theory for parabolic pseudo-differential equations via Calder\'on--Zygmund methods \cite{KKL16}, the $L_q(L_p)$ theory for equations with space-time nonlocal operators \cite{KPR21nonlocal}, the anisotropic nonlocal theory \cite{CKP24anisotropic}, the generalized Bessel potential framework with O-regularly varying L\'evy profiles \cite{JSS25}, and the low-intensity kernel theory \cite{KP25}.
Yang and Zhou \cite{YZ25} established weighted $L_q(L_p)$ estimates for time-fractional equations driven by the time-independent operator $\phi(\Delta)$, where $\phi$ is a Bernstein function satisfying a lower scaling condition.

The contribution of Theorem \ref{thm 11.24.15:58} is to combine, in a single weighted mixed-norm setting, three features which have not previously been treated simultaneously: jumping kernels governed by a Bernstein function without exact order of scaling, time-dependent coefficients with a spatial variable, and weighted measure spaces (with Muckenhoupt weights) in both space and time. 
In contrast to the exact fractional case, the order of differentiability is determined by $\phi$ and need not correspond to a fixed fractional power.
The use of Muckenhoupt weights is compatible with the weighted Hardy--Littlewood maximal estimates, Rubio de Francia extrapolation, and weighted Calder\'on--Zygmund estimates used in the proof.

\subsection{Weighted Sobolev--H\"older embeddings}

We next establish a pointwise embedding theorem for the weighted Sobolev spaces generated by $\phi(\Delta)$. 
Here, the term Sobolev--H\"older is used in the sense of an embedding from a weighted Sobolev spaces into a space of locally continuous functions with a specified modulus of continuity.

Classical Sobolev--H\"older embeddings show that sufficiently high Sobolev regularity implies pointwise H\"older regularity. 
A closely related pointwise approach is given by the inequalities of Bojarski and Haj{\l}asz \cite{BH93} and to the Haj{\l}asz--Sobolev theory on metric measure spaces \cite{H96}.
Weighted Sobolev embeddings with Muckenhoupt weights have also been studied in various forms; see, for instance, \cite{GU09,GW91} and the function-space embedding results of Meyries and Veraar \cite{MeyriesVeraar14}.

The embedding proved in this paper is different from the above classical weighted embedding results. 
The admissible modulus of continuity depends on the Bernstein function $\phi$ and, the Muckenhoupt weights. 
The estimate also keeps track of the local weighted mass of balls, which will be used for the later Krylov-type occupation bounds.

The precise statement is the following weighted Sobolev--H\"older embedding.

\begin{thm}[Weighted Sobolev--H\"older embedding]
\label{thm 05.19.17:17}
Let $\phi$ be a Bernstein function satisfying Assumption \ref{asm 07.21.13:06}$(\delta_0,\delta_1)$ with $0<\delta_0\le\delta_1<1$.
Let $\gamma\in(0,1]$ satisfy $2\delta_1\gamma<1$.
Suppose that $w_1\in A_\infty(\mathbb R^d)$ and $w_2\in A_\infty(\mathbb R)$.
Let $p,q\in(1,\infty)$, $p_0\in(R_{w_1},p)$, and $q_0\in(R_{w_2},q)$ be such that $w_1\in A_p(\mathbb R^d)$, $w_2\in A_q(\mathbb R)$, and
\begin{equation}
\label{26.02.03.21.14}
2\delta_0\gamma
-
2\delta_0\gamma\frac{q_0}{q}
-
\frac{dp_0}{p}
>0.
\end{equation}
Then the space $\mathbb H^{\phi,\gamma}_{q,p}(w_2,w_1,T)$ embeds into $C^{\omega_{\phi,\gamma}}_{\rm loc}([0,T]\times\mathbb R^d)$, where
$$
\|u\|_{C^{\delta_{\phi,\gamma}}(Q)}:=\|u\|_{L_{\infty}(Q)}+[u]_{C^{\delta_{\phi,\gamma}}(Q)},\quad [u]_{C^{\delta_{\phi,\gamma}}(Q)}:\sup_{\substack{(t,x)\in Q \\ (s,y)\in Q}}\frac{|u(t,x)-u(s,y)|}{\delta_{\phi,\gamma}(|x-y|)+\delta_{\phi,\gamma}(d_{\phi,\gamma}(t,s))},
$$
and
$$
\delta_{\phi,\gamma}(\rho):=\phi(\rho^{-2})^{-\gamma\left(1-\frac{q_0}{q}\right)}\rho^{-\frac{dp_0}{p}},\quad d_{\phi,\gamma}(t,s):=(\phi^{-1}(|t-s|^{-1/\gamma}))^{-1/2}.
$$
In particular, if $w_{1} \equiv 1$, $w_{2} \equiv 1$, then $p_{0}=q_{0}=1$, and the embedding is global.
\end{thm}
The condition $p_0>R_w$ in Theorem \ref{thm 05.19.17:17} reflects the loss caused by the weights. 
This loss is not an artifact of the proof. 
Even for the elliptic case, the endpoint corresponding to the critical index $R_w$ cannot be reached in general; see Example \ref{rem:endpoint-Rw-counterexample}.

\subsection{Martingale problems and Krylov-type estimates}

We now turn to the probabilistic consequences of the weighted parabolic regularity theory. 
For time-dependent coefficients, the natural object is the time-inhomogeneous martingale problem associated with the family of operators $\{\mathcal L_{\phi,t}^a\}_{t\ge0}$, or equivalently with the space-time operator $\partial_t+\mathcal L_{\phi,t}^a$.
The weighted parabolic estimates give solvability and a priori estimates for the associated backward equations. 
The weighted Sobolev--H\"older embedding then converts these Sobolev estimates into pointwise estimates.
This pointwise control yields the Krylov-type estimates.

The Krylov-type estimates should be viewed as weighted occupation-time bounds. 
Precisely, through those estimates, we control quantities of the form
$$
    \mathbb E_{s,x}\int_s^T |f(t,\mathcal X_t)|\,\mathrm{d}t
$$
by weighted mixed norms of $f$, with local normalizations determined by the weighted mass of parabolic cylinders or spatial balls.

\begin{defn}[Time-dependent martingale problem]
\label{def:time-dependent-martingale-problem}
Let $0\leq s<T<\infty$, and let $\mathbb{D}_{s,T}:=\mathbb{D}([s,T];\mathbb{R}^d)$ be the Skorokhod space of $\mathbb{R}^{d}$-valued c\`adl\`ag functions, equipped with the Skorokhod $J_1$-topology.
We also define $\mathbb{D}_{s,\infty}:=\mathbb{D}([s,\infty);\mathbb{R}^{d})$ in the same way. Let $\mathcal{X}_t(\omega):=\omega(t)$ be the canonical process on $\mathbb{D}_{s,T}$.
Define the raw canonical filtration by
$$
\mathcal{D}_t^0:=
\sigma(\mathcal{X}_r:s\leq r\leq t),
\quad
t\in[s,T]  \quad \left( \mathcal{D}^{0}_{\infty} = \bigcup_{t \geq s} \mathcal{D}^{0}_{t} \right),
$$
and its right-continuous modification by
$$
\mathcal{D}_t
:=
\bigcap_{\varepsilon>0}
\mathcal{D}_{(t+\varepsilon)\wedge T}^0,
\quad
t\in[s,T].
$$
\begin{enumerate}[(i)]
\item Let $\{\mathcal{L}_{\phi,t}^a\}_{t\geq0}$ be the family of nonlocal operators under consideration. 
For $x\in\mathbb R^d$, a probability
measure $\mathbb{P}_{s,T,x}$ is called a solution to the time-dependent martingale problem for $(\mathcal{L}_{\phi,t}^a,s,x)$ on $[s,T]$ if $\mathbb{P}_{s,T,x}(\mathcal{X}_s=x)=1$ and, for every $F\in C^{1,2}_{b}([s,\infty)\times\mathbb{R}^{d})$, the process
\begin{equation*}
\label{eq:def-time-dependent-martingale}
\begin{aligned}
\mathfrak M_t^F
:=&
F(t,\mathcal X_t)-F(s,x)
-
\int_s^t
\left[
\partial_rF(r,\mathcal X_r)
+
\mathcal L_{\phi,r}^a
\bigl(F(r,\cdot)\bigr)(\mathcal X_r)
\right]
\,\mathrm{d}r,
\end{aligned}
\end{equation*}
is a martingale with respect to $(\mathcal{D}_t)_{t\in[s,T]}$ under $\mathbb{P}_{s,T,x}$.

\item We say that the time-dependent martingale problem for $(\mathcal{L}_{\phi,t}^a,s,x)$ is well posed on $[s,T]$ if it admits a unique solution $\mathbb P_{s,T,x}$.

\item We say that the time-dependent martingale problem for $\{\mathcal{L}_{\phi,t}^a\}_{t\geq0}$ is globally well posed if, for every $s\geq0$ and $x\in\mathbb R^d$, there exists a unique probability measure $\mathbb{P}_{s,x}
$
such that, for every $T>s$, its restriction
$$
\mathbb P_{s,x}
\circ
\left(
\mathcal X|_{[s,T]}
\right)^{-1}
$$
is a solution to the time-dependent martingale problem for
$(\mathcal L_{\phi,t}^a,s,x)$ on $[s,T]$.
\end{enumerate}
\end{defn}

Our final result establishes the global well-posedness of the
time-dependent martingale problem, the strong Markov property of
the associated time-inhomogeneous family, and a weighted
parabolic Krylov-type estimate.

\begin{thm}[Global well-posedness, strong Markov property, and a Krylov-type estimate]
\label{thm 02.26.19:05}

Let $\phi$ be a Bernstein function satisfying Assumption \ref{asm 07.21.13:06}$(\delta_0,\delta_1)$ with $0<\delta_0\leq\delta_1<1$. 
Assume that
$$
V(y)=j_d^\phi(|y|),
$$
and that $a=a(t,x,y)$ satisfies Assumption \ref{asm 11.10.16:43}.

\begin{enumerate}[(i)]
\item
The time-dependent martingale problem for
$\{\mathcal L_{\phi,t}^a\}_{t\geq0}$ is globally well posed in the
sense of
Definition \ref{def:time-dependent-martingale-problem}.

\item
Let $\mathbb P_{s,x}$ denote the unique martingale solution
starting from $(s,x)$. Then the family
$$
\left\{
\mathbb P_{s,x}
:
s\geq0,\ x\in\mathbb R^d
\right\}
$$
is a time-inhomogeneous strong Markov family. More precisely, if $\tau\geq s$ is a finite (a.s.) $\mathcal{D}_{t}$ stopping time,
then for every $\vartheta\geq0$ and every bounded Borel function
$g:\mathbb R^d\to\mathbb R$,
\begin{equation*}
\label{eq:main-strong-Markov}
\mathbb E_{s,x}
\left[
g(\mathcal X_{\tau+\vartheta})
\mid
\mathcal D_\tau
\right]
=
P_{\tau,\tau+\vartheta}g(\mathcal X_\tau)
\qquad
\mathbb P_{s,x}\text{-a.s.},
\end{equation*}
where
$$
P_{r,t}g(z)
:=
\mathbb E_{r,z}[g(\mathcal X_t)],
\qquad
0\leq r\leq t.
$$

\item
Suppose that $w_{1} \in A_{\infty}(\mathbb{R}^{d})$, $w_{2} \in A_{\infty}(\mathbb{R})$. Choose $1<p_{0}<p<\infty$, $1<q_{0}<q<\infty$, and $\gamma \in (0,1]$ such that
$$
w_1\in A_p(\mathbb R^d),
\quad
w_2\in A_q(\mathbb R),\quad R_{w_1}<p_0<p,
\quad
R_{w_2}<q_0<q,
$$
and
\begin{equation*}
\label{eq:main-Krylov-exponents}
2\delta_1\gamma<1,
\qquad
2\delta_0\gamma
\left(
1-\frac{q_0}{q}
\right)
>
\frac{dp_0}{p}.
\end{equation*}
For $0<s<T$, and $x\in\mathbb R^d$, set
$\rho_{s,T}:=\frac{T-s}{2}\wedge1$, 
and define
$$
\mathcal K_{s,T,x}(f)
:=
\mathbb E_{s,x}
\int_s^T
f(r,\mathcal X_r)\,\mathrm{d}r,
\qquad
f\in C_c^\infty((0,T)\times\mathbb R^d).
$$
Then $\mathcal K_{s,T,x}$ extends uniquely to a bounded linear
functional on
$L_{q,p}(w_2,w_1,T)$.
Moreover, for every
$f\in L_{q,p}(w_2,w_1,T)$,
\begin{equation}
\label{eqn 01.26.15:17}
\left|
\mathcal K_{s,T,x}(f)
\right|
\leq
C
\frac{
\|f\|_{L_{q,p}(w_2,w_1,T)}
}{
[w_2(\mathbb I_{\phi,\gamma}(s,\rho_{s,T}))]^{1/q}
[w_1(B_{\rho_{s,T}}(x))]^{1/p}
}.
\end{equation}
The constant $C>0$ is independent of $s$, $x$, and $f$, and 
 $\mathbb{I}_{\phi,\gamma}(t,\rho):=[t,t+\phi(\rho^{-2})^{-\gamma})$.
\end{enumerate}
\end{thm}

Martingale problems for nonlocal operators go back to the works of Komatsu \cite{K73} and Stroock \cite{S75}. 
For pure jump and stable-like processes, fundamental uniqueness results were obtained by Bass \cite{Bass88} and by Bass and Tang \cite{BassTang09}.
Analytic approaches connecting Cauchy problems with martingale problems were developed, among others, by Abels and Kassmann \cite{AK09} and by Mikulevi\v{c}ius and Pragarauskas \cite{MP14,MP14-2}. 
In this paper, well-posedness of the martingale problem is investigated by using the weighted Sobolev regularity theory.

The second probabilistic consequence is a weighted Krylov-type estimate. 
Through such estimates, we control the expected time spent by the process in a set, or more generally, the integral of an observable along the sample path, in terms of its norm.
Occupation time densities for stable-like and pure jump Markov processes were studied by Bass \cite{Bass88Occupation}. 
We also refer to Krylov \cite{K74,K87,Krylov80} and Bahlali, Eddahbi, Ouknine \cite{BEY17} for classical formulations of Krylov-type estimates, and to the recent work of Hu and Zhao \cite{HZ25} for related estimates in generalized Orlicz--Besov spaces.

Our contribution is to obtain well-posedness of the martingale problem and Krylov-type estimates. Especially, we use the weighted mixed norm for the control of expected time spend by processes.

\section{Proof of weighted boundedness of nonlocal operators: Theorem \ref{thm 11.10.18:01}}
\subsection{Preliminaries}

In this subsection, we collect several estimates used throughout the proof of Theorem \ref{thm 11.10.18:01}. 
We first record some elementary consequences of the weak scaling condition. 
Then we recall the Green function estimate for $\phi(\Delta)$ and a weighted $L_{p}$ boundedness of Calder\'on--Zygmund operators.

\begin{lem}
\label{lem:scaling-integrals}
Let $\phi$ satisfy Assumption \ref{asm 07.21.13:06}$(\delta_0,\delta_1)$ with $0<\delta_0\leq\delta_1<1$, and let $b>0$.
Then the following estimates hold.

\begin{equation}
\label{eqn 07.01.15:42}
\begin{aligned}
\int_{r}^{\infty} \lambda^{-a}\left( \phi(\lambda^{-2}) \right)^{b} \mathrm{d}\lambda 
 \leq 
\frac{c_{0}^{-b}}{a+2\delta_{0}b-1} \left( \phi(r^{-2})  \right)^{b} r^{1-a} \quad \forall \, a+2\delta_{0}b >1,
\end{aligned}
\end{equation}
\begin{equation}
\label{eqn 02.13.11:44}
\begin{aligned}
\int_{r}^{\infty} \lambda^{-a}\left( \phi(\lambda^{-2}) \right)^{-b} \mathrm{d}\lambda 
 \leq  \frac{c_{0}^{b}}{a-2\delta_{1}b-1} \left( \phi(r^{-2})  \right)^{-b} r^{1-a} \quad \forall\,  a-2\delta_{1}b >1,
\end{aligned}
\end{equation}
\begin{equation}
\label{eqn 07.01.16:07}
\begin{aligned}
\int_{0}^{r} \lambda^{a} \left( \phi(\lambda^{-2}) \right)^{b} \mathrm{d}\lambda   
 \leq  \frac{c_{1}^{b}}{a-2\delta_{1}b+1}\left( \phi(r^{-2}) \right)^{b}r^{a+1} \quad \forall \, a-2\delta_{1}b >-1,
\end{aligned}
\end{equation}
\begin{equation}
\label{eqn 07.01.16:30}
\begin{aligned}
\int_{0}^{r} \lambda^{a} \left( \phi(\lambda^{-2}) \right)^{-b} \mathrm{d}\lambda   
 \leq \frac{c_{0}^{-b}}{a+2\delta_{0}b+1}\left( \phi(r^{-2}) \right)^{-b}r^{a+1} \quad \forall \,  a+2\delta_{0}b >-1.
\end{aligned}
\end{equation}
\end{lem}

\begin{proof}
We only prove the estimates briefly, since they are direct consequences of Assumption \ref{asm 07.21.13:06}. 
If $\lambda\geq r$, then
$\lambda^{-2}\leq r^{-2}$, and hence
$$
\phi(\lambda^{-2})
\leq
C\left(\frac r\lambda\right)^{2\delta_0}\phi(r^{-2}),
\qquad
\phi(\lambda^{-2})
\geq
C\left(\frac r\lambda\right)^{2\delta_1}\phi(r^{-2}).
$$
These two inequalities directly imply \eqref{eqn 07.01.15:42} and \eqref{eqn 02.13.11:44}.
If $0<\lambda\leq r$,
then $\lambda^{-2}\geq r^{-2}$, and
$$
\phi(\lambda^{-2})
\leq
C\left(\frac r\lambda\right)^{2\delta_1}\phi(r^{-2}),
\qquad
\phi(\lambda^{-2})
\geq
C\left(\frac r\lambda\right)^{2\delta_0}\phi(r^{-2}).
$$
These give \eqref{eqn 07.01.16:07} and \eqref{eqn 07.01.16:30}.
\end{proof}

We next recall the Green function estimate for $\phi(\Delta)$.

\begin{thm}[Green function estimate; {\cite[Theorem 1]{Z09Riesz}}]\label{thm 07.21.11:32}
Suppose that $\phi$ is a Bernstein function satisfying Assumption \ref{asm 07.21.13:06}$(\delta_{0},\delta_{1})$ with $\delta_{1}<1/2$.
Then there exists a radial function $K_d^\phi(x)=k_d^\phi(|x|)$ such that
$$
\mathcal F_d[K_d^\phi](\xi)
=
\frac{1}{\phi(|\xi|^2)}
$$
in the sense of distributions. 
Moreover,
$$
C^{-1}
\frac{1}{r^d\phi(r^{-2})}
\leq
k_d^\phi(r)
\leq
C
\frac{1}{r^d\phi(r^{-2})},
\qquad r>0,
$$
and $k_d^\phi$ is differentiable with
$$
\frac{\mathrm d}{\mathrm dr}k_d^\phi(r)
=
-2\pi r\, k_{d+2}^\phi(r),
\qquad r>0.
$$
\end{thm}

For later use, we introduce the difference kernel
\begin{equation}
\label{24.07.04.17.46}
\mathfrak K_d^\phi(z,y)
:=
K_d^\phi(z)-K_d^\phi(z+y).
\end{equation}

\begin{prop}
\label{24.07.04.17.53}
Let $\phi$ be a Bernstein function satisfying
Assumption \ref{asm 07.21.13:06}$(\delta_0,\delta_1)$ with $0<\delta_0\leq\delta_1<1$.

\begin{enumerate}
\item[\emph{(i)}] 
Suppose that $\delta_{1} < 1/2$. Then, for every $y\in\mathbb R^d\setminus\{0\}$,
\begin{align}
\label{eqn 06.30.14:58}
\int_{\mathbb R^d}
|\mathfrak K_d^\phi(z,y)|\,\mathrm dz
\leq
\frac{C}{\phi(|y|^{-2})}.
\end{align}
Moreover,
$$
\int_{\mathbb R^d}
\mathfrak K_d^\phi(z,y)\,\mathrm dz
=
0.
$$

\item[\emph{(ii)}]
For $\varepsilon\in(0,1/2]$, set $\phi^\varepsilon(r):=(\phi(r))^\varepsilon$. 
Then, for every $z\in\mathbb R^d\setminus\{0\}$,
$$
\int_{\mathbb R^d}
|\mathfrak K_d^{\phi^\varepsilon}(z,y)|
j_d^\phi(|y|)\,\mathrm dy
\leq
C
\frac{(\phi(|z|^{-2}))^{1-\varepsilon}}{|z|^d}.
$$
\end{enumerate}
The constants $C$ depend only on $c_0,c_1,\delta_0,\delta_1,d$ and, in \emph{(ii)}, also on $\varepsilon$.
\end{prop}

\begin{proof}
We prove only the estimates needed later. 
For (i), split
$$
\int_{\mathbb R^d}
|\mathfrak K_d^\phi(z,y)|\,\mathrm dz
=
\int_{|z|\leq2|y|}\cdots\,\mathrm dz
+
\int_{|z|>2|y|}\cdots\,\mathrm dz.
$$
By Theorem \ref{thm 07.21.11:32},
\begin{align*}
\int_{|z|\leq2|y|}
|\mathfrak K_d^\phi(z,y)|\,\mathrm dz
&\leq \int_{ |z| \leq 2|y|} (K^{\phi}_{d}(z+y) + K^{\phi}_{d}(z))\mathrm{d}z \\
&\leq
C
\int_{|z|\leq3|y|}
\frac{1}{|z|^d\phi(|z|^{-2})}\,\mathrm dz
\leq
\frac{C}{\phi(|y|^{-2})},
\end{align*}
where we used \eqref{eqn 07.01.16:30}. For the second part, the
fundamental theorem of calculus and Theorem \ref{thm 07.21.11:32} give
$$
|\mathfrak K_d^\phi(z,y)|
\leq
C|y|
\int_0^1
\frac{1}{|z+ry|^{d+1}\phi(|z+ry|^{-2})}\,\mathrm dr.
$$
Since $|z|>2|y|$, we have $|z+ry|\simeq |z|$, and hence
$$
\int_{|z|>2|y|}
|\mathfrak K_d^\phi(z,y)|\,\mathrm dz
\leq
C|y|
\int_{2|y|}^\infty
\frac{1}{\rho^2\phi(\rho^{-2})}\,\mathrm d\rho
\leq
\frac{C}{\phi(|y|^{-2})},
$$
where the last inequality follows from \eqref{eqn 02.13.11:44}; here the
condition $\delta_1<1/2$ is used.

It remains to show the cancellation. 
Let $\eta\in C_c^\infty(\mathbb R^d)$ satisfy $\eta(0)=1$. 
Since $\mathfrak K_d^\phi(\cdot,y)\in L_1(\mathbb R^d)$, by the dominated convergence theorem,
$$
\int_{\mathbb R^d}\mathfrak K_d^\phi(z,y)\,\mathrm dz
=
\lim_{n\to\infty}
\int_{\mathbb R^d}
\mathfrak K_d^\phi(z,y)\eta(z/n)\,\mathrm dz.
$$
Using the Fourier transform of $K_d^\phi$,
$$
\int_{\mathbb R^d}
\mathfrak K_d^\phi(z,y)\eta(z/n)\,\mathrm dz
=
\int_{\mathbb R^d}
\left(\mathrm{e}^{i y\cdot\xi/n}-1\right)
\frac{\mathcal{F}[\eta](\xi)}{\phi(n^{-2}|\xi|^2)}\,\mathrm d\xi.
$$
By the weak scaling condition and the fact that $\delta_1<1/2$, the absolute value of the last integral is bounded by
$$
C n^{-1+2\delta_1}|y|
\int_{\mathbb R^d}
|\xi|
\frac{|\mathcal{F}[\eta](\xi)|}{\phi(|\xi|^2)}\,\mathrm d\xi,
$$
which tends to zero as $n\to\infty$. 
This proves the cancellation.

For (ii), split the integral into the two regions
$|y|\leq |z|/2$ and $|y|>|z|/2$. 
In the first region, the fundamental theorem of calculus gives
$$
|\mathfrak K_d^{\phi^\varepsilon}(z,y)|
\leq 
C|y|
\int_0^1
\frac{1}{|z+ry|^{d+1}
(\phi(|z+ry|^{-2}))^\varepsilon}\,\mathrm dr.
$$
Since $|z+ry|\simeq |z|$ and $j_d^\phi(r)\simeq r^{-d}\phi(r^{-2})$, we obtain by \eqref{eqn 07.01.16:07}
$$
\int_{|y|\leq |z|/2}
|\mathfrak K_d^{\phi^\varepsilon}(z,y)|
j_d^\phi(|y|)\,\mathrm dy
\leq
C
\frac{(\phi(|z|^{-2}))^{1-\varepsilon}}{|z|^d}.
$$
In the region $|y|>|z|/2$, Theorem \ref{thm 07.21.11:32} yields
$$
|\mathfrak K_d^{\phi^\varepsilon}(z,y)|
\leq
C
\left(
\frac{1}{|z+y|^d(\phi(|z+y|^{-2}))^\varepsilon}
+
\frac{1}{|z|^d(\phi(|z|^{-2}))^\varepsilon}
\right).
$$
Splitting further into $|z|/2<|y|\leq2|z|$ and $|y|>2|z|$, and using \eqref{eqn 07.01.15:42} and \eqref{eqn 07.01.16:30}, we obtain the same bound.
This proves (ii).
\end{proof}

Finally, we recall a weighted $L_{p}$ boundedness of Calder\'on--Zygmund operators used in the case when $a=a(y)$. We bring the statement from \cite[Lemma 4.2]{DJK23nonlocal}.

\begin{thm}[Weighted Calder\'on--Zygmund theorem]\label{CZ}
Let $T$ be a linear operator of the form
$$
Tf(x)
=
\int_{\mathbb R^d}
K(z)f(x-z)\,\mathrm dz.
$$
Assume that
$$
\|Tf\|_{L_2(\mathbb{R}^d)}
\leq
C_1\|f\|_{L_2(\mathbb{R}^d)},
\qquad f\in L_2(\mathbb{R}^d),
$$
and that there exist constants $C_2>0$ and $\gamma>0$ such that
\begin{align}
\label{eqn 02.19.21:49}
|K(x+z)-K(x)|
\leq
C_2
\frac{|z|^\gamma}{|x|^{d+\gamma}},
\qquad |x|\geq4|z|.
\end{align}
Then, for every $p\in(1,\infty)$ and $w\in A_p(\mathbb R^d)$,
$$
\|Tf\|_{L_p(\mathbb{R}^d,w\,\mathrm{d}x)}
\leq
C\|f\|_{L_p(\mathbb{R}^d,w\,\mathrm{d}x)},
$$
where $C$ depends only on
$C_1,C_2,d,\gamma,p$, and $[w]_{A_p(\mathbb R^d)}$.
\end{thm}

\subsubsection*{Organization of the proof of Theorem \ref{thm 11.10.18:01}}

We provide a structure of the proof of 
Theorem \ref{thm 11.10.18:01}. The proof consists of two parts:
first the spatially homogeneous case, and then the localization
argument for spatially dependent coefficients.

{\small
$$
\begin{aligned}
(A)+(B) 
&\Longrightarrow (C),\\
(C)+\text{Theorem \ref{CZ}}
&\Longrightarrow (D)
\Longrightarrow (E),\\
(F)
&\Longrightarrow (G)
\Longrightarrow (H),\\
(E)+\text{Assumption \ref{24.08.18.13.17}}
&\Longrightarrow (I),\\
(E)+(H)+(I)
&\Longrightarrow \text{Theorem \ref{thm 11.10.18:01}.}
\end{aligned}
$$
}
Here,
{\small$$
\begin{array}{lll}
(A) \, \text{Lemma \ref{lem:scaling-integrals}: basic estimates from weak scaling}, &&  (B) \, \text{Theorem \ref{thm 07.21.11:32}: Green function estimate},\\[1mm]
(C) \, \text{Proposition \ref{24.07.04.17.53}: estimates on the difference kernel}, && (D) \, \text{Lemma \ref{24.07.04.18.22}: low-order $x$-independent estimate},\\[1mm]
(E) \, \text{Theorem \ref{thm 01.15.15:07}: $x$-independent coefficient case},&& (F) \, \text{Proposition \ref{lem 10.30.10:35}: Haj{\l}asz-type pointwise estimate},\\[1mm]
(G) \, \text{Lemma \ref{lem 08.12.18:53}: small-jump estimate},&& (H) \, \text{Lemma \ref{lem 07.17.19:30}: cutoff commutator estimate},\\[1mm]
(I) \, \text{Lemma \ref{lem 07.17.19:49}: localized estimate}, &&\\[1mm]
\end{array}
$$}

\subsection{The $x$-independent coefficient case}

In this subsection, we consider the case $a=a(y)$.

\begin{thm}
\label{thm 01.15.15:07}
Let $\phi$ be a Bernstein function satisfying Assumption \ref{asm 07.21.13:06}$(\delta_0,\delta_1)$ with
$0<\delta_0\leq\delta_1<1$, and let $V$ be a nonnegative function satisfying \eqref{eqn 02.11.14:19}. Suppose that $a=a(y)$ is a bounded
measurable function. 
Assume that one of the following:
\begin{enumerate}[(i)]
    \item $0<\delta_0\leq\delta_1<1/2$,
    \item $1/2<\delta_0\leq\delta_1<1$,
    \item The coefficient $a:\bR^{d} \to \bR$ satisfies Assumption \ref{24.07.02.14.41}.
\end{enumerate}
Then, for every $p\in(1,\infty)$ and $w\in A_p(\mathbb R^d)$,
$$
\|\mathcal L_V^a u\|_{L_p(w)}
\leq
C
\|\phi(\Delta)u\|_{L_p(w)}.
$$
Here, $C$ depends only on
$\|a\|_{L_\infty},c_0,c_1,C_V,d,\delta_0,\delta_1,p$, and
$[w]_{A_p(\mathbb R^d)}$.
\end{thm}

\begin{rem}
\label{rem:x-independent-cancellation-necessity}
We explain why a cancellation condition (Assumption \ref{24.07.02.14.41}) is necessary to cover the full range $0<\delta_0\leq\delta_1<1$. 
The proof of Theorem \ref{thm 01.15.15:07} shows that the critical situation arises when the scaling exponent meet the threshold $1/2$. 
Therefore, it is enough to consider the model case
$$
d=1,
\qquad
\phi(\lambda)=\lambda^{1/2},
\qquad
\phi(\Delta)=\Delta^{1/2}.
$$
In this case, the corresponding jumping kernel satisfies
$$
j_1^\phi(y)\simeq |y|^{-2}.
$$

(i) First, consider the operator without the compensation term,
$$
\mathcal L_1^a u(x)
:=
\int_{\mathbb R}
\bigl(u(x+y)-u(x)\bigr)
\frac{a(y)}{|y|^2}\,\mathrm dy.
$$
Let $u\in C_c^\infty(\mathbb R)$. 
By Taylor's theorem,
$$
u(x+y)
=
u(x)+yu'(x)+O(|y|^2)
\qquad
\text{as } y\to0.
$$
Thus, it follows that
$$
\mathcal L_1^a u(x)
=
u'(x)
\int_{|y| \leq 1}
\frac{y\,a(y)}{|y|^2}\,\mathrm dy
+
\int_{|y| \leq 1}
\frac{O(|y|^2)a(y)}{|y|^2}\,\mathrm dy + \int_{|y| >1} \bigl( u(x+y) - u(x) \bigr) \frac{a(y)}{|y|^{2}} \, \mathrm{d}y.
$$
The second and third integrals are finite. 
However, the first integral may diverge. 
Indeed, let
$$
a(y):=\operatorname{sgn}(y)\mathbf 1_{(-1,1)}(y).
$$
Then $a$ is bounded, but it does not satisfy the cancellation condition, since for every $0<r<R<1$,
$$
\int_{-R}^{-r}
y\,a(y)\frac{\mathrm dy}{|y|^2}
+
\int_r^R
y\,a(y)\frac{\mathrm dy}{|y|^2}
=
2\int_r^R\frac{\mathrm dt}{t}
=
2\log(R/r)
\neq0.
$$
Consequently, if $u'(x)\neq0$, then for such $a$, it follows that
$$
\left|u'(x)
\int_{|y|<1}
\frac{y a(y)}{|y|^{2}}
 \mathrm dy\right|=
\infty.
$$
Thus $\mathcal L_1^a u(x)$ is not well-defined  even as in the singular integral. 
This shows that, in the critical case, the compensation term $y\cdot\nabla u(x) \chi_{\delta_{0},\delta_{1}}(y)$ cannot be omitted.

(ii) Even after adding the compensation term, a cancellation condition is still needed in order to obtain the natural operator bound. 
Consider the Fourier symbol
$$
\psi(\xi)
:=
\int_{\mathbb R}
\left(
\mathrm{e}^{i\xi y}-1-i\xi y
\right)
a(y)\frac{\mathrm dy}{|y|^2}.
$$
Its imaginary part is
$$
\Im[\psi(\xi)]
=
\int_{\mathbb R}
\bigl(\sin(\xi y)-\xi y\bigr)
a(y)\frac{\mathrm dy}{|y|^2}.
$$
Now take
$$
a(y):=\mathbf 1_{(0,1)}(y).
$$
Again, $a$ does not satisfy the cancellation condition. 
By the change of variables $\lambda \mapsto|\xi|y$, for $\xi\neq0$ we get
$$
\Im[\psi(\xi)]
=
\operatorname{sgn}(\xi)|\xi|
\int_0^{|\xi|}
\bigl(\sin\lambda-\lambda\bigr)
\frac{\mathrm d\lambda}{\lambda^2}.
$$
Hence, for $|\xi|>2$,
\begin{align}
\label{eqn 01.12.17:27}
|\Im[\psi(\xi)]|
&\geq
|\xi|
\int_2^{|\xi|}
\frac{\lambda-\sin\lambda}{\lambda^2}\,\mathrm d\lambda
\nonumber\\
&\geq
\frac{|\xi|}{2}
\int_2^{|\xi|}
\frac{\mathrm d\lambda}{\lambda}
=
\frac12|\xi|\log(|\xi|/2).
\end{align}
Therefore, the multiplier associated with
$$
\mathcal L_2^a u(x)
:=
\int_{\mathbb R}
\bigl(u(x+y)-u(x)-yu'(x)\bigr)
\frac{a(y)}{|y|^2}\,\mathrm dy
$$
grows faster than  $|\xi|$. 
In particular, $\mathcal L_2^a$ cannot be bounded from $W_2^1(\mathbb R)$ to $L_2(\mathbb R)$. 
Indeed, if such a bound were true, then for every $u\in\mathcal S(\mathbb R)$ with
$$
\operatorname{supp}\mathcal{F}_1[u]
\subset
\{\xi\in\mathbb R:|\xi|>2\},
$$
Plancherel's theorem and \eqref{eqn 01.12.17:27} would imply
$$
\int_{|\xi|>2}
|\xi|^2\log^2(|\xi|/2)|\mathcal{F}_1[u](\xi)|^2\,\mathrm d\xi
\leq
C
\int_{|\xi|>2}
|\xi|^2|\mathcal{F}_1[u](\xi)|^2\,\mathrm d\xi.
$$
Since $u$ is arbitrary, this would force $\log(|\xi|/2)$ to be bounded for large $|\xi|$, which is impossible.
Thus, the natural estimate fails without cancellation.
\end{rem}

We first prove the estimate in the range $0<\delta_0\leq\delta_1<1/2$. 
In this case, no compensation term appears in the definition of $\mathcal L_V^a$.

\begin{lem}
\label{24.07.04.18.22}
Let $\phi$ be a Bernstein function satisfying Assumption \ref{asm 07.21.13:06}$(\delta_0,\delta_1)$ with $0<\delta_0\leq\delta_1<1/2$. 
Let $V$ be a nonnegative function satisfying \eqref{eqn 02.11.14:19}, and let $a=a(y)$ be a bounded measurable function. 
Then, for every $p\in(1,\infty)$, $w\in A_p(\mathbb R^d)$, and $u\in\mathcal S(\mathbb R^d)$,
$$
\sup_{\varepsilon\in(0,1)}
\left\|
\int_{\varepsilon<|y|<\varepsilon^{-1}}
\bigl(u(\cdot+y)-u(\cdot)\bigr)a(y)V(y)\,\mathrm dy
\right\|_{L_p(w)}
\leq
C\|\phi(\Delta)u\|_{L_p(w)}.
$$
Here $C$ depends only on
$c_0,c_1,\delta_0,\delta_1,d,p,\|a\|_{L_\infty},C_V$, and
$[w]_{A_p(\mathbb R^d)}$.
\end{lem}

\begin{proof}
Since $\delta_{1}<1/2$, we apply Theorem \ref{thm 07.21.11:32} to rewrite
\begin{align*}
u(x) = \phi(\Delta)^{-1}\phi(\Delta)u(x) = - \int_{\bR^{d}} K^{\phi}_{d}(z) \phi(\Delta)u(x-z) \mathrm{d}z.
\end{align*}
Hence, for any $y\in \bR^{d}$,
\begin{equation}\label{eqn 01.16.12:30}
\begin{aligned}
u(x+y) - u(x) &= \int_{\bR^{d}} K^{\phi}_{d}(z) \left[\phi(\Delta)u(x-z) - \phi(\Delta)u(x+y-z)\right] \mathrm{d}z 
\\
&=\int_{\bR^{d}} \mathfrak{K}^{\phi}_{d}(z,y) \, \phi(\Delta)u(x-z) \mathrm{d}z,
\end{aligned}
\end{equation}
where $\mathfrak{K}_d^{\phi}$ is a function defined in \eqref{24.07.04.17.46}. Therefore, we have
\begin{align*}
&\int_{\varepsilon<|y|<\varepsilon^{-1}} \left( u(x+y) - u(x)  \right) a(y) V(y)\mathrm{d}y
\\
&= \int_{\varepsilon<|y|<\varepsilon^{-1}} \left( \int_{\bR^{d}} \mathfrak{K}^{\phi}_{d}(z,y) \phi(\Delta)u(x-z) \mathrm{d}z \right)  a(y) V(y)\mathrm{d}y.
\end{align*}
One can check that the function
$$
\mathfrak{K}^{\phi}_{d}(z,y)\phi(\Delta)u(x-z)a(y)V(y)
$$
is integrable in $(y,z)$ by exploiting Proposition \ref{24.07.04.17.53}. Hence, by Fubini's theorem, we have
\begin{equation}
\label{eqn 10.23.11:27}
\begin{aligned}
&\int_{\varepsilon<|y|<\varepsilon^{-1}} \left( u(x+y) - u(x)  \right) a(y) V(y)\mathrm{d}y 
\\
& = \int_{\bR^{d}} \left( \int_{\varepsilon<|y|<\varepsilon^{-1}} \mathfrak{K}^{\phi}_{d}(z,y) a(y)V(y)\mathrm{d}y  \right) \phi(\Delta)u(x-z) \mathrm{d}z   
\\
&:= \int_{\bR^{d}} \mathcal{K}^{a,\varepsilon,\phi,V}_{d}(z) \phi(\Delta)u(x-z) \mathrm{d}z.
\end{aligned}
\end{equation}
To complete the proof, it suffices to check that the integral operator in \eqref{eqn 10.23.11:27} and its kernel $\mathcal{K}_d^{a,\varepsilon,\phi,V}$ satisfy the assumptions in Theorem \ref{CZ}.

\textbf{Step 1.} $L_2$-estimation.

By Plancherel's theorem,
\begin{equation}
\label{24.06.25.19.12}
\left\|\int_{\bR^d}\mathcal{K}_d^{a,\varepsilon,\phi,V}(\cdot-z)f(z)\mathrm{d}z\right\|_{L_2}^2=\int_{\bR^d}\left|\mathcal{F}_d\left[\mathcal{K}_d^{a,\varepsilon,\phi,V}\right](\xi)\right|^2\left|\mathcal{F}_d[f](\xi)\right|^2\mathrm{d}\xi.
\end{equation}
If we have
\begin{equation}
\label{24.06.25.19.11}
\sup_{\xi\in\bR^d}\left|  \mathcal{F}_{d} \left[\mathcal{K}^{a,\varepsilon,\phi,V}_{d} \right] (\xi) \right| \leq C(c_0,c_1,C_{V},d,\delta_0,\delta_1,\|a\|_{L_{\infty}(\bR^d)}),
\end{equation}
then by combining \eqref{24.06.25.19.12}, \eqref{24.06.25.19.11} and Plancherel's theorem, we obtain the $L_2$-estimation.

Now we prove \eqref{24.06.25.19.11}.
Using \eqref{eqn 06.30.14:58} and the fact that $a$ is bounded, we can check that
\begin{align*}
\left| \mathcal{F}_{d} \left[\mathcal{K}^{a,\varepsilon,\phi,V}_{d} \right] (\xi) \right|
 \leq C_{V}\|a\|_{L_{\infty}(\bR^d)}\int_{\bR^{d}} \left|   \mathcal{F}_{d}\left[ \mathfrak{K}^{\phi}_{d}(\cdot,y) \right](\xi) \right| j_d^{\phi}(|y|) \mathrm{d}y.
\end{align*}
Since
\begin{align*}
\mathcal{F}_{d}\left[\mathfrak{K}^{\phi}_{d}(\cdot,y)  \right](\xi)
&= \mathcal{F}_{d}\left[K^{\phi}_{d}\right](\xi) - \mathcal{F}_{d}\left[K^{\phi}_{d}(\cdot+y)\right](\xi)= \left( 1 - \mathrm{e}^{i \xi\cdot y} \right) \mathcal{F}_{d}\left[ K^{\phi}_{d} \right](\xi)=  \frac{1- \mathrm{e}^{i \xi\cdot y}}{\phi(|\xi|^{2})},
\end{align*}
we have
\begin{align}\label{eqn 08.07.13:43}
\left|\mathcal{F}_{d}\left[\mathfrak{K}^{\phi}_{d}(\cdot,y)\right](\xi) \right|  \leq 
\begin{cases} |\xi||y|\phi(|\xi|^2)^{-1} & \quad  \text{if} \quad |\xi||y|\leq 1,
\\
2 \phi(|\xi|^2)^{-1} & \quad \text{if} \quad |\xi||y| \geq 1.
\end{cases}
\end{align}
Using \eqref{eqn 08.07.13:43}, we have
\begin{align*}
\left| \mathcal{F}_{d}\left[ \mathcal{K}^{a,\epsilon,\phi,V}_{d} \right](\xi) \right| &\leq C \left(|\xi|\phi(|\xi|^2)^{-1}\int_{|y|\leq |\xi|^{-1}}|y|  j_d^{\phi}(|y|) \mathrm{d}y
 +\phi(|\xi|^{2})^{-1}\int_{|y| \geq |\xi|^{-1}} j_d^{\phi}(|y|) \mathrm{d}y\right)\\
&\leq C \left(|\xi| \int_{0}^{|\xi|^{-1}} \frac{\phi(\rho^{-2})}{\phi(|\xi|^{2})} \mathrm{d}\rho +  \int_{|\xi|^{-1}}^{\infty} \frac{\phi(\rho^{-2})}{\phi(|\xi|^2)} \frac{\mathrm{d}\rho}{\rho}\right) \leq C,
\end{align*}
where for the last inequality we used \eqref{eqn 07.01.16:07} and \eqref{eqn 07.01.15:42}. Thus \eqref{24.06.25.19.11} is proved.

\textbf{Step 2.} Verifying the H\"ormander condition \eqref{eqn 02.19.21:49}.

In this step, we prove that
\begin{equation}
\label{24.06.25.19.21}
\begin{aligned}
\left| \mathcal{K}^{a,\varepsilon,\phi,V}_{d}(x+z) - \mathcal{K}^{a,\varepsilon,\phi,V}_{d}(x) \right|  \leq C \frac{|z|^{\bar{\delta}}}{|x|^{d+\bar{\delta}}},\quad |x|\geq 4|z|,
\end{aligned}
\end{equation}
where $\bar{\delta} = (1-2\delta_1) \wedge 2\delta_{0}$, and $C=C(c_0,c_1,C_{V},d,\delta_{0},\delta_1)$.
From the definition of $\mathcal{K}^{a,\varepsilon,\phi,V}_{d}$, we have
\begin{align*}
&\left| \mathcal{K}^{a,\varepsilon,\phi,V}_{d}(x+z) - \mathcal{K}^{a,\varepsilon,\phi,V}_{d}(x) \right|
\\
&\quad \leq C \int_{\bR^{d}} \left|   K^{\phi}_{d}(x+y+z) - K^{\phi}_{d}(x+y) - K^{\phi}_{d}(x+z) + K^{\phi}_{d}(x) \right| j_{d}^{\phi}(|y|) \mathrm{d}y
\\
&\quad \leq C \left( \int_{A_1(x,z)} \cdots  + \int_{A_2(x,z)} \cdots  + \int_{A_3(x,z)} \cdots + \int_{A_4(x,z)} \cdots   \right)
\\
&\quad := C \left( \mathcal{I}_{1}(x,z) + \mathcal{I}_{2}(x,z) + \mathcal{I}_{3}(x,z) + \mathcal{I}_{4}(x,z) \right),
\end{align*}
where
\begin{align*}
    &A_1(x,z):=\{y\in\bR^d:|y|\leq|z|\},\\
    &A_2(x,z):=\{y\in\bR^d:|z|<|y|,\quad |x|/2< |x+y|\},\\
    &A_3(x,z):=\{y\in\bR^d:|z|<|y|,\quad |x+y|\leq 2|z|\}\\
    &A_4(x,z):=\{y\in\bR^d:|z|<|y|,\quad 2|z|<|x+y|\leq |x|/2\}.
\end{align*}
To complete the proof, it suffices to prove that
$$
\max_{1\leq j\leq 4}\mathcal{I}_j(x,z)  \leq C\frac{|z|^{\bar{\delta}}}{|x|^{d+\bar{\delta}}}.
$$

\textbf{Step 2-1.} Estimation of $\mathcal{I}_1$.

By the fundamental theorem of calculus and Theorem \ref{thm 07.21.11:32},
\begin{equation}
\label{24.06.24.23.42}
\begin{aligned}
    &\left|   K^{\phi}_{d}(x+y+z) - K^{\phi}_{d}(x+y) - K^{\phi}_{d}(x+z) + K^{\phi}_{d}(x) \right|\\
    &\leq |K_d^{\phi}(x+y+z)-K_d^{\phi}(x+z)|+|K_d^{\phi}(x+y)-K_d^{\phi}(x)|\\
    &\leq\int_0^1\left(|\nabla_yK_d^{\phi}(x+ry+z)|+|\nabla_yK_d^{\phi}(x+ry)|\right)|y|\mathrm{d}r\\
    &\leq C|y|\int_0^1\left(\frac{1}{|x+ry+z|^{d+1}\phi(|x+ry+z|^{-2})}+\frac{1}{|x+ry|^{d+1}\phi(|x+ry|^{-2})}\right)\mathrm{d}r.
\end{aligned}
\end{equation}
Since $y\in A_1(x,z)$, for all $r\in[0,1]$,
\begin{align*}
    |y|\leq \frac{|x|}{2}\leq |x+z|-|y|\leq |x+ry+z|,\quad |y|\leq \frac{|x|}{2}\leq |x|-|y|\leq |x+ry|.
\end{align*}
Hence,
\begin{align*}
    &\left|   K^{\phi}_{d}(x+y+z) - K^{\phi}_{d}(x+y) - K^{\phi}_{d}(x+z) + K^{\phi}_{d}(x) \right|\\
    &\leq C\frac{|y|}{|x|^{d+1}}\int_0^1\left(\frac{1}{\phi(|x+ry+z|^{-2})}+\frac{1}{\phi(|x+ry|^{-2})}\right)\mathrm{d}r.
\end{align*}
Therefore,
\begin{align*}
\mathcal{I}_{1}(x,z)
&\leq C|x|^{-d-1} \int_{|y|\leq |z|}  \int_{0}^{1}\left(   \frac{\phi(|y|^{-2})}{\phi(|x+ry+z|^{-2})} + \frac{\phi(|y|^{-2})}{\phi(|x+ry|^{-2})}\right)  \mathrm{d}r  \frac{\mathrm{d}y}{|y|^{d-1}}
\\
&\leq C |x|^{-d-1} \int_{|y|\leq |z|}\int_{0}^{1}\left(   \frac{|x+ry+z|^{2\delta_1}+|x+ry|^{2\delta_1}}{|y|^{2\delta_1}} \right)  \mathrm{d}r \frac{\mathrm{d}y}{|y|^{d-1}}
\\
&\leq C \frac{1}{|x|^{d+(1-2\delta_1)}} \int_{0}^{|z|} \rho^{-2\delta_1} \mathrm{d}\rho = C \frac{|z|^{1-2\delta_1}}{|x|^{d+(1-2\delta_1)}}\leq C\frac{|z|^{\bar{\delta}}}{|x|^{d+\bar{\delta}}}.
\end{align*}

\textbf{Step 2-2.} Estimation of $\mathcal{I}_2$.

Using the similar argument in \eqref{24.06.24.23.42}, we have
\begin{equation}
\label{24.06.24.23.43}
\begin{aligned}
    &\left|   K^{\phi}_{d}(x+y+z) - K^{\phi}_{d}(x+y) - K^{\phi}_{d}(x+z) + K^{\phi}_{d}(x) \right|\\
    & \leq C|z|\int_0^1\left(\frac{1}{|x+y+rz|^{d+1}\phi(|x+y+rz|^{-2})}+\frac{1}{|x+rz|^{d+1}\phi(|x+rz|^{-2})}\right)\mathrm{d}r.
\end{aligned}
\end{equation}
Since $y\in A_2(x,z)$, for all $r\in[0,1]$, we have $|x+y+rz|,|x+rz|\geq |x|/4$.
Hence,
\begin{align*}
    &\left|   K^{\phi}_{d}(x+y+z) - K^{\phi}_{d}(x+y) - K^{\phi}_{d}(x+z) + K^{\phi}_{d}(x) \right|\\
    &\leq  C\frac{|z|}{\phi(16|x|^{-2})}\int_0^1\left(\frac{\phi(16|x|^{-2})}{|x+y+rz|^{d+1}\phi(|x+y+rz|^{-2})}+\frac{\phi(16|x|^{-2})}{|x+rz|^{d+1}\phi(|x+rz|^{-2})}\right)\mathrm{d}r\\
    & \leq C\frac{|z|}{\phi(|x|^{-2})}\int_0^1\left(\frac{|x|^{-2\delta_1}}{|x+y+rz|^{d+(1-2\delta_1)}}+\frac{|x|^{-2\delta_1}}{|x+rz|^{d+(1-2\delta_1)}}\right)\mathrm{d}r\\
    & \leq C\frac{|z|}{\phi(|x|^{-2})|x|^{d+1}}\leq C\frac{|z|^{1-2\delta_1}}{\phi(|z|^{-2})|x|^{d+(1-2\delta_1)}} \leq C\frac{|z|^{\bar{\delta}}}{\phi(|z|^{-2})|x|^{d+\bar{\delta}}}.
\end{align*}
Therefore, by \eqref{eqn 07.01.15:42},
\begin{equation}
\label{eqn 08.05.18:45}
\begin{aligned}
\mathcal{I}_{2}(x,z)
&\leq C\frac{|z|^{\bar{\delta}}}{\phi(|z|^{-2})|x|^{d+\bar{\delta}}}\int_{|z|\leq |y|} j_d^{\phi}(|y|)\mathrm{d}y \\
&\leq C\frac{|z|^{\bar{\delta}}}{|x|^{d+\bar{\delta}}\phi(|z|^{-2})}\int_{|z|}^{\infty}\rho^{-1}\phi(\rho^{-2})\mathrm{d}\rho\leq C\frac{|z|^{\bar{\delta}}}{|x|^{d+\bar{\delta}}}.
\end{aligned}
\end{equation}

\textbf{Step 2-3.} Estimation of $\mathcal{I}_3$. 

We decompose $\mathcal{I}_3$ into two parts;
$
\mathcal{I}_3(x,z)\leq \mathcal{I}_{3,1}(x,z)+\mathcal{I}_{3,2}(x,z),
$
where
\begin{align*}
    \mathcal{I}_{3,1}(x,z)&:=\int_{A_3(x,z)}\left(|K_d^{\phi}(x+y+z)|+|K_d^{\phi}(x+y)|\right) j_d^{\phi}(|y|)\mathrm{d}y,\\
    \mathcal{I}_{3,2}(x,z)&:=\int_{A_3(x,z)}|K_d^{\phi}(x+z)-K_d^{\phi}(x)|j_d^{\phi}(|y|)\mathrm{d}y.
\end{align*}
Since $y\in A_3(x,z)$, we see that $|y|\geq |x|/2$, and $|x+y+z|\leq 3|z|\leq |x|\leq 2|y|$.
Hence,
\begin{align*}
&\mathcal{I}_{3,1}(x,z)\\
&\leq C\int_{A_3(x,z)} \left(\frac{\phi(|y|^{-2})}{\phi(|x+y+z|^{-2})|x+y+z|^{d}}+\frac{\phi(|y|^{-2})}{\phi(|x+y|^{-2})|x+y|^{d}}\right) \frac{\mathrm{d}y}{|y|^{d}}
\\
&\leq \frac{C}{|x|^{d}} \left(\int_{|x+y+z|\leq 3|z|,\, |y|\geq |x|/2}   \frac{|y|^{-2\delta_{0}}}{|x+y+z|^{d-2\delta_{0}}}  \mathrm{d}y+\int_{|x+y|\leq 2|z|,\,|y|\geq|x|/2}   \frac{|y|^{-2\delta_{0}}}{|x+y|^{d-2\delta_{0}}}  \mathrm{d}y\right)
\\
&\leq \frac{C}{|x|^{d+2\delta_{0}}} \int_{0}^{3|z|} \rho^{-1+2\delta_{0}} \mathrm{d}\rho = C \frac{|z|^{2\delta_{0}}}{|x|^{d+2\delta_{0}}}\leq C\frac{|z|^{\bar{\delta}}}{|x|^{d+\bar{\delta}}}.
\end{align*}
Using \eqref{24.06.24.23.42} and \eqref{eqn 08.05.18:45}, we also have $\mathcal{I}_{3,2}(x,z) \leq C\frac{|z|^{\bar{\delta}}}{|x|^{d+\bar{\delta}}}$.
Therefore, $\mathcal{I}_{3}(x,z) \leq C\frac{|z|^{\bar{\delta}}}{|x|^{d+\bar{\delta}}}$.

\textbf{Step 2-4.} Estimation of $\mathcal{I}_4$.

If $y\in A_4(x,z)$, then for all $r\in[0,1]$, $|x+y+rz|>|z|$ and $|x+rz|,|y|>|x|/2$.
Using \eqref{24.06.24.23.43}, it follows that
\begin{align*}
    &\left|   K^{\phi}_{d}(x+y+z) - K^{\phi}_{d}(x+y) - K^{\phi}_{d}(x+z) + K^{\phi}_{d}(x) \right|\\
    &\leq  C|z|\int_0^1\frac{\phi(|z|^{-2})}{\phi(|z|^{-2})|x+y+rz|^{d+1}\phi(|x+y+rz|^{-2})} \mathrm{d}r\\
   &\quad + C|z| \int_{0}^{1} \frac{\phi(|x|^{-2})}{\phi(|x|^{-2})|x+rz|^{d+1}\phi(|x+rz|^{-2})})\mathrm{d}r\\
    &\leq C\int_0^1\left(\frac{|z|^{1-2\delta_1}}{\phi(|z|^{-2})|x+y+rz|^{d+(1-2\delta_1)}}+\frac{|z||x|^{-2\delta_1}}{\phi(|x|^{-2})|x+rz|^{d+(1-2\delta_1)}}\right)\mathrm{d}r\\
    &\leq C\int_0^1\frac{|z|^{1-2\delta_1}}{\phi(|z|^{-2})|x+y+rz|^{d+(1-2\delta_1)}}\mathrm{d}r+C\frac{|z|^{1-2\delta_1}}{\phi(|z|^{-2})|x|^{d+(1-2\delta_1)}}.
\end{align*}
Therefore, using \eqref{eqn 07.01.15:42},
\begin{align*}
\mathcal{I}_{4}(x,z) &\leq C|z|^{1-2\delta_1}  \int_{0}^{1}\int_{|z| < |y|, \, |z|<|x+y+rz|}\frac{1}{|x+y+rz|^{d+(1-2\delta_1)}}\frac{\phi(|y|^{-2})}{\phi(|z|^{-2})} \frac{\mathrm{d}y}{|y|^d}\mathrm{d}r
\\
&\quad + C\frac{|z|^{1-2\delta_1}}{\phi(|z|^{-2})|x|^{d+(1-2\delta_1)}} \int_{|z|< |y|} j_{d}^{\phi}(|y|)  \mathrm{d}y
\\
&\leq C |z|^{1-2\delta_1+2\delta_0}\int_{0}^{1} \int_{|z|<|y|,\,|z|<|x+y+rz|} \frac{|y|^{-d-2\delta_0}}{|x+y+rz|^{d+(1-2\delta_1)}}  \mathrm{d}y \mathrm{d}r
\\
&\quad + C\frac{|z|^{1-2\delta_1}}{\phi(|z|^{-2})|x|^{d+(1-2\delta_1)}} \int_{|z|}^{\infty} \rho^{-1} \phi(\rho^{-2})  \mathrm{d}\rho
\\
&\leq C \left(\frac{|z|^{1-2\delta_1+2\delta_0}}{|x|^{d+2\delta_0}} \int_{|z|}^{\infty} \rho^{-1-(1-2\delta_1)} \mathrm{d}\rho  +   \frac{|z|^{1-2\delta_1}}{|x|^{d+(1-2\delta_1)}}\right)
\\
&\leq C \left(\frac{|z|^{2\delta_0}}{|x|^{d+2\delta_0}} +   \frac{|z|^{1-2\delta_1}}{|x|^{d+(1-2\delta_1)}}\right) \leq C \frac{|z|^{\bar{\delta}}}{|x|^{d+\bar{\delta}}}.
\end{align*}
\eqref{24.06.25.19.21} is proved.
By Theorem \ref{24.07.04.17.53}, we obtain the desired result.
\end{proof}

\begin{lem}
\label{lem 10.20.12:56}
Let $f:(0,\infty)\to(0,\infty)$ satisfy
\begin{equation}
\label{eqn 10.20.12:58}
c_1
\left(\frac Rr\right)^{\nu_1}
\leq
\frac{f(R)}{f(r)}
\leq
c_2
\left(\frac Rr\right)^{\nu_2},
\qquad
0<r\leq R<\infty,
\end{equation}
where $0<\nu_1\leq\nu_2<1$. 
Then there exists a complete Bernstein function $\Phi_f$ such that
\begin{equation}
\label{eqn 10.20.13:25}
C^{-1}f(r)
\leq
\Phi_f(r)
\leq
Cf(r),
\qquad r>0.
\end{equation}
The constant $C$ depends only on $c_1,c_2,\nu_1,\nu_2$.
\end{lem}
\begin{proof}
We adopt the proof of \cite[Lemma 3.2]{CKKW18}. 
Define
$$
\tilde{f}(r) := \int_{0}^{\infty} \mathrm{e}^{-rs}f(s) \mathrm{d}s,
\quad
\Phi_{f}(r) := \int_{0}^{\infty} \frac{r}{r+s} \frac{f(s)}{s} \mathrm{d}s.
$$
Then, using the relation
$$
\int_{0}^{\infty} (1-\mathrm{e}^{-ru})\mathrm{e}^{-su} \mathrm{d}u = \frac{r}{s(r+s)},
$$
and by Fubini's theorem, we have
\begin{align*}
\Phi_{f}(r) &= \int_{0}^{\infty}  \frac{r}{r+s} \frac{f(s)}{s} \mathrm{d}s
 = \int_{0}^{\infty} \int_{0}^{\infty} (1-\mathrm{e}^{-ru})\mathrm{e}^{-su} \mathrm{d}u f(s) \mathrm{d}s
 = \int_{0}^{\infty} (1-\mathrm{e}^{-ru}) \tilde{f}(u) \mathrm{d}u.
\end{align*}
Note that $\tilde{f}$ is a completely monotone function (see \cite[Theorem 1.4]{schbern}). 
Also, using \eqref{eqn 10.20.12:58} and the change of variable $rs \mapsto s$, one can observe that
\begin{align*}
\int_{0}^{\infty} (1\wedge r) \tilde{f}(r) \mathrm{d}r
&= \int_{0}^{1} r \int_{0}^{\infty} \mathrm{e}^{-rs} f(s) \mathrm{d}s \mathrm{d}r + \int_{1}^{\infty} \int_{0}^{\infty} \mathrm{e}^{-rs}f(s) \mathrm{d}s \mathrm{d}r
\\
& \leq C \int_{0}^{1} r \left( \int_{0}^{1} \mathrm{e}^{-rs} s^{\nu_{1}} \mathrm{d}s  +  \int_{1}^{\infty} \mathrm{e}^{-rs}s^{\nu_{2}} \mathrm{d}s \right)\mathrm{d}r
\\
& \quad + C \int_{1}^{\infty} \left( \int_{0}^{1} \mathrm{e}^{-rs}s^{\nu_{1}} \mathrm{d}s  + \int_{1}^{\infty} \mathrm{e}^{-rs}s^{\nu_{2}} \mathrm{d}s \right)\mathrm{d}r
\\
& \leq C \int_{0}^{1} \left( r^{-\nu_{1}}\int_{0}^{r} \mathrm{e}^{-s}s^{\nu_{1}} \mathrm{d}s  + r^{-\nu_{2}} \int_{r}^{\infty} \mathrm{e}^{-s}s^{\nu_{2}} \mathrm{d}s \right)\mathrm{d}r
\\
& \quad + C \int_{1}^{\infty} \left( r^{-\nu_{1}-1}\int_{0}^{r} \mathrm{e}^{-s}s^{\nu_{1}} \mathrm{d}s  + r^{-\nu_{2}-1} \int_{r}^{\infty} \mathrm{e}^{-s}s^{\nu_{2}} \mathrm{d}s \right)\mathrm{d}r
\\
& \leq C  \left( \int_{0}^{1} (r^{-\nu_{1}} + r^{-\nu_{2}}) \mathrm{d}r + \int_{1}^{\infty} (r^{-\nu_{1}-1} + r^{-\nu_{2}-1} )\mathrm{d}r\right) < \infty,
\end{align*}
which implies that the measure $\tilde{f}(r)\mathrm{d}r$ is a L\'evy measure with the completely monotone density function $\tilde{f}$ with respcet to the Lebesgue measure. 
Therefore, $\Phi_{f}$ is a complete Bernstein function by definition (see \textit{e.g.} \cite[Definition 6.1]{schbern}).

Now it only remains to show \eqref{eqn 10.20.13:25}. Using \eqref{eqn 10.20.12:58},
\begin{align*}
\Phi_{f}(r) = \int_{0}^{\infty} \frac{r}{r+s} \frac{f(s)}{s} \mathrm{d}s \geq \frac{1}{2} \int_{0}^{r} \frac{f(s)}{s} \mathrm{d}s
\geq \frac{f(r)}{2c_{2}}\int_{0}^{r} r^{-\nu_{2}} s^{\nu_{2}-1} \mathrm{d}s = \frac{1}{2c_{2}\nu_{2}}f(r).  
\end{align*}
Similarly, we have
\begin{align*}
\Phi_{f}(r) &\leq  \int_{0}^{r} \frac{f(s)}{s} \mathrm{d}s + r\int_{r}^{\infty} \frac{f(s)}{s^{2}}
\\
& \leq \frac{f(r)}{c_{1}}\int_{0}^{r} r^{-\nu_{1}} s^{\nu_{1}-1} \mathrm{d}s + c_{2}f(r) \int_{r}^{\infty} r^{1-\nu_{2}}s^{-2+\nu_{2}} \mathrm{d}s 
\\
&\leq \left( \frac{1}{c_{1}\nu_{1}} + c_{2}(1-\nu_{2}) \right) f(r).
\end{align*}
The lemma is proved.
\end{proof}

\begin{proof}[Proof of Theorem \ref{thm 01.15.15:07}]

\textbf{Case 1.} $0<\delta_{0} \leq \delta_1<1/2$.

By Fatou's lemma and Lemma \ref{24.07.04.18.22},
\begin{align*}
    &\left\|\int_{\bR^d}(u(\cdot+y)-u(\cdot))a(y)V(y)\mathrm{d}y\right\|_{L_p(w)}
\\
    &\leq \liminf_{\varepsilon\downarrow0}\left \| \int_{\varepsilon<|y|<\varepsilon^{-1}} \left( u(\cdot+y) - u(\cdot)\right) a(y) V(y)\mathrm{d}y \right\|_{L_{p}(w)}
\\
    &\leq\sup_{\varepsilon\in(0,1)}\left \| \int_{\varepsilon<|y|<\varepsilon^{-1}} \left( u(\cdot+y) - u(\cdot)\right) a(y) V(y)\mathrm{d}y \right\|_{L_{p}(w)}\leq C\|\phi(\Delta)u\|_{L_p(w)}.
\end{align*}

\textbf{Case 2.} $1/2 < \delta_{0} \leq \delta_{1} <1$.

By the fundamental theorem of calculus, it follows that
\begin{align*}
&\int_{\bR^{d}} \left( u(x+y) - u(x) - y\cdot \nabla u(x) \right)  a(y) V(y) \mathrm{d}y 
\\
&= \int_{0}^{1} \left( \int_{\bR^{d}} \left( \nabla u(x+sy) - \nabla u(x) \right) \cdot y \, a(y) V(y) \mathrm{d}y \right) \mathrm{d}s  
\\
&= \int_{0}^{1} s^{-d} \left( \int_{\bR^{d}} \left( \nabla u(x+y) - \nabla u(x) \right) \cdot (s^{-1}y) \, a(s^{-1}y) V(s^{-1}y) \mathrm{d}y \right) \mathrm{d}s  
\\
&=   \int_{\bR^{d}} \left( \nabla u(x+y) - \nabla u(x) \right) \cdot y \,\left( \int_{0}^{1} s^{-d-1} a(s^{-1}y) \frac{V(s^{-1}y)}{j_{d}^{\phi}(|y|)}   \mathrm{d}s  \right) j_{d}^{\phi}(|y|)\mathrm{d}y  
\\
&= \sum_{i=1}^{d}  \int_{\bR^{d}} \left( u_{x^{i}}(x+y) -  u_{x^{i}}(x) \right) \left( \frac{y^{i}}{|y|} \int_{0}^{1} s^{-d-1} a(s^{-1}y) \frac{V(s^{-1}y)}{j_{d}^{\phi}(|y|)}   \mathrm{d}s  \right) |y| j_{d}^{\phi}(|y|)    \mathrm{d}y  
\\
&:= \sum_{i=1}^{d} \int_{\bR^{d}} \left( u_{x^{i}}(x+y) -  u_{x^{i}}(x) \right) a^{i}(y) \tilde{j}_{d}^{\phi}(|y|)   \mathrm{d}y,
\end{align*}
where $\tilde{j}_{d}^{\phi}(r) := rj_{d}^{\phi}(r)$.
By \eqref{eqn 02.25.17:23}, we have
\begin{equation*}\label{eqn 11.13.15:17}
\begin{aligned}
 |a^{i}(y)|& \leq  C \int_{0}^{1} s^{-d-1} |a(y/s)| \frac{s^{d}|y|^{-d}}{|y|^{-d}}\frac{\phi(s^{2}|y|^{-2})}{\phi(|y|^{-2})}  \mathrm{d}s\leq C \int_{0}^{1} s^{-1+2\delta_{0}} \mathrm{d}s \leq C,
\end{aligned}
\end{equation*}
and
$\tilde{j}_{d}^{\phi}(r) \simeq r^{-d} \tilde{\phi}(r^{-2})$, where $\tilde{\phi}(r):=r^{-1/2}\phi(r)$.
Observe that $\tilde{\phi}$ satisfies
$$
c_0 \left(\frac{R}{r}\right)^{\delta_{0}-1/2}  \leq \frac{\tilde{\phi}(R)}{\tilde{\phi}(r)}  \leq c_1 \left( \frac{R}{r} \right)^{\delta_1-1/2} \quad \forall \, 0<r\leq R<\infty.
$$
Using Lemma \ref{lem 10.20.12:56} to  $\tilde{\phi}$, there exists a complete Bernstein function $\Phi_{\tilde{\phi}}$ that satisfying $\tilde{\phi}(r) \simeq \Phi_{\tilde{\phi}}(r)$ for all $r>0$ and
\begin{align*}
\tilde{c}_0 \left(\frac{R}{r}\right)^{\delta_{0}-1/2}  \leq \frac{\Phi_{\tilde{\phi}}(R)}{\Phi_{\tilde{\phi}}(r)}  \leq \tilde{c}_1\left( \frac{R}{r} \right)^{\delta_{1}-1/2} \quad \forall \, 0<r\leq R<\infty,
\end{align*}
where the constants $\tilde{c}_{0},\tilde{c}_{1}$ depend only on $c_{0},c_{1},\delta_{0},\delta_{1}$. Therefore, we have $\tilde{j}_{d}^{\phi}(r) \simeq r^{-d} \Phi_{\tilde{\phi}}(r^{-2})$ for all $r>0$.
Since $\delta_{0}-1/2 \leq \delta_{1} -1/2 <1/2$, we can apply Lemma \ref{24.07.04.18.22} so that
\begin{align*}
&\left \| \int_{\bR^{d}} \left( u(\cdot+y) - u(\cdot) -y\cdot \nabla_x u(\cdot) \right) a(y) V(y) \mathrm{d}y \right\|_{L_{p}(w)}\\
&\leq\liminf_{\varepsilon\to0}\sum_{i=1}^{d}\left\|\int_{\varepsilon<|y|<\varepsilon^{-1}} \left( u_{x^{i}}(\cdot+y) -  u_{x^{i}}(\cdot) \right) a^{i}(y) \tilde{j}_{d}^{\phi}(|y|)   \mathrm{d}y\right\|_{L_{p}(w)} \leq C \|(\Phi_{\tilde{\phi}})(\Delta)\nabla u\|_{L_{p}(w)},
\end{align*}
where the constant $C$ depends only on $c_{0},c_{1}\delta_{0},\delta_{1},d,[w]_{A_p(\bR^d)}$.
Using \eqref{eqn 01.28.13:23}, we can check that 
$$
m^{i}(\xi):= \frac{\Phi_{\tilde{\phi}}(|\xi|^{2})\xi^{i}}{\phi(|\xi|^{2})} \quad i=1,\dots,d
$$
satisfies
$$
\left| D^{n}_{\xi}m^{i}(\xi) \right| \leq C |\xi|^{-n} \quad \forall\, \xi\in \bR^{d},  \quad \forall\, n\in \bN_{0}, \quad  \forall\, i.
$$
Therefore, by Mikhlin's multiplier theorem (\textit{e.g.} \cite[Theorem 4]{K80}) it follows that
\begin{equation*}\label{24.08.18.12.39}
\begin{aligned}
&\left \| \int_{\bR^{d}} \left( u(\cdot+y) - u(\cdot) -y\cdot \nabla u(\cdot) \right) a(y) V(y) \mathrm{d}y \right\|_{L_{p}(w)}
\\
&\leq C  \|\Phi_{\tilde{\phi}}(\Delta)\nabla u\|_{L_{p}(w)}
\\
&= C \|\Phi_{\tilde{\phi}}(\Delta)\nabla [ \phi^{-1}(\Delta) \phi(\Delta)] u\|_{L_{p}(w)}
\\
&\leq C \sum_{i=1}^{d} \|\Phi_{\tilde{\phi}}(\Delta)D_{x^{i}} \phi^{-1}(\Delta) \phi(\Delta)u\|_{L_{p}(w)} \leq C \|\phi(\Delta)u\|_{L_{p}(w)}.
\end{aligned}
\end{equation*}

\textbf{Case 3.} $a$ satisfies Assumption \ref{24.07.02.14.41}.

It suffices to show that
\begin{align}\label{eqn 07.10.15:47}
\sup_{\varepsilon \in (0,1)} \| \mathcal{L}^{a}_{V,\varepsilon} u \|_{L_{p}(w)} \leq C \|\phi(\Delta) u\|_{L_{p}(w)},
\end{align}
where 
$$
\mathcal{L}_{V,\varepsilon}^{a}u(x):=\int_{\varepsilon < |y| < \varepsilon^{-1}} \left( u(x+y) - u(x) -y\cdot \nabla u(x)\chi_{\delta_{0},\delta_{1}}(y) \right) a(y) V(y)\mathrm{d}y,
$$
and the constant $C=C(\|a\|_{L_{\infty}(\bR^d)},c_0,c_1,C_{V,\phi},d,\delta_{0},\delta_1,p,[w]_{A_{p}(\bR^d)})$ independent of $\varepsilon$. Let $\phi^{1/2}(r):=\left( \phi(r) \right)^{1/2}$.
One can check that $\phi^{1/2}$ satisfies
$$
c_0^{1/2}\left(\frac{R}{r}\right)^{\delta_0/2}\leq\frac{\phi^{1/2}(R)}{\phi^{1/2}(r)}\leq c_{1}^{1/2}\left(\frac{R}{r}\right)^{\delta_{1}/2},
$$
By Theorem \ref{thm 07.21.11:32}, there exists $K_d^{\phi^{1/2}}$ such that $\cF_d\left[K_d^{\phi^{1/2}}\right]=1/\phi(|\cdot|^2)^{1/2}$.  Hence, by following the argument from \eqref{eqn 01.16.12:30} to \eqref{eqn 10.23.11:27}, with Assumption \ref{24.07.02.14.41}, we have
\begin{align*}
\mathcal{L}^{a}_{V,\varepsilon}u(x) =\int_{\bR^{d}} \left(\int_{\varepsilon<|y|<\varepsilon^{-1}}\mathfrak{K}^{\phi^{1/2}}_{d}(z,y)a(y)V(y)\mathrm{d}y\right)(\phi^{1/2}(\Delta)u(x-z)-\phi^{1/2}(\Delta)u(x)) \mathrm{d}z,
\end{align*}
where $\mathfrak{K}_d^{\phi^{1/2}}(z,y):=K_d^{\phi^{1/2}}(z+y)-K_d^{\phi^{1/2}}(z)$.
By exploiting Proposition \ref{24.07.04.17.53} (ii), we define a function $m_{\varepsilon}\in L_{\infty}(\bR^d)$ with $\sup_{\varepsilon \in (0,1)} \|m_{\varepsilon}\|_{L_{\infty}(\mathbb{R}^{d})} <\infty$ satisfying
$$
\int_{\varepsilon<|y|<\varepsilon^{-1}}\mathfrak{K}^{\phi^{1/2}}_{d}(z,y)a(y)V(y)\mathrm{d}y=m_{\varepsilon}(z)j_d^{\phi^{1/2}}(|z|) := m_{\varepsilon}(z) \tilde{V}(z).
$$
Therefore, we have
\begin{align*}
\mathcal{L}^{a}_{V,\varepsilon}u(x)
 = \int_{\mathbb{R}^{d}} \left( \phi^{1/2}(\Delta)u(x-z) - \phi^{1/2}(\Delta)u(x) \right)  m_{\varepsilon}(z)\tilde{V}(z)\mathrm{d}z
:= \mathcal{L}^{m_{\varepsilon}}_{\tilde{V}} \left( \phi^{1/2}(\Delta)u \right)(x)
\end{align*}
Hence, if we define $\mathcal{L}^{m_{\varepsilon}}_{\tilde{V},\theta}$ ($\theta\in(0,1)$) like $\mathcal{L}^{a}_{V,\varepsilon}$, the following inequality
\begin{align}\label{eqn 07.10.16:18}
\sup_{\varepsilon,\theta \in (0,1)} \|\mathcal{L}^{m_{\varepsilon}}_{\tilde{V},\theta} (\phi^{1/2}(\Delta) u) \|_{L_{p}(w)} \leq C \|\phi (\Delta) u\|_{L_{p}(w)}
\end{align}
shows \eqref{eqn 07.10.15:47}. Using Theorem \ref{thm 07.21.11:32} again,

Since $\delta_{1}/2 < 1/2$, we can apply Lemma \ref{24.07.04.18.22} again to obtain 
\begin{align*}
\sup_{\varepsilon,\theta \in (0,1)} \|\mathcal{L}^{m_{\varepsilon}}_{\tilde{V},\theta}(\phi^{1/2}(\Delta)u)\|_{L_{p}(w)} \leq C \|\phi^{1/2}(\Delta)[\phi^{1/2}(\Delta)u] \|_{L_{p}(w)}  = C \|\phi(\Delta)u\|_{L_{p}(w)}.
\end{align*}
Hence, \eqref{eqn 07.10.16:18} follows, and the theorem is proved.
\end{proof}

\subsection{Localization and commutator estimates}

In this subsection we collect the estimates needed to pass from the $x$-independent coefficient case to the general case $a=a(x,y)$. 
The main point is to control the commutator terms produced by localization. 
The first ingredient is a Haj\l asz-type pointwise estimate adapted to our nonlocal operator.

\begin{prop}[Haj\l asz-type pointwise estimate]
\label{lem 10.30.10:35}
Let $\phi$ be a Bernstein function satisfying
Assumption \ref{asm 07.21.13:06}$(\delta_0,\delta_1)$ with
$0<\delta_0\leq\delta_1<1$. Let $\gamma\in(0,1]$ satisfy
$2\delta_1\gamma<1$. Then, for every $f\in C_c^\infty(\mathbb R^d)$ and
every $x,y\in\mathbb R^d$ with $x\neq y$,
\begin{align}
\label{eqn 01.28.22:16}
|f(x)-f(y)|
\leq
C
\bigl(\phi(|x-y|^{-2})\bigr)^{-\gamma}
\left(
\mathbb M[\phi^\gamma(\Delta)f](x)
+
\mathbb M[\phi^\gamma(\Delta)f](y)
\right),
\end{align}
where $C=C(c_0,c_1,\delta_0,\delta_1,d,\gamma)$.
\end{prop}
\begin{proof}
For $x,y \in \mathbb{R}^{d}$, and a function $g \in C^{\infty}_{c}(\mathbb{R}^{d})$, define
\begin{align*}
\mathcal{D}[g](x,y)
=:\int_{\bR^{d}} (\phi(|x-y|^{-2}))^{\gamma} \left|K^{\phi^{\gamma}}_{d}(x-z) - K^{\phi^{\gamma}}_{d}(y-z)  \right| \left| g(z) \right|\mathrm{d}z.
\end{align*}
Then using \eqref{eqn 01.16.12:30} one can check that
\begin{align*}
(\phi(|x-y|^{-2}))^{\gamma}|| f(x) - f(y) | \leq \mathcal{D}[\phi(\Delta)^{\gamma}f](x,y).
\end{align*}
Therefore, to prove \eqref{eqn 01.28.22:16}, it suffices to show
\begin{align}\label{eqn 01.28.22:53}
\mathcal{D}[g](x,y) \leq C \left( \mathbb{M}g (x) + \mathbb{M}g(y) \right).
\end{align}
We split $\mathcal{D}[g](x,y)$ into the following two parts
\begin{align*}\label{eqn 10.30.10:37}
\int_{  |x-z| \geq 2|x-y|} \dots \mathrm{d}z + \int_{ |x-z| \leq 2|x-y|} \dots \mathrm{d}z.
:= \mathcal{D}_{1}[g](x,y) + \mathcal{D}_{2}[g](x,y).
\end{align*}
For $\mathcal{D}_{1}[g](x,y)$, letting $\theta(s,x,y) = (x-z) + s(y-x)$, we can observe that
\begin{align*}
\mathcal{D}_{1}[g](x,y) &\leq C \int_{|x-z| \geq 2|x-y|} \int_{0}^{1} |\theta(s,x,y)|^{-d-1}\left( \frac{\phi(|x-y|^{-2})}{\phi(|\theta(s,x,y)|^{-2})} \right)^{\gamma} \mathrm{d}s |x-y| |g(z)| \mathrm{d}z
\end{align*}
Since 
\begin{align*}
|\theta(s,x,y)| \geq |x-z| - |y-x| \geq |x-z|/2 \quad \forall \, |x-z| \geq 2|x-y| ,\quad s\in[0,1],
\end{align*} 
we have
\begin{align*}
\mathcal{D}_{1}[g](x,y)  &\leq C \int_{|x-z| \geq 2|x-y|} |x-z|^{-d-1} \left(\frac{\phi(|x-y|^{-2})}{\phi(|x-z|^{-2})} \right)^{\gamma}|x-y| |g(z)| \mathrm{d}z
\\
& \leq  C \int_{|x-z| \geq 2|x-y|} |x-z|^{-d-1} \left( \frac{|x-y|^{-2}}{|x-z|^{-2}} \right)^{\delta_{1}\gamma} |x-y| |g(z)| \mathrm{d}z
\\
&\leq C \int_{|z| \geq 2|x-y|} |z|^{-d-1} \left( \frac{|x-y|^{-2}}{|z|^{-2}} \right)^{\delta_{1}\gamma} |x-y| |g(x-z)| \mathrm{d}z
\\
&\leq C |x-y|^{1-2\delta_{1}\gamma} \int_{|z|\geq 2|x-y|} |z|^{-d-1+2\delta_{1}\gamma} |g(x-z)| \mathrm{d}z.
\end{align*}
Using the integration by parts formula
\begin{equation}
\label{eqn 08.13.17:03}
\begin{aligned}
&\int_{r \leq |z| \leq R} F(|z|) G(z) \mathrm{d}z
\\
=& \int_{r}^{R} F(\rho) \int_{\partial B_{\rho}} G(\sigma)\mathrm{d}\sigma \mathrm{d}\rho = \int_{r}^{R} F(\rho) \frac{\mathrm{d}}{\mathrm{d}\rho}\left(\int_{\partial B_{\rho}} G(z)\mathrm{d}z \right) \mathrm{d}\rho   
\\
=& -\int_{r}^{R} F'(\rho) \left( \int_{B_{\rho}}G(z) \mathrm{d}z \right) \mathrm{d}\rho
 +  F(R) \int_{B_{R}}G(z) \mathrm{d}z  -F(r) \int_{B_{r}}G(z) \mathrm{d}z
\end{aligned}
\end{equation}
we see that (put $F(\rho)=\rho^{-d-1+2\delta_{1}\gamma}$, $r=2|x-y|$, $R\to\infty$)
\begin{equation*}\label{eqn 10.29.22:25}
\begin{aligned}
& |x-y|^{1-2\delta_{1}\gamma} \int_{|z|\geq 2|x-y|} |z|^{-d-1+2\delta_{1}\gamma} |g(x-z)| \mathrm{d}z 
\\
&\leq C |x-y|^{1-2\delta_{1}\gamma} \int_{2|x-y|}^{\infty} \rho^{-d-2+2\delta_{1}\gamma} \int_{B_{\rho}} |g(x-z)| \mathrm{d}z \mathrm{d}\rho   
\\
&\quad + C|x-y|^{1-2\delta_{1}\gamma}|x-y|^{-d-1+2\delta_{1}\gamma} \int_{B_{2|x-y|}} |g(x-z)|\mathrm{d}z  
\\
&\leq C |x-y|^{1-2\delta_{1}\gamma}\int_{2|x-y|}^{\infty} \rho^{-d-2+2\delta_{1}\gamma} \int_{B_{\rho}(x)} |g(z)| \mathrm{d}z \mathrm{d}\rho  
\\
&\quad + C|x-y|^{1-2\delta_{1}\gamma}|x-y|^{-d-1+2\delta_{1}\gamma} \int_{B_{2|x-y|}(x)} |g(z)|\mathrm{d}z   
\\
&\leq C|x-y|^{1-2\delta_{1}\gamma} \mathbb{M}g(x) \int_{2|x-y|}^{\infty} \rho^{-2+2\delta_{1}\gamma}\mathrm{d}\rho + C\mathbb{M}g(x) \leq C \mathbb{M}g(x).
\end{aligned}
\end{equation*}
Thus we have 
\begin{align*}
\mathcal{D}_{1}[g](x,y) \leq C \mathbb{M}g(x).
\end{align*}
For $\mathcal{D}_{2}[g](x,y)$ observe that $\mathcal{D}_{2}[g](x,y)$ is bounded by sum of the following two terms
\begin{align*}
\mathcal{D}_{2,1}[g](x,y)=:\int_{|x-z|\leq 2|x-y|} (\phi(|x-y|^{-2}))^{\gamma}  |K^{\phi^{\gamma}}_{d}(x-z)| |g(z)| \mathrm{d}z,
\end{align*}
and
\begin{align*}
\mathcal{D}_{2,2}[g](x,y)=:\int_{|x-z|\leq 2|x-y|} (\phi(|x-y|^{-2}))^{\gamma}  |K^{\phi^{\gamma}}_{d}(y-z)| |g(z)| \mathrm{d}z.
\end{align*}
For $\mathcal{D}_{2,1}[g](x,y)$, using \eqref{eqn 08.13.17:03} with $F(\rho) = \rho^{-d+2\delta_{0}\gamma}$, $R=2|y|$,  $r\downarrow 0$,
\begin{align*}
\mathcal{D}_{2,1}[g](x,y)
&\leq C \int_{|x-z|\leq 2|x-y|} |x-z|^{-d} \left( \frac{\phi(|x-y|^{-2})}{\phi(|x-z|^{-2})} \right)^{\gamma} |g(z)| \mathrm{d}z
\\
&\leq C \int_{|x-z| \leq 2|x-y|} |x-z|^{-d} \left( \frac{|x-y|^{-2}}{|x-z|^{-2}}  \right)^{\delta_{0}\gamma} |g(z)| \mathrm{d}z
\\
&\leq C \int_{|z| \leq 2|x-y|} |z|^{-d} \left( \frac{|x-y|^{-2}}{|z|^{-2}}  \right)^{\delta_{0}\gamma} |g(x-z)| \mathrm{d}z
\\
&\leq C  |x-y|^{-2\delta_{0}\gamma}\int_{0}^{2|x-y|} \rho^{-d-1+2\delta_{0}\gamma} \int_{B_{\rho}}|g(x-z)| \mathrm{d}z \mathrm{d}\rho
\\
&\quad + C |x-y|^{-2\delta_{0}\gamma} |x-y|^{-d+2\delta_{0}\gamma} \int_{B_{2|x-y|}} |g(x-z)|\mathrm{d}z
\\
&\leq C  |x-y|^{-2\delta_{0}\gamma} \mathbb{M}g(x) \int_{0}^{2|x-y|} \rho^{-1+2\delta_{0}\gamma} \mathrm{d}\rho + C\mathbb{M}f(x)\leq C \mathbb{M}g(x).
\end{align*}
Doing the similar calculation we have the following control of $\mathcal{D}_{2,2}[g](x,y)$;
\begin{align*}
\mathcal{D}_{2,2}[g](x,y)
&\leq C \int_{|x-z|\leq 2|x-y|} |y-z|^{-d} \left( \frac{\phi(|x-y|^{-2})}{\phi(|y-z|^{-2})} \right)^{\gamma} |g(z)| \mathrm{d}z
\\
& \leq C \int_{|y-z|\leq 3|x-y|} |y-z|^{-d} \left( \frac{\phi(|x-y|^{-2})}{\phi(|y-z|^{-2})} \right)^{\gamma} |g(z)| \mathrm{d}z
\\
&\leq C \int_{|y-z| \leq 3|x-y|} |y-z|^{-d} \left( \frac{|x-y|^{-2}}{|y-z|^{-2}}  \right)^{\delta_{0}\gamma} |g(z)| \mathrm{d}z
\\
&\leq C \int_{|z| \leq 3|x-y|} |z|^{-d} \left( \frac{|x-y|^{-2}}{|z|^{-2}}  \right)^{\delta_{0}\gamma} |g(y-z)| \mathrm{d}z
\\
&\leq C  |x-y|^{-2\delta_{0}\gamma}\int_{0}^{3|x-y|} \rho^{-d-1+2\delta_{0}\gamma} \int_{B_{\rho}}|g(y-z)| \mathrm{d}z \mathrm{d}\rho
\\
&\quad + C |x-y|^{-2\delta_{0}\gamma}|x-y|^{-d+2\delta_{0}\gamma}\int_{B_{3|x-y|}}|g(y-z)|\mathrm{d}z
\\
&\leq C |x-y|^{-2\delta_{0}\gamma} \mathbb{M}g(y) \int_{0}^{3|x-y|} \rho^{-1+2\delta_{0}\gamma} \mathrm{d}\rho +  C\mathbb{M}g(y)\leq C \mathbb{M}g(y).
\end{align*}
Hence, we have \eqref{eqn 01.28.22:53}, and the proposition is proved.
\end{proof}

\begin{corollary}
\label{thm 07.17.15:37}
Let $\phi$ and $\gamma$ be as in Proposition \ref{lem 10.30.10:35}. Let
$p\in(1,\infty)$ and $w\in A_p(\mathbb R^d)$. Then, for every
$u\in C_c^\infty(\mathbb R^d)$ and $y\in\mathbb R^d\setminus\{0\}$,
\begin{align*}
\|u(\cdot+y)-u(\cdot)\|_{L_p(w)}
\leq
C
\bigl(\phi(|y|^{-2})\bigr)^{-\gamma}
\left(
\|\phi^\gamma(\Delta)u\|_{L_p(w)}
+
\|\phi^\gamma(\Delta)u(\cdot+y)\|_{L_p(w)}
\right).
\end{align*}
In particular, if $w\equiv1$, then
$$
\|u(\cdot+y)-u(\cdot)\|_{L_p}
\leq
C
\bigl(\phi(|y|^{-2})\bigr)^{-\gamma}
\|\phi^\gamma(\Delta)u\|_{L_p}.
$$
\end{corollary}

\begin{proof}
Apply Proposition \ref{lem 10.30.10:35} with $x$ and $x+y$, and then use the boundedness of the Hardy--Littlewood maximal operator on $L_p(w)$.
For the case $w \equiv 1$, use the translation invariance of the Lebesgue measure.
\end{proof}
\begin{rem}
Note that, in Proposition \ref{lem 10.30.10:35} and Corollary \ref{thm 07.17.15:37}, one can take $\gamma=1$ when $\delta_{1} <1/2$. Hence,  it follows that
\begin{align*}
\|u(\cdot+y) - u(\cdot) \|_{L_{p}} \leq C (\phi(|y|^{-2}))^{-1} \|\phi(\Delta) u\|_{L_{p}} \quad \text{if} \quad \delta_{1}<1/2.
\end{align*}
In particular,
\begin{align*}
\|u(\cdot+y) - u(\cdot) \|_{L_{p}} \leq C |y|^{2\delta_{1}} \|\Delta^{\delta_{1}} u\|_{L_{p}} \quad \forall\, \delta_{1}<1/2
\end{align*}
holds.
\end{rem}

We shall also use the following small-jump estimate.

\begin{lem}
\label{lem 08.12.18:53}
Let $\phi$ be a Bernstein function satisfying
Assumption \ref{asm 07.21.13:06}$(\delta_0,\delta_1)$ with
$0<\delta_0\leq\delta_1<1$. Let $p\in(1,\infty)$ and
$w\in A_p(\mathbb R^d)$. Suppose that
$$
0\vee \left(1-\frac{1}{2\delta_1}\right)<\gamma<\frac12.
$$
Then, for every $u\in C_c^\infty(\mathbb R^d)$,
\begin{align*}
\left\|
\int_{|y|\leq1}
|u(\cdot+y)-u(\cdot)|
|y|j_d^\phi(|y|)\,\mathrm dy
\right\|_{L_p(w)}
\leq
C
\|\phi^\gamma(\Delta)u\|_{L_p(w)}.
\end{align*}
\end{lem}
\begin{proof}

From \eqref{eqn 01.16.12:30}, and Proposition \ref{lem 10.30.10:35} we observe that
\begin{align*}
|u(x+y) - u(x)| &\leq \left( \phi(|y|^{-2})^{-\gamma} \right)\mathcal{D}[\phi^{\gamma}(\Delta)u](x+y,y)
\\
& \leq C \left( \phi(|y|^{-2})^{-\gamma} \right) \left( \mathbb{M}\phi^{\gamma}(\Delta)u(x) + \mathbb{M}\phi^{\gamma}(\Delta)u(x+y)   \right),
\end{align*}
where the constant $C$ is independent of $x,y$. Therefore, we have
\begin{align*}
&\int_{|y|\leq 1} \left| u(x+y) - u(x) \right| |y| j^{\phi}_{d}(|y|)\mathrm{d}y
\\
&\leq C \int_{|y|\leq 1} |y| \left( \phi(|y|^{-2}) \right)^{1-\gamma} \left( \mathbb{M}\phi^{\gamma}(\Delta)u(x) + \mathbb{M}\phi^{\gamma}(\Delta)u(x+y)   \right) |y|^{-d} \mathrm{d}y
\\
&\leq C \sum_{j=0}^{\infty} \int_{2^{-j-1} < |y| \leq 2^{-j}} |y| \left( \phi(|y|^{-2}) \right)^{1-\gamma} \left( \mathbb{M}\phi^{\gamma}(\Delta)u(x) + \mathbb{M}\phi^{\gamma}(\Delta)u(x+y)   \right) |y|^{-d} \mathrm{d}y
\\
&\leq C \sum_{j=0}^{\infty} \int_{2^{-j-1} <|y| \leq 2^{-j}} 2^{-j} 2^{2\delta_{1}(1-\gamma)j} |y|^{-d}  \left( \mathbb{M}\phi^{\gamma}(\Delta)u(x) + \mathbb{M}\phi^{\gamma}(\Delta)u(x+y)   \right)  \mathrm{d}y
\\
&\leq C \sum_{j=0}^{\infty} 2^{[2\delta_{1}(1-\gamma)-1]j} \frac{1}{|B_{2^{-j}}|}\int_{B_{2^{-j}}}  \left( \mathbb{M}\phi^{\gamma}(\Delta)u(x) + \mathbb{M}\phi^{\gamma}(\Delta)u(x+y)   \right)  \mathrm{d}y
\\
&\leq C \sum_{j=0}^{\infty} 2^{[2\delta_{1}(1-\gamma)-1]j} \left(  \mathbb{M}\phi^{\gamma}(\Delta)u(x) + \mathbb{M}\mathbb{M}\phi^{\gamma}(\Delta)u(x)  \right)
\\
&\leq C \left(  \mathbb{M}\phi^{\gamma}(\Delta)u(x) + \mathbb{M}\mathbb{M}\phi^{\gamma}(\Delta)u(x)  \right),
\end{align*}
where the last summation converges due to the choice of $\gamma \in (1-1/2\delta_{1}, 1/2)$ (recall $\delta_{1} <1$). Therefore, it follows that
\begin{align*}
&\int_{\mathbb{R}^{d}} \left( \int_{|y| \leq 1}  \left|  u(x+y) - u(x) \right| |y| j^{\phi}_{d}(|y|) \mathrm{d}y  \right)^{p} w(x) \mathrm{d}x
\\
&\leq C  \int_{\mathbb{R}^{d}} \left(  \mathbb{M}\phi^{\gamma}(\Delta)u(x) + \mathbb{M}\mathbb{M}\phi^{\gamma}(\Delta)u(x)  \right)^{p} w(x) \mathrm{d}x
\\
&\leq C \|\phi^{\gamma}(\Delta)u\|^{p}_{L_{p}(w)}.
\end{align*}
The lemma is proved.
\end{proof}

\begin{rem}\label{rmk 01.28.14:36}
Since the functions
\begin{align*}
m^{\phi}_{\gamma}(\xi) = \frac{(1+\phi(|\xi|^{2}))^{\gamma/2}}{1 + ( \phi(|\xi|^{2}) )^{\gamma/2}}, \quad \frac{(\phi(|\xi|^{2}))^{\gamma/2}}{(1+\phi(|\xi|^{2}))^{\gamma/2}}
\end{align*}
satisfy the Mikhlin multiplier condition (see the proof of \cite[Lemma 6.1 (iv)]{KKK13}, it follows that
\begin{align*}
\|u\|_{H^{\phi,\gamma}_{p}(w)} \simeq \|u\|_{L_{p}(w)} + \|\phi(\Delta)^{\gamma/2}u\|_{L_{p}(w)}.
\end{align*}
In particular, $\|u\|_{H^{\phi,2}_{p}(w)} \simeq \|u\|_{L_{p}(w)} + \|\phi(\Delta)u\|_{L_{p}(w)}$. Similarly, 
$$
\|u\|_{H^{\phi,\nu}_{p}(w)} \leq \|u\|_{H^{\phi,\gamma}_{p}(w)} \quad  \forall\, \nu \leq \gamma.
$$
\end{rem}
The next lemma controls the commutator between the nonlocal difference of
$u$ and a smooth cutoff function.

\begin{lem}
\label{lem 07.17.19:30}
Let $\phi$ be a Bernstein function satisfying
Assumption \ref{asm 07.21.13:06}$(\delta_0,\delta_1)$ with
$0<\delta_0\leq\delta_1<1$. Let $p\in(1,\infty)$ and
$w\in A_p(\mathbb R^d)$. Let $\gamma\in(0,1/2)$ satisfy
$$
1-\frac{1}{2\delta_1}<\gamma<\frac12.
$$
Then, for every $\varphi\in C_c^\infty(\mathbb R^d)$ and
$u\in C_c^\infty(\mathbb R^d)$,
\begin{equation*}
\begin{aligned}
&\int_{\mathbb R^d}\int_{\mathbb R^d}
\left(
\int_{\mathbb R^d}
|u(x+y)-u(x)|
|\varphi(x-z+y)-\varphi(x-z)|
j_d^\phi(|y|)\,\mathrm dy
\right)^p
w(x)\,\mathrm dx\,\mathrm dz\\
&\leq
C
\left(
\|\varphi\|_{L_p}^p+\|\nabla\varphi\|_{L_p}^p
\right)
\|u\|_{H_p^{\phi,2\gamma}(w)}^p .
\end{aligned}
\end{equation*}
\end{lem}
\begin{proof}
We decompose the integral with respect to $y$ into two parts
\begin{align*}
& \int_{\mathbb{R}^{d}}  |u(x+y) - u(x)|  | \varphi(x-z+y) - \varphi(x-z)| j^{\phi}_{d}(|y|)\mathrm{d}y 
\\
& = \int_{|y| > 1}   |u(x+y) - u(x)|  | \varphi(x-z+y) - \varphi(x-z)| j^{\phi}_{d}(|y|)\mathrm{d}y 
\\
& \quad + \int_{|y| \leq 1}  |u(x+y) - u(x)|  | \varphi(x-z+y) - \varphi(x-z)| j^{\phi}_{d}(|y|)\mathrm{d}y
\\
&:= I_{1}(x,z) + I_{2}(x,z).
\end{align*}
Observe that, by Minkowski's inequality
\begin{align*}
&\int_{\mathbb{R}^{d}} \int_{\mathbb{R}^{d}} I_{1}^{p}(x,z) w(x) \mathrm{d}x \mathrm{d}z
\leq 2\|\varphi\|^{p}_{L_{p}} \int_{\mathbb{R}^{d}} \left( \int_{|y| >1} |u(x+y) - u(x)| j^{\phi}_{d}(|y|)\mathrm{d}y \right)^{p} w(x)\mathrm{d}x.
\end{align*}
Also, by \eqref{eqn 07.01.15:42} we have
\begin{equation}
\label{eqn 07.17.16:46}
\begin{aligned}
&\int_{|y| >1} |u(x+y) - u(x)| j^{\phi}_{d}(|y|)\mathrm{d}y   
\\
&\leq C \left( \int_{|y| >1} |u(x+y)| \phi(|y|^{-2})|y|^{-d} \mathrm{d}y +   |u(x)| \int_{|y|>1} \phi(|y|^{-2})|y|^{-d}\mathrm{d}y \right)  
\\
&\leq C \left( \int_{|y| >1} |u(x+y)| \phi(|y|^{-2})|y|^{-d} \mathrm{d}y + |u(x)|\int_{1}^{\infty} \rho^{-1} \phi(\rho^{-2}) \mathrm{d}\rho \right)   
\\
&\leq C \left( \int_{|y| >1} |u(x+y)| \phi(|y|^{-2})|y|^{-d} \mathrm{d}y + |u(x)| \right).
\end{aligned}
\end{equation}
For each $j =0,1,\dots$,
\begin{equation}
\label{eqn 11.21.18:44}
\begin{aligned}
&\int_{2^{j} < |y| \leq 2^{j+1}} |u(x+y)||y|^{-d} \phi(|y|^{-2}) \mathrm{d}y   
\\
&\leq C\int_{2^{j} < |y| \leq 2^{j+1}} |u(x+y)||y|^{-d} |y|^{-2\delta_{1}} 2^{2(j+1)\delta_{1}}\phi(2^{-2(j+1)}) \mathrm{d}y  
\\
&\leq C \frac{1}{|B_{2^{j+1}}|} \int_{2^{j} < |y| \leq 2^{j+1}} |u(x+y)| \phi(2^{-2(j+1)}) \mathrm{d}y   
\\
&\leq C 2^{-2\delta_{0}(j+1)}\mathbb{M}u(x).
\end{aligned}
\end{equation}
Therefore, it follows that
\begin{align*}
\int_{|y| >1} |u(x+y)| \phi(|y|^{-2})|y|^{-d} \mathrm{d}y \leq C \left[\sum_{j=0}^{\infty} 2^{-2\delta_{0}(j+1)}\right] \mathbb{M}u(x) \leq C \mathbb{M}u(x).
\end{align*}
Plugging this to \eqref{eqn 07.17.16:46},
\begin{align*}
\int_{\mathbb{R}^{d}} \int_{\mathbb{R}^{d}} I_{1}(x,z)^{p} w(x) \mathrm{d}x\mathrm{d}z
\leq C  \|\varphi\|^{p}_{L_{p}} \int_{\mathbb{R}^{d}} \left( |u(x)|^{p} + |\mathbb{M}u(x)|^{p} \right) w(x) \mathrm{d}x
\leq C  \|\varphi\|^{p}_{L_{p}}  \|u\|^{p}_{L_{p}(w)}
\end{align*}
follows. Now, we consider $I_{2}$. By Minkowski's inequality, fundamental theorem of calculus and Lemma \ref{lem 08.12.18:53} (recall the choice of $\gamma$)
\begin{align*}
&\int_{\mathbb{R}^{d}} \int_{\mathbb{R}^{d}} I^{p}_{2}(x,z) w(x)\mathrm{d}x \mathrm{d}z
\\
&\leq  \int_{\mathbb{R}^{d}} \left(\int_{|y| \leq 1} \left( \int_{\mathbb{R}^{d}} |u(x+y) - u(x)|^{p} |\varphi(x-z+y) - \varphi(x-z)|^{p} \mathrm{d}z  \right)^{1/p} j^{\phi}_{d}(|y|)\mathrm{d}y \right)^{p} w(x)\mathrm{d}x
\\
&\leq  \|\nabla \varphi\|^{p}_{L_{p}} \int_{\mathbb{R}^{d}} \left( \int_{|y| \leq 1} |u(x+y) - u(x)| |y|   j^{\phi}_{d}(|y|)\mathrm{d}y \right)^{p} w(x) \mathrm{d}x
\\
&\leq C  \|\nabla \varphi\|^{p}_{L_{p}} \|\phi^{\gamma}(\Delta)u\|^{p}_{L_{p}(w)}.
\end{align*}
Therefore, we have
\begin{align*}
&\int_{\mathbb{R}^{d}} \int_{\mathbb{R}^{d}} \left( \int_{\mathbb{R}^{d}}  |u(x+y) - u(x)|  | \varphi(x-z+y) - \varphi(x-z)| j^{\phi}_{d}(|y|)\mathrm{d}y    \right)^{p} w(x)\mathrm{d}x \mathrm{d}z
\\
&\leq C (\|\varphi\|^{p}_{L_{p}}+ \|\nabla \varphi\|^{p}_{L_{p}})\left( \|u\|^{p}_{L_{p}(w)} +  \|\phi^{\gamma}(\Delta)u\|^{p}_{L_{p}(w)} \right).
\end{align*}
Finally, using Remark \ref{rmk 01.28.14:36}, we get the desired estimate. The lemma is proved.
\end{proof}

The following localized estimate is the only place where the continuity assumption on the coefficient in the spatial variable is used.

\begin{lem}
\label{lem 07.17.19:49}
Let $\phi$ be a Bernstein function satisfying
Assumption \ref{asm 07.21.13:06}$(\delta_0,\delta_1)$ with
$0<\delta_0\leq\delta_1<1$. Let $V$ be a nonnegative function satisfying
\eqref{eqn 02.11.14:19}, and let
$a:\mathbb R^d\times\mathbb R^d\to\mathbb C$ satisfy
Assumption \ref{24.08.18.13.17}. If
$\delta_0\leq1/2\leq\delta_1$, assume in addition  Assumption \ref{24.07.02.14.41}. Let
$\varphi\in C_c^\infty(B_\varepsilon)$ for some $\varepsilon\in(0,1)$.
Then, for every $p\in(1,\infty)$ and $w\in A_p(\mathbb R^d)$,
\begin{align}
\label{eqn 07.16.19:54}
\|\varphi\mathcal L_V^a u\|_{L_p(w)}
\leq
C(\mathcal C_1+\mathcal C_2)
\|\phi(\Delta)u\|_{L_p(w)},
\end{align}
where
\begin{align*}
\mathcal C_1
&:=
\|\varphi\|_{L_\infty}
\left(
\sup_{\substack{x\in B_{2\varepsilon}\\ y\in\mathbb R^d}}
|a(x,y)|
+
\varepsilon^\beta
\int_{|h|\leq\varepsilon}
\frac{\nu(|h|)}{|h|^{d+\beta}}\,\mathrm dh
\right)
+
\varepsilon\|\nabla\varphi\|_{L_\infty}
\sup_{\substack{x\in B_{2\varepsilon}\\ y\in\mathbb R^d}}
|a(x,y)|,
\end{align*}
and $
\mathcal{C}_2:=
\|\varphi\|_{L_\infty}
\sup_{x\in B_{\varepsilon}, y\in\mathbb R^d}
|a(x,y)|$.
Here, the constant $C$ depends only on
$c_0,c_1,C_V,d,\delta_0,\delta_1,p,[w]_{A_p(\mathbb R^d)}$ and the
structural constants in Assumption \ref{24.08.18.13.17}.
\end{lem}
\begin{proof}
Due to the Rubio de Francia extrapolation theorem, it suffices to prove the case $ \beta p >d$. 

\textbf{Case 1.} Assume that $a(x,y)$ is infinitely differentiable in $x$.
\\
 Set $\phi_{\beta}(r) =: r^{\beta/2}$ (\textit{i.e.} $\phi_{\beta}(\Delta)= \Delta^{\beta/2}$).   Since $\varphi$ does not depend on $y$, we have
\begin{equation}
\label{eqn 07.16.16:43}
\begin{aligned}
\varphi(x) \mathcal{L}^{a}_{V}u (x)   
&= \int_{\mathbb{R}^{d}} \left(  u(x+y) - u(x) - y \cdot \nabla u(x) \chi_{\delta_{0},\delta_{1}}  \right)  \varphi(x)a(x,y) V(y)\mathrm{d}y    
\\
&:= \int_{\mathbb{R}^{d}} \left(  u(x+y) - u(x) - y \cdot \nabla u(x) \chi_{\delta_{0},\delta_{1}}  \right)  \tilde{a}(x,y) V(y)\mathrm{d}y     
\end{aligned}
\end{equation}
Since $\tilde{a}(x,y) \in C^{\infty}_{c}(B_{\varepsilon})$ (in $x$), we rewrite
\begin{equation}
\label{eqn 07.16.14:44}
\begin{aligned}
\tilde{a}(x,y) =  [(\Delta)^{-\beta/2}\Delta^{\beta/2}\tilde{a}(\cdot,y)](x)
&= \int_{\mathbb{R}^{d}} K^{\phi_{\beta}}_{d}(x-z) [\Delta^{\beta/2}\tilde{a}(\cdot,y)](z) \mathrm{d}z    
\\
& := \int_{\mathbb{R}^{d}} K^{\phi_{\beta}}_{d}(x-z) \bar{a}(z,y) \mathrm{d}z
\end{aligned}
\end{equation}
Combining \eqref{eqn 07.16.16:43}, and \eqref{eqn 07.16.14:44}, we obtain
\begin{align*}
\varphi(x) \mathcal{L}^{a}_{V}u (x)
 &= \int_{\mathbb{R}^{d}} \left(  u(x+y) - u(x) - y \cdot \nabla u(x) \chi_{\delta_{0},\delta_{1}}  \right) \left(  \int_{\mathbb{R}^{d}} K^{\phi_{\beta}}_{d}(x-z) \bar{a}(z,y) \mathrm{d}z  \right) V(y) \mathrm{d}y
\\
& = \int_{\mathbb{R}^{d}} K^{\phi_{\beta}}_{d}(x-z) \left( \int_{\mathbb{R}^{d}}  \left(  u(x+y) - u(x) - y \cdot \nabla u(x) \chi_{\delta_{0},\delta_{1}}  \right)  \bar{a}(z,y) V(y) \mathrm{d}y \right) \mathrm{d}z
\\
& := \int_{\mathbb{R}^{d}} K^{\phi_{\beta}}_{d}(x-z) \mathcal{L}^{\bar{a}(z)}_{V}u(x) \mathrm{d}z.
\end{align*}
This implies that
\begin{align*}
\int_{\mathbb{R}^{d}} \left| \varphi(x)\mathcal{L}^{a}_{V}u(x)  \right|^{p} w(x) \mathrm{d}x  &= \int_{B_{\varepsilon}} \left| \varphi(x)\mathcal{L}^{a}_{V}u(x)  \right|^{p} w(x) \mathrm{d}x 
\\
& \leq C \int_{B_{\varepsilon}} \left|  \int_{\mathbb{R}^{d}} |x-z|^{\beta-d} \mathcal{L}^{\bar{a}(z)}_{V} u(x) \mathrm{d}z  \right|^{p} w(x) \mathrm{d}x
\\
&\leq C \int_{B_{\varepsilon}} \left|  \int_{B_{2\varepsilon}} |x-z|^{\beta-d} \mathcal{L}^{\bar{a}(z)}_{V} u(x) \mathrm{d}z  \right|^{p} w(x) \mathrm{d}x
\\
&\quad + C \int_{B_{\varepsilon}} \left|  \int_{\mathbb{R}^d\setminus B_{2\varepsilon}} |x-z|^{\beta-d} \mathcal{L}^{\bar{a}(z)}_{V} u(x) \mathrm{d}z  \right|^{p} w(x) \mathrm{d}x
 := C\left( I_{1} +I_{2}  \right).
\end{align*}
By H\"older's inequality (recall we assume $\beta p >d$),
\begin{align*}
I_{1} & \leq C \int_{B_{\varepsilon}} \left( \int_{B_{2\varepsilon}} \left( |x-z|^{\beta-d} \right)^{p/(p-1)} \mathrm{d}z \right)^{(p-1)} \left( \int_{B_{2\varepsilon}} \left|  \mathcal{L}^{\bar{a}(z)}_{V} u(x) \right|^{p} \mathrm{d}z \right) w(x) \mathrm{d}x
\\
& \leq C \int_{B_{\varepsilon}} \left( \int_{B_{3\varepsilon}} \left( |z|^{\beta-d} \right)^{p/(p-1)} \mathrm{d}z \right)^{(p-1)} \left( \int_{B_{2\varepsilon}} \left|  \mathcal{L}^{\bar{a}(z)}_{V} u(x) \right|^{p} \mathrm{d}z \right) w(x) \mathrm{d}x
\\
& \leq C  \int_{B_{\varepsilon}}   \left( \int_{0}^{3\varepsilon} \rho^{\frac{\beta p - d }{p-1}} \rho^{-1}\mathrm{d}\rho \right)^{(p-1)} \left( \int_{B_{2\varepsilon}} \left|  \mathcal{L}^{\bar{a}(z)}_{V} u(x) \right|^{p} \mathrm{d}z \right) w(x) \mathrm{d}x
\\
& \leq C \varepsilon^{\beta p -d} \int_{B_{\varepsilon}} \left( \int_{B_{2\varepsilon}} \left|  \mathcal{L}^{\bar{a}(z)}_{V}u(x) \right|^{p} \mathrm{d}z \right) w(x) \mathrm{d}x
\\
 & = C \varepsilon^{\beta p -d} \int_{B_{2\varepsilon}}   \int_{B_{\varepsilon}} \left|  \mathcal{L}^{\bar{a}(z)}_{V}u(x) \right|^{p}   w(x) \mathrm{d}x \mathrm{d}z.
\end{align*}
Recall that
\begin{align*}
\mathcal{L}^{\bar{a}(z)}_{V}u(x) = \int_{\mathbb{R}^{d}} \left( u(x+y) - u(x) - y \cdot \nabla u(x) \chi_{\delta_{0},\delta_{1}}  \right) [\Delta^{\beta/2}\varphi(\cdot)a(\cdot,y)](z) V(y) \mathrm{d}y.
\end{align*}
By inspecting the proof of \cite[Lemma 5.1]{DJK23nonlocal} (see p. 19 therein) $\sup_{|z| \leq 2\varepsilon, \, y\in \mathbb{R}^{d}} |\bar{a}(z,y)| \leq  C_{1}$,
where 
\begin{align*}
C_{1} = C_{1}(\varepsilon,\varphi,\nu)  &= C  \sup_{z\in \mathbb{R}^{d}} |\varphi(z)| \left( \sup_{\substack{z\in B_{2\varepsilon} \\ y\in\mathbb{R}^{d}}} |a(z,y)| \varepsilon^{-\beta}  + \int_{|y|\leq \varepsilon} \frac{\nu(|y|)}{|y|^{d+\beta}} \mathrm{d}y \right)
\\
&\quad + C \sup_{z\in \mathbb{R}^{d}} |\nabla \varphi (z)| \sup_{\substack{z\in B_{2\varepsilon} \\ y\in\mathbb{R}^{d}}} |a(z,y)| \varepsilon^{1-\beta}.  
\end{align*}
Observe that
\begin{align*}
I_{1} &\leq C\varepsilon ^{\beta p - d} \int_{B_{2\varepsilon}} \int_{B_{\varepsilon}} \left|  \mathcal{L}^{\bar{a}(z)}_{V}u(x) \right|^{p} w(x) \mathrm{d}x \mathrm{d}z
\\
& = C\varepsilon ^{\beta p - d} \int_{B_{2\varepsilon}} \int_{B_{\varepsilon}} \left| C_{1}(\varepsilon,\varphi,\nu) \mathcal{L}^{\bar{a}(z)/C_{1}}_{V}u(x) \right|^{p} w(x) \mathrm{d}x \mathrm{d}z
\end{align*}
Therefore, applying Theorem \ref{thm 01.15.15:07} (with $\bar{a}(z)/C_{1}$ in place of $a$), we have
\begin{align*}
I_{1} &\leq  C_{1}(\varepsilon,\varphi,\nu)^{p} \varepsilon^{\beta p -d} C \int_{B_{2\varepsilon}} \|\phi(\Delta)u\|^{p}_{L_{p}(w)} \mathrm{d}z \\
&= C C_{1}(\varepsilon,\varphi,\nu)^{p} \varepsilon^{\beta p} \|\phi(\Delta)u\|^{p}_{L_{p}(w)}
= C \mathcal{C}_{1}(\varepsilon,\varphi,\nu)^{p} \|\phi(\Delta)u\|^{p}_{L_{p}(w)},
\end{align*}
where
\begin{align*}
 \mathcal{C}_{1}(\varepsilon,\varphi,\nu)  
=&   \sup_{z\in \mathbb{R}^{d}} |\varphi(z)| \left( \sup_{z\in B_{2\varepsilon},y\in\mathbb{R}^{d}} |a(z,y)| +  \varepsilon^{\beta}   \int_{|y|\leq \varepsilon} \frac{\nu(|y|)}{|y|^{d+\beta}} \mathrm{d}y   \right)\\  
 &+\varepsilon \sup_{z\in \mathbb{R}^{d}} |\nabla \varphi (z)| \sup_{z\in B_{2\varepsilon},y\in\mathbb{R}^{d}} |a(z,y)|,
\end{align*}
and the constant $C$ is independent of $\varepsilon,\varphi$. For $I_{2}$, by Minkowski's inequality, 
\begin{equation}
\label{eqn 11.20.16:45}
\begin{aligned}
I_{2} &\leq \left(  \int_{B_{2\varepsilon}^{c}} \left( \int_{B_{\varepsilon}}  \left| |x-z|^{\beta - d} \mathcal{L}^{\bar{a}(z)}_{V} u(x) \right|^{p} w(x) \mathrm{d}x \right)^{1/p}  \mathrm{d}z \right) ^{p}   \nonumber
\\
& \leq C \left(  \int_{B_{2\varepsilon}^{c}} |z|^{\beta-d} \left( \int_{B_{\varepsilon}}  \left|  \mathcal{L}^{\bar{a}(z)}_{V} u(x) \right|^{p} w(x) \mathrm{d}x \right)^{1/p} \mathrm{d}z \right) ^{p}.
\end{aligned}
\end{equation}
Again inspecting the proof of \cite[Lemma 5.1]{DJK23nonlocal} (see p. 20 therein), we have
\begin{align*}
\sup_{|z| > 2\varepsilon, y\in \mathbb{R}^{d}} |\bar{a}(z,y)| \leq C \int_{|h| >\varepsilon} \sup_{y\in \mathbb{R}^{d}} \left| \varphi(z+h) a(y,z+h) \right| \frac{1}{|h|^{d+\beta}} \mathrm{d}h := C_{2}(\varepsilon,\varphi)
\end{align*}
Hence, applying Theorem \ref{thm 01.15.15:07} to \eqref{eqn 11.20.16:45} (with $\bar{a}/C_{2}$ in place of $a$), we have
\begin{align*}
I_{2} & \leq C \left(  \int_{B_{2\varepsilon}^{c}} |z|^{\beta-d} \left( \int_{B_{\varepsilon}}  \left|  \mathcal{L}^{\bar{a}(z)}_{V} u(x) \right|^{p} w(x) \mathrm{d}x \right)^{1/p} \mathrm{d}z \right) ^{p}
\\
& = C \left(  \int_{B_{2\varepsilon}^{c}} |z|^{\beta-d} \left( \int_{B_{\varepsilon}}  \left| C_{2} \mathcal{L}^{\bar{a}(z)/C_{2}}_{V} u(x) \right|^{p} w(x) \mathrm{d}x \right)^{1/p} \mathrm{d}z \right) ^{p}
\\
& \leq C \left(  \int_{B_{2\varepsilon}^{c}} |z|^{\beta-d} C_{2}(\varepsilon,\varphi) \|\phi(\Delta)u\|^{p}_{L_{p}(w)} \mathrm{d}z \right)^{p}.
\end{align*}
Observe that, since $\varphi(z+h)=0$ if $|z+h| > \varepsilon$,
\begin{align*}
&\int_{B_{2\varepsilon}^{c}} |z|^{\beta-d} C_{2}(\varepsilon,\varphi) \mathrm{d}z
\\
&= \int_{|z|>2\varepsilon} |z|^{\beta-d} \int_{|h|>\varepsilon} \sup_{y \in \mathbb{R}^{d}} |\varphi(z+h) a(z+h,y)| \frac{1}{|h|^{d+\beta}} \mathrm{d}h \mathrm{d}z
\\
& \leq \int_{|h| > \varepsilon} \int_{|z| > 2\varepsilon, |z+h|\leq \varepsilon} |z|^{\beta-d} \sup_{y \in \mathbb{R}^{d}} |\varphi(z+h) a(z+h,y)| \frac{1}{|h|^{d+\beta}}  \mathrm{d}z \mathrm{d}h 
\\
&\leq \sup_{z\in \mathbb{R}^{d}}|\varphi(z)| \sup_{z\in B_{\varepsilon},y\in \mathbb{R}^{d}}|a(z,y)| \int_{|h| > \varepsilon} \int_{|z| > 2\varepsilon, |z+h|\leq \varepsilon} |z|^{\beta-d} \frac{1}{|h|^{d+\beta}}    \mathrm{d}z  \mathrm{d}h
\\
&\leq \sup_{z\in \mathbb{R}^{d}}|\varphi(z)| \sup_{z\in B_{\varepsilon},y\in \mathbb{R}^{d}}|a(z,y)| \int_{|h| > \varepsilon} \varepsilon^{\beta-d} |B_{\varepsilon}(h)| \frac{1}{|h|^{d+\beta}}\mathrm{d}h
\\
&\leq C \sup_{z\in \mathbb{R}^{d}}|\varphi(z)| \sup_{z\in B_{\varepsilon},y\in \mathbb{R}^{d}}|a(z,y)|.
\end{align*}
Therefore, if we let 
\begin{align*}
\mathcal{C}_{2} = \mathcal{C}_{2}(\varepsilon,\varphi)= \sup_{z\in \mathbb{R}^{d}}|\varphi(z)| \sup_{z\in B_{\varepsilon},y\in \mathbb{R}^{d}}|a(z,y)|,
\end{align*}
it follows that
\begin{align*}
I_{2}
\leq  C \mathcal{C}_{2}(\varepsilon,\varphi)^{p} \|\phi(\Delta)u \|^{p}_{L_{p}(w)}.
\end{align*}
Hence, \eqref{eqn 07.16.19:54} follows.

\textbf{Case 2.} Now we consider the general case of $a$.
\\
Take a standard mollifier $\zeta \in C^{\infty}_{c}(\mathbb{R}^{d})$, and define 
$$
\zeta_{n}(z) = n^{d}\zeta(nz), \quad a_{n}(x,y) = [a(\cdot,y) \ast \zeta_{n}](x).
$$
Then $a_{n}(x,y) \to a(x,y)$ as $n\to \infty$, and $a_{n}$ satisfies Assumption \ref{24.08.18.13.17} without any changes to the function $\nu$. Therefore, we can use Fatou's lemma to obtain
\begin{align*}
\|\varphi \mathcal{L}^{a}_{V}u\|_{L_{p}(w)}  
\leq \liminf_{n\to\infty} \|\varphi \mathcal{L}^{a_{n}}_{V}u\|_{L_{p}(w)}
 \leq C \left( \mathcal{C}_{1}+\mathcal{C}_{2} \right)\|\phi(\Delta)u\|_{L_{p}(w)}.
\end{align*}
The lemma is proved.
\end{proof}

\subsection{Proof of Theorem \ref{thm 11.10.18:01}}

We now prove Theorem \ref{thm 11.10.18:01}.

\begin{proof}[Proof of Theorem \ref{thm 11.10.18:01}]
Let $\eta\in C_c^\infty(B_1(0))$ be a nonnegative function such that $\|\eta\|_{L_{p}}=1$.
For $z\in\mathbb R^d$, write
$$
\eta_z(x):=\eta(x-z).
$$
Since $\int_{\mathbb R^d}|\eta(x-z)|^p\,\mathrm dz
=
\|\eta\|_{L_p}^p$,
we have
\begin{align}
\label{eqn:translation-average}
\|\mathcal L_V^a u\|_{L_p(w)}^p
=
\int_{\mathbb R^d}\int_{\mathbb R^d}
|\eta_z(x)\mathcal L_V^a u(x)|^p
\,\mathrm dz\,w(x)\,\mathrm dx .
\end{align}
Observe that
\begin{equation}
\label{eqn:product-decomposition-main}
\begin{aligned}
\eta_z(x)\mathcal L_V^a u(x)
&=
\mathcal L_V^a(\eta_z u)(x)
-
u(x)\mathcal L_V^a\eta_z(x)\\
&\quad
-
\int_{\mathbb R^d}
(u(x+y)-u(x))
(\eta_z(x+y)-\eta_z(x))
a(x,y)V(y)\,\mathrm dy\\
&=:I_1(x,z)+I_2(x,z)+I_3(x,z).
\end{aligned}
\end{equation}
We estimate the three terms separately.

First, by Lemma \ref{lem 07.17.19:30}, we obtain (here $\gamma \in (1-(1/2\delta_{1}),1/2)$)
\begin{align}\label{eqn:I3-final}
\int_{\mathbb R^d}\int_{\mathbb R^d}
|I_3(x,z)|^p\,\mathrm dz\,w(x)\,\mathrm dx
\leq
C\|u\|_{H_p^{\phi,2\gamma}(w)}^p
\leq
C\|u\|_{H_p^{\phi,2}(w)}^p .
\end{align}
Next, we estimate $I_2$. 
Observe that $|I_2(x,z)|^p=|u(x)|^p|\mathcal{L}^{a}_{V}\eta_{z}(x)|^p$, and hence,
\begin{align}
\label{eqn:I2-reduction-final}
\int_{\mathbb R^d}\int_{\mathbb R^d}
|I_2(x,z)|^p\,\mathrm dz\,w(x)\,\mathrm dx
=
\int_{\mathbb R^d}
|u(x)|^p
\left(
\int_{\mathbb R^d}|\mathcal{L}^{a}_{V}\eta_{z}(x)|^p\,\mathrm dz
\right)
w(x)\,\mathrm dx .
\end{align}
We claim that
\begin{equation}
\label{eqn:I2-cutoff-bound}
\sup_{x\in\mathbb R^d}
\int_{\mathbb R^d}
|\mathcal{L}^{a}_{V}\eta_{z}(x)|^p\,\mathrm dz
\leq C.
\end{equation}
Indeed, since $a$ is bounded and $V\leq C_Vj_d^\phi$, it suffices to estimate
$$
\tilde{I}_{2}:=\int_{\mathbb R^d}
\left(
\int_{\mathbb R^d}
\left|
\eta_z(x+y)-\eta_z(x)
-
\nabla\eta_z(x)\cdot y\,\chi_{\delta_0,\delta_1}(y)
\right|
j_d^\phi(|y|)\,\mathrm dy
\right)^p
\mathrm dz .
$$
By the fundamental theorem of calculus and Minkowski's inequality, and Lemma \ref{lem:scaling-integrals}, it follows that
\begin{equation}\label{eqn 07.03.16:05}
\begin{aligned}
\tilde{I}_{2}
& \leq \int_{\mathbb{R}^{d}} \left(  \int_{|y| \geq 1} |\eta(x-z+y) - \eta(x-z)| j^{\phi}_{d}(|y|) \mathrm{d}y \right)^{p} \mathrm{d}z 
\\
& \quad + \mathbf{1}_{\delta_{0}>1/2} \int_{\mathbb{R}^{d}} \left( \int_{|y| \geq 1} \left| \int_{0}^{1} \nabla \eta(x-z+ry) \cdot y  \mathrm{d}r - \nabla \eta(x-z) \cdot y \right| j^{\phi}_{d}(|y|) \mathrm{d}y \right)^{p} \mathrm{d}z
\\
&\quad+ \int_{\mathbb{R}^{d}} \left(  \int_{|y| \leq 1} \int_{0}^{1} \int_{0}^{1}|y|^{2} |D^{2}\eta(x-z+sry)| \mathrm{d}r \mathrm{d}s j^{\phi}_{d}(|y|) \mathrm{d}y \right)^{p} \mathrm{d}z    
\\
& \quad + \mathbf{1}_{\delta_{1}< 1/2} \int_{\mathbb{R}^{d}} \left( \int_{|y|\leq 1} \left| \int_{0}^{1} \nabla \eta(x-z+ry) \cdot y \mathrm{d}r \right| j^{\phi}_{d}(|y|)\mathrm{d}y \right)^{p} \mathrm{d}z
\\
&\leq  2^{p}\|\eta\|^{p}_{L_{p}} \left(\int_{|y|\geq1} |y|^{-d}\phi(|y|^{-2}) \mathrm{d}y \right)^{p}  + \mathbf{1}_{\delta_{0}>1/2}   2^{p}\|\nabla \eta\|^{p}_{L_{p}} \left(  \int_{|y|\leq 1} |y|    j^{\phi}_{d}(|y|) \mathrm{d}y   \right)^{p}
\\
& \quad +\left(\int_{|y| \leq 1} |y|^{2} \int_{0}^{1}\int_{0}^{1}\left(  \int_{\mathbb{R}^{d}} |D^{2}\eta(x-z+sry)|^{p} \mathrm{d}z \right)^{1/p} \mathrm{d}r\mathrm{d}s j^{\phi}_{d}(|y|) \mathrm{d}y \right)^{p}  
\\
& \quad + \mathbf{1}_{\delta_{1}<1/2} \|\nabla \eta\|^{p}_{L_{p}} \left( \int_{|y|\leq 1} |y|  j^{\phi}_{d}(|y|) \mathrm{d}y \right)^{p}
\\
&\leq C \|\eta\|^{p}_{H^{2}_{p}}  \left( \int_{1}^{\infty} \rho^{-1}\phi(\rho^{-2})\mathrm{d}\rho  \mathrm{d}\rho+ \int_{0}^{1} \rho \phi(\rho^{-2}) \mathrm{d}\rho  \right) 
\\
& \quad + C \|\eta\|^{p}_{H^{2}_{p}} \left( \mathbf{1}_{\delta_{1}<1/2} \int_{0}^{1} \phi(\rho^{-2})\mathrm{d}\rho + \mathbf{1}_{\delta_{0}>1/2} \int_{1}^{\infty}\phi(\rho^{-2}) \mathrm{d}\rho \right)
\\
& \leq C \|\eta\|_{H^{2}_{p}}.
\end{aligned}
\end{equation}
This proves \eqref{eqn:I2-cutoff-bound}. 
Combining \eqref{eqn:I2-reduction-final} and \eqref{eqn:I2-cutoff-bound}, we get
\begin{align}
\label{eqn:I2-final}
\int_{\mathbb R^d}\int_{\mathbb R^d}
|I_2(x,z)|^p\,\mathrm dz\,w(x)\,\mathrm dx
\leq
C\|u\|_{L_p(w)}^p
\leq
C\|u\|_{H_p^{\phi,2}(w)}^p .
\end{align}

It remains to estimate $I_1$. 
Let $\varphi \in C_c^\infty(B_4(0))$ be such that $\varphi=1$ on $B_2(0)$, and put $\varphi_z(x):=\varphi(x-z)$.
We decompose $I_1=\varphi_z I_1+(1-\varphi_z)I_1$.
Applying Lemma \ref{lem 07.17.19:49}, we obtain
\begin{align}
\label{eqn:I1-local-final}
\int_{\mathbb R^d}\int_{\mathbb R^d}
|\varphi_z(x)I_1(x,z)|^p\,\mathrm dz\,w(x)\,\mathrm dx
\leq
C
\int_{\mathbb R^d}\int_{\mathbb R^d}
|\phi(\Delta)(\eta_z u)(x)|^p
w(x)\,\mathrm dx\,\mathrm dz .
\end{align}
Here we emphasize that the constant $C$ is independent of $z$ since the underlying one in Lemma \ref{lem 07.17.19:49} is independent of $z$. 
Note that $1-\varphi_z(x) \neq 0$ only when $|x-z|\geq2$. 
Since $\eta_z$ is supported in $B_1(z)$, it follows that
$$
(1-\varphi_{z})(x)\eta_z(x)=0,
\quad
(1-\varphi_{z})(x) \nabla\eta_z(x)=0,
$$
and also
$$
(1-\varphi_{z})(x)\eta_z(x+y)=0 \quad \forall\, |y|\leq1.
$$
Therefore, one can deduce that
$$
|(1-\varphi_z(x))I_1(x,z)|
\leq
C
\int_{|y|>1}
|\eta_z(x+y)u(x+y)|j_d^\phi(|y|)\,\mathrm dy .
$$
By Minkowski's inequality,
\begin{align*}
&\left(
\int_{\mathbb R^d}
|(1-\varphi_z(x))I_1(x,z)|^p\,\mathrm dz
\right)^{1/p}
\\
&\quad\leq
C
\int_{|y|>1}
|u(x+y)|j_d^\phi(|y|)
\left(
\int_{\mathbb R^d}
|\eta(x+y-z)|^p\,\mathrm dz
\right)^{1/p}
\mathrm dy
\\
&\quad =
C\|\eta\|_{L_p}
\int_{|y|>1}
|u(x+y)|j_d^\phi(|y|)\,\mathrm dy .
\end{align*}
Using \eqref{eqn 11.21.18:44}, we have
$$
\int_{|y|>1}
|u(x+y)|j_d^\phi(|y|)\,\mathrm dy
\leq
C\mathbb M u(x).
$$
Hence,
\begin{align}
\label{eqn:I1-far-final}
\int_{\mathbb R^d}\int_{\mathbb R^d}
|(1-\varphi_z(x))I_1(x,z)|^p\,\mathrm dz\,w(x)\,\mathrm dx
\leq
C
\|\mathbb M u\|_{L_p(w)}^p
\leq
C\|u\|_{L_p(w)}^p .
\end{align}
Combining \eqref{eqn:I1-local-final} and \eqref{eqn:I1-far-final}, we get
\begin{align}
\label{eqn:I1-final-reduction}
\int_{\mathbb R^d}\int_{\mathbb R^d}
|I_1(x,z)|^p\,\mathrm dz\,w(x)\,\mathrm dx
\leq
C
\int_{\mathbb R^d}\int_{\mathbb R^d}
|\phi(\Delta)(\eta_z u)(x)|^p
w(x)\,\mathrm dx\,\mathrm dz
+
C\|u\|_{L_p(w)}^p .
\end{align}
It only remains to bound the double integral in \eqref{eqn:I1-final-reduction}. 
Like \eqref{eqn:product-decomposition-main}, we write
\begin{equation*}
\begin{aligned}
\phi(\Delta)(\eta_z u)(x)
&=
\eta_z(x)\phi(\Delta)u(x)
+
u(x)\phi(\Delta)\eta_z(x)\\
&\quad
+
\int_{\mathbb R^d}
(u(x+y)-u(x))
(\eta_z(x+y)-\eta_z(x))
j_d^\phi(|y|)\,\mathrm dy\\
&=:
\mathcal I_0(x,z)+\mathcal I_2(x,z)+\mathcal I_3(x,z).
\end{aligned}
\end{equation*}
It is easy to see that
$$
\int_{\mathbb R^d}\int_{\mathbb R^d}
|\mathcal I_0(x,z)|^p\,\mathrm dz\,w(x)\,\mathrm dx
\leq
C
\|\phi(\Delta)u\|_{L_p(w)}^p .
$$
For $\mathcal{I}_{2}$, we have the following correspondence to  \eqref{eqn:I2-final}
$$
\int_{\mathbb R^d}\int_{\mathbb R^d}
|\mathcal I_2(x,z)|^p\,\mathrm dz\,w(x)\,\mathrm dx
\leq
C
\|u\|_{L_p(w)}^p .
$$
Lemma \ref{lem 07.17.19:30} gives
$$
\int_{\mathbb R^d}\int_{\mathbb R^d}
|\mathcal I_3(x,z)|^p\,\mathrm dz\,w(x)\,\mathrm dx
\leq
C
\|u\|_{H_p^{\phi,2\gamma}(w)}^p
\leq
C
\|u\|_{H_p^{\phi,2}(w)}^p .
$$
Therefore,
\begin{align}
\label{eqn:localized-phi-final}
\int_{\mathbb R^d}\int_{\mathbb R^d}
|\phi(\Delta)(\eta_z u)(x)|^p
w(x)\,\mathrm dx\,\mathrm dz
\leq
C
\|u\|_{H_p^{\phi,2}(w)}^p .
\end{align}
Finally, combining
\eqref{eqn:translation-average},
\eqref{eqn:product-decomposition-main},
\eqref{eqn:I3-final},
\eqref{eqn:I2-final},
\eqref{eqn:I1-final-reduction}, and
\eqref{eqn:localized-phi-final}, we obtain $
\|\mathcal L_V^a u\|_{L_p(w)}^p\leq C \|u\|_{H_p^{\phi,2}(w)}^p$.
The theorem is proved.
\end{proof}

\section{Proof of weighted Sobolev regularity estimates: Theorem \ref{thm 11.24.15:58}}

\begin{lem}\label{lem 11.21.17:27}\label{thm 11.13.15:42}
Let  $\phi$ be a Bernstein function satisfying Assumption \ref{asm 07.21.13:06}($\delta_0,\delta_1$) with $0<\delta_0\leq\delta_1<1$, $1<p,q<\infty$, $w_{1} \in A_{p}(\mathbb{R}^{d})$, $w_{2} \in A_{q}(\mathbb{R})$, and $0<T<\infty$.
Suppose that $V(y)=j_d^{\phi}(|y|)$, and the function $a=a(t,y)$ satisfies Assumption \ref{asm 11.10.16:43}.  
Then for any  $f\in L_{q,p}(w_{2},w_{1},T)$, there exists a unique solution $u \in \mathbb{H}^{\phi,2}_{q,p,0}(w_{2},w_{1},T)$ to the equation
\begin{align}\label{eqn 11.10.16:36}
\partial_{t}u = \mathcal{L}_{V}^{a}u +f, \quad u(0,\cdot)=0.
\end{align}
Moreover, the solution $u$ satisfies the following estimates
\begin{align}\label{eqn 08.19.16:37}
\|u\|_{\mathbb{H}^{\phi,2}_{q,p,0}(w_{2},w_{1},T)} \leq C_{0} \|f\|_{L_{q,p}(w_{2},w_{1},T)},
\end{align}
where the constant $C_{0}=C_{0}(c_{0},c_{1},\delta_{0},\delta_{1},a,d,p,q,[w_{1}]_{A_{p}(\mathbb{R}^{d})},[w_{2}]_{A_{q}(\mathbb{R})},T)$ independent of $f$, and
\begin{align}\label{eqn 08.19.16:37-2}
\| \mathcal{L}^{a}_{V}u\|_{L_{q,p}(w_{2},w_{1},T)} + \| \phi(\Delta)u\|_{L_{q,p}(w_{2},w_{1},T)} \leq C_{1} \|f\|_{L_{q,p}(w_{2},w_{1},T)},
\end{align}
where the constant $C_{1}=C_{1}(c_{0},c_{1},\delta_{0},\delta_{1},a,d,p,q,[w_{1}]_{A_{p}(\mathbb{R}^{d})},[w_{2}]_{A_{q}(\mathbb{R})})$ independent of $f$ and $T$.
\end{lem}

\begin{proof}
If $a(t,y) \equiv 1$, then $\mathcal{L}^{a}_{V}=\phi(\Delta)$, and the result is a direct consequence of Theorem \ref{thm 11.13.15:42}. Hence, we can assume that $a$ is not a constant. We invoke a theory of additive processes \cite{CZ23additive,CKP24anisotropic}. 
\\
\textbf{Step 1.} Uniqueness of solutions.
\\
We first prove the uniqueness of solutions. Suppose that $u,v \in \mathbb{H}^{\phi,2}_{q,p}(w_{2},w_{1},T)$ are solutions to \eqref{eqn 11.10.16:36}. Then $w=u-v$ satisfies
\begin{align*}
\partial_{t}w = \mathcal{L}^{a}_{V}w, \quad w(0,\cdot)=0.
\end{align*}
Let $w_{n}$ be a defining sequence of $w$ in $\mathbb{H}^{\phi,2}_{q,p}(w_{2},w_{1},T)$ and let $f_{n} = \partial_{t}w_{n}  - \mathcal{L}^{a}_{V}w_{n}$. Then applying Theorem \ref{thm 11.10.18:01}, we have
\begin{align}\label{eqn 11.10.19:00}
\|f_{n}\|_{L_{q,p}(w_{2},w_{1},T)} &= \|\partial_{t}w_{n} - \mathcal{L}^{a}_{V}w_{n}\|_{L_{q,p}(w_{2},w_{1},T)}   \nonumber
\\
&\leq \|(\partial_{t}-\mathcal{L}^{a}_{V})(w_{n}-w)\|_{L_{q,p}(w_{2},w_{1},T)} + \|(\partial_{t} - \mathcal{L}^{a}_{V})w\|_{L_{q,p}(w_{2},w_{1},T)}    \nonumber
\\
&\leq C \|w_{n} - w\|_{\mathbb{H}^{\phi,2}_{q,p,0}(w_{2},w_{1},T)}.
\end{align}
Moreover, due to the probabilistic representation
\begin{align}\label{eqn 08.13.16:23}
w_{n}(t,x) = \int_{0}^{t} \mathbb{E}[f_{n}(s,x+Z_{s,t})] \mathrm{d}s,
\end{align}
where $Z_{s,t}$ is a additive process with the characteristic function (see \cite[Lemma 2.12]{CKP24anisotropic})
\begin{align*}
\mathbb{E}[\mathrm{e}^{i \xi Z_{s,t}}] = \exp{\left( \int_{s}^{t} \int_{\mathbb{R}^{d}} \left( \mathrm{e}^{i \xi \cdot y} -1 -i \xi \cdot y \chi_{\delta_0,\delta_1}(y) \right) a(r,y) j^{\phi}_{d}(y) \mathrm{d}y \mathrm{d}r \right)}.
\end{align*}
\textit{Claim 1}: We claim that the transition density $q^{a}_{s,t}(x)$ of $Z_{s,t}$ exists and satisfies
\begin{align*}
c^{-1}\varrho_{\phi}(t-s,|x|) \leq q^{a}_{s,t}(x) \leq c \varrho_{\phi}(t-s,|x|) \quad \forall\, 0<t-s<\infty, \quad \forall\, x,y \in \mathbb{R}^{d},
\end{align*}
where
\begin{align*}
\varrho_\phi(\tau,r)
:=
\left( \phi^{-1}(\tau^{-1}) \right)^{d/2}
\wedge
\frac{\tau\phi(r^{-2})}{r^d}.
\end{align*}
Define 
\begin{align*}
\Phi(r):=\frac{1}{\phi(r^{-2})},
\qquad r>0.
\end{align*}
Then by Assumption \ref{asm 07.21.13:06}, it follows that
\begin{align}\label{eq:Phi-global-scaling}
c_0
\left(\frac{R}{r}\right)^{2\delta_0}
\leq
\frac{\Phi(R)}{\Phi(r)}
\leq
c_1
\left(\frac{R}{r}\right)^{2\delta_1}
\quad
\forall\, 0<r\leq R<\infty.
\end{align}
Using  \eqref{eqn 02.25.17:23}, we can rewrite
\begin{align*}
a(t,y)j_d^\phi(|y|)
=
\frac{a(t,y)|y|^d\Phi(|y|)}{|y|^{d}\Phi(|y|)} j^{\phi}_{d}(|y|) := \frac{m(t,y)}{|y|^{d}\Phi(|y|)},
\end{align*}
where $m(t,y)$ also satisfies Assumption \ref{asm 11.10.16:43} (with the same $\delta_{0},\delta_{1}$).  If
$\delta_0\leq 1/2\leq\delta_1$, from Assumption \ref{asm 11.10.16:43}
\begin{align*}
0
=
\int_{r<|y|<R}
y\,a(t,y)j_d^\phi(|y|)\,\mathrm{d}y
=
\int_r^R
\rho^{d} j_d^\phi(\rho)
\left(
\int_{\mathbb S^{d-1}}
\sigma \, a(t,\rho \sigma)\,\mathrm{d}\sigma
\right)
\mathrm{d}\rho \quad \forall\, 0<r<R<\infty.
\end{align*}
Hence,
\begin{align*}
\int_{\mathbb S^{d-1}}
\sigma \, a(t,\rho\sigma)\,\mathrm{d}\theta=0
\end{align*}
for almost every $\rho>0$. Thus, $m$ also satisfies
\begin{align*}
\int_{|y| \leq r} y m(t,y) \mathrm{d}y
 = \int_{0}^{r} \rho^{d} j_{d}(\rho) \rho^{d} \Phi(\rho) \left( \int_{\mathbb{S}^{d-1}} \sigma \, a(t, \rho \sigma) \mathrm{d}\sigma \right) \mathrm{d}\rho = 0.
\end{align*}
for any $r>0$.  Furthermore, one can observe that
the functions
\begin{align*}
\Phi_\lambda(r)
:=
\frac{\Phi(\Phi^{-1}(\lambda)r)}{\lambda}
\quad (\lambda>0),
\end{align*}
satisfy \eqref{eq:Phi-global-scaling} with constants
independent of any $\lambda>0$. From this, one can check that \cite[Proposition 2.1]{CZ23additive} holds without the restriction $0<t-s<1$, by following the argument therein (put $m$ in place of $\kappa$). Combining this with
\cite[Theorem 2.2]{CZ23additive}, we deduce that the transition density $q^{a}_{s,t}$ exists and satisfies
\begin{align*}
q^{a}_{s,t}(x,y) &\asymp \left( \Phi^{-1}(t-s) \right)^{-d}\left( 1 \wedge \frac{(t-s)(\Phi^{-1}(t-s))^{d}}{|x|^{d}\Phi(|x|)} \right)
\\
& =  (\phi^{-1}((t-s)^{-1}))^{d/2} \left( 1 \wedge (\phi^{-1}((t-s)^{-1}))^{-d/2}\frac{(t-s)\phi(|x-y|^{-2})}{|x-y|^{d}} \right)
\\
& = \varrho_{\phi}(t-s,|x-y|).
\end{align*}
\textit{Claim 2}: We argue that
\begin{align}\label{eq:additive-kernel-L1}
\sup_{\tau>0}
\int_{\mathbb R^d}
\varrho_\phi(\tau,|x|)\,\mathrm{d}x
\leq C.
\end{align}
Observe that
\begin{align}\label{eqn 08.12.16:35}
\int_{|x|\leq (\phi^{-1}(\tau^{-1}))^{-1/2}} \varrho_{\phi}(\tau,x) \mathrm{d}x = \int_{|x|\leq (\phi^{-1}(\tau^{-1}))^{-1/2}} \left( \phi^{-1}(\tau^{-1}) \right)^{d/2} \mathrm{d}x = C.
\end{align}
Using \eqref{eqn 07.01.15:42}, we obtain
\begin{align}\label{eqn 08.12.16:35-2}
 \int_{|x|>(\phi^{-1}(\tau^{-1}))^{-1/2}}
\varrho_\phi(\tau,x)\,\mathrm{d}x
=
C
\int_{(\phi^{-1}(\tau^{-1}))^{-1/2}}^\infty
\tau \phi(r^{-2}) r^{-1} \,\mathrm{d}r
\leq
C.
\end{align}
From \eqref{eqn 08.12.16:35} and \eqref{eqn 08.12.16:35-2}, we obtain \eqref{eq:additive-kernel-L1}.

\textit{Claim 1}, and \textit{Claim 2} imply that (recall \eqref{eqn 08.13.16:23})
\begin{align*}
w_{n}(t,x) = \int_{0}^{t} \left( \int_{\mathbb{R}^{d}} q^{a}_{s,t}(y-x) f_{n}(s,y) \mathrm{d}y \right) \mathrm{d}s,
\end{align*}
where the kernel $q^{a}_{s,t}$ has a radially decreasing majorant which is integrable in $x$ uniformly in $s,t$. Therefore, we apply \textit{e.g.} \cite[Corollary 2.1.12]{Grafakos} to derive 
\begin{align}\label{eqn 12.18.14:18}
\left|  w_{n}(t,x) \right|  \leq C \int_{0}^{t} \mathbb{M}f_{n}(s,x) \mathrm{d}s.
\end{align}
Applying Minkowski's inequality to \eqref{eqn 12.18.14:18}, we have
\begin{align*}
\left( \int_{\mathbb{R}^{d}} \left| w_{n}(t,x) \right|^{p} w_{1}(x) \mathrm{d}x \right)^{1/p} &\leq C(T) \int_{0}^{t} \left( \int_{\mathbb{R}^{d}} \left| \mathbb{M}f_{n}(s,x) \right|^{p} w_{1}(x) \mathrm{d}x \right)^{1/p} \mathrm{d}s   \nonumber
\\
&\leq C(T) \int_{0}^{t} \left( \int_{\mathbb{R}^{d}} \left| f_{n}(s,x) \right|^{p} w_{1}(x) \mathrm{d}x \right)^{1/p} \mathrm{d}s   \nonumber
\\
&:= C(T) \int_{0}^{t} F_{n}(s) \mathrm{d}s   \nonumber
\\
&\leq C(T) \mathcal{M}F_{n}(t),
\end{align*}
where $\mathcal{M}$ is the Hardy-Littlewood maximal function in $t$. This yields that
\begin{align*}
\|w_{n}\|^{q}_{L_{q,p}(w_{2},w_{1},T)} &\leq C(T) \int_{0}^{T} \left| \mathcal{M}F_{n}(t)  \right|^{q} w_{2}(t)\mathrm{d}t   \nonumber
\\
& \leq C(T) \int_{0}^{T} \left| F_{n}(t) \right|^{q} w_{2}(t) \mathrm{d}t   \nonumber
\\
& = C(T) \|f_{n}\|^{q}_{L_{q,p}(w_{2},w_{1},T)}.
\end{align*}
Combining \eqref{eqn 11.10.19:00}, we have
\begin{align*}
\|w_{n}\|^{q}_{L_{q,p}(w_{2},w_{1},T)} \leq C(T) \|w_{n}-w\|^{q}_{\mathbb{H}^{\phi,2}_{q,p,0}(w_{2},w_{1},T)}.
\end{align*}
Since the right hand side of the above inequality goes to 0 as $n\to \infty$ (recall that $w_{n}$ is a defining sequence of $w$), $w_{n}$ converges to $0$ in $L_{q,p}(w_{2},w_{1},T)$, which implies that $u-v=w=0$ in $\mathbb{H}^{\phi,2}_{q,p,0}(w_{2},w_{1},T)$. Hence, we prove the uniqueness of solutions. 
\\
\textbf{Step 2.} Existence and estimations of solutions.
\\
For simplicity, we denote
\begin{align*}
\bar a:=\frac{a_0}{2},
\quad
b(t,y):=a(t,y)-\bar{a} \in [\bar{a} , a_{1}-\bar{a}].
\end{align*}
Also, let $p^{\bar{a}}_{t-s}(x)$ be the transition density of SBM generating the operator $(\bar{a}\phi)(\Delta)$, and $q^{b}_{s,t}(x)$ be the transition density whose existence is guaranteed in \textit{Claim 1} of \textbf{Step 1}. Let $f_{n} \in C^{\infty}_{c}(\mathbb{R}^{d}_{+})$ be a sequence of functions which converges to $f$ in $L_{q,p}(w_{2},w_{1},T)$. By following the argument in \cite[Theorem 2.8]{CKP24anisotropic}, 
we further observe that
\begin{align}\label{eqn 08.20.17:25}
u_{n}(t,x)  & = \int_{0}^{t} \mathbb{E}[f(s,x+Z_{s,t})] \mathrm{d}s  
\\
&= \int_{0}^{t} \int_{\mathbb{R}^{d}} p^{\bar{a}}_{t-s}(x-y) \left( \int_{\mathbb{R}^{d}} q^{b}_{s,t}(y-z) f_{n}(s,z) \mathrm{d}z \right) \mathrm{d}y \mathrm{d}s  \nonumber
\end{align}
is a unique classical solution to 
\begin{align*}
\partial_{t}u_{n}(t,x) = \mathcal{L}^{a}_{V}u_{n}(t,x) + f_{n}(t,x) \quad u_{n}(0,x)=0.
\end{align*}

To prove the estimate 
\begin{align}\label{eqn 08.13.17:56}
\| \phi(\Delta) u_{n}\|_{L_{q,p}(w_{2},w_{1},T)} \leq C \|f_{n}\|_{L_{q,p}(w_{2},w_{1},T)},
\end{align}
We use the standard Littlewood-Paley projection. Take $\Psi \in \mathcal{S}(\mathbb R^d)$, whose Fourier transform $\hat{\Psi}$ satisfies
\begin{align*}
\operatorname{supp}\hat{\Psi}
\subset
\{
\xi\in\mathbb{R}^d:
1/2 \leq |\xi| \leq 2
\},
\quad
\sum_{j\in\mathbb Z}\hat{\Psi}_{j}(\xi)=:\sum_{j\in\mathbb Z}\hat{\Psi}(2^{-j}\xi)=1
\quad
\text{for} \quad \xi\neq0.
\end{align*}
Let $\Delta_j$ denote the corresponding 
Littlewood-Paley operator
\begin{align*}
\Delta_{j}f = \mathcal{F}^{-1}[\hat{\Psi}_{j}\mathcal{F}(f)] \quad j \in \mathbb{Z}.
\end{align*}
We also choose
$\Phi \in \mathcal{S}(\mathbb{R}^{d})$
such that the Fourier transform $\hat{\Phi}$ is supported on a annulus and satisfies $\hat{\Phi}=1$ on the support of $\hat{\Psi}$, and denote the corresponding operator by $\widetilde\Delta_j$.
For each $j \in \mathbb{Z}$, set
\begin{align*}
\lambda_j:=\phi(2^{2j}).
\end{align*}
For $r>0$, consider the multiplier
\begin{align*}
m_{j,r}(\xi)
:=
- \hat{\Psi}_{j}(\xi)
\phi(|\xi|^2)
\mathrm e^{-\bar a \phi(|\xi|^2)r}.
\end{align*}
The weak scaling condition gives
\begin{align*}
 \phi(|\xi|^2) \asymp  \lambda_j \quad 2^{j-1} \leq |\xi| \leq 2^{j+1}.
\end{align*}
Using this and \eqref{eqn 01.28.13:23}, for every multi-index $\alpha$, it follows that
\begin{align*}
\left|
D_\xi^\alpha\phi(|\xi|^2)
\right|
\leq
C_\alpha
2^{-j|\alpha|}
\lambda_j \quad \quad 2^{j-1} \leq |\xi| \leq 2^{j+1}.
\end{align*}
Therefore, by the chain rule,
\begin{align*}
|D_\xi^\alpha m_{j,r}(\xi)|
\leq
C_\alpha
2^{-j|\alpha|}
\lambda_j
\mathrm e^{-c r\lambda_j}\quad 2^{j-1} \leq |\xi| \leq 2^{j+1}
\end{align*}
holds, where the constant $c$ does not depend on $j$. Therefore, we have
\begin{align*}
& \left| \int_{\mathbb{R}^{d}} \mathrm{e}^{i \xi \cdot x}  (1+\mathbf{1}_{|2^{j}x| >1}|x|^{-2N}|x|^{2N}) m_{j,r}(\xi) \mathrm{d}\xi \right|
\\
& \leq  \int_{2^{j} \leq |\xi| \leq 2^{j+1}} \left| m_{r,j}(\xi) \right| \mathrm{d}\xi + C\mathbf{1}_{|2^{j}x| > 1} |x|^{-2N} \int_{2^{j} \leq |\xi| \leq 2^{j+1}} \left| \Delta_{\xi}^{N}m_{j,r}(\xi) \right| \mathrm{d}\xi
\\
&\leq C \lambda_{j}\mathrm{e}^{-cr\lambda_{j}} 2^{jd} + C \mathbf{1}_{|2^{j}x| > 1}|x|^{-2N} \lambda_{j}\mathrm{e}^{-cr\lambda_{j}} 2^{jd}2^{-2Nj}
\\
&=
C
\lambda_j
\mathrm e^{-cr\lambda_j}
2^{jd}
(1+\mathbf{1}_{|2^{j}x|>1}|2^{j}x|^{-2N}) \leq C_{N} \lambda_{j} \mathrm{e}^{-c r \lambda_{j}}2^{jd}(1+|2^{j}x|^{2})^{-N}
\end{align*}
for any $N>d$. 
This implies that the operator
\begin{gather}\label{eqn 08.13.20:11}
f \mapsto \Delta_{j}\phi(\Delta)\mathcal{P}^{\bar{a}}_{r}f = \Delta_{j}\phi(\Delta)[p^{\bar{a}}_{r}\ast f] = \Delta_{j}\phi(\Delta)[p^{\bar{a}}_{r}\ast \widetilde{\Delta}_{j}f]
\\
\left( \text{i.e. } \mathcal{F}_{d}[\Delta_{j}\phi(\Delta)\mathcal{P}^{\bar{a}}_{r}f](\xi) = m_{j,r}(\xi)\mathcal{F}_{d}[f](\xi)   \right)  \nonumber
\end{gather}
has a kernel $K_{j,r}$($=\mathcal{F}_{d}^{-1}[m_{j,r}]$), which has a radially decreasing majorant integrable in $x$ uniformly in $r>0$. Note that the last equality in \eqref{eqn 08.13.20:11} holds since $\hat{\Phi}=1$ on the support of $\hat{\Psi}$. Using \cite[Corollary 2.1.12]{Grafakos} again, we have 
\begin{align}\label{eq:dyadic-P-estimate}
\left|
\Delta_j\phi(\Delta)\mathcal P^{\bar{a}}_{r} f_{n}(x)
\right|
\leq
C
\lambda_j
\mathrm e^{-c r \lambda_j}
\mathbb M
\bigl(\widetilde\Delta_j f_{n} \bigr)(x).
\end{align}
Combining \eqref{eq:dyadic-P-estimate}, and using the argument in \textit{Claim 1} and \textit{Claim 2} in \textbf{Step 1} for $q^{b}_{s,t}$ it follows that
\begin{align}\label{eq:dyadic-solution-time-convolution}
\left|
\Delta_j\phi(\Delta)u_{n}(t,x)
\right|
&\leq
C
\int_0^t
\lambda_j
\mathrm e^{-c (t-s)\lambda_j}
\mathbb M_x^2
\bigl(\widetilde\Delta_jf_{n}(s,\cdot)\bigr)(x)
\,\mathrm{d}s.
\end{align}
Note that
\begin{align*}
\sup_{j\in\mathbb Z}
\int_0^\infty \lambda_j
\mathrm e^{-c\lambda_jr}
\mathbf 1_{\{r>0\}} \mathrm{d}r
<\infty.
\end{align*}
Therefore, using \cite[Corollary 2.1.12]{Grafakos} to the right hand side of  \eqref{eq:dyadic-solution-time-convolution}, we have
\begin{align}\label{eq:dyadic-space-time-maximal}
\left|
\Delta_j\phi(\Delta)u_{n}(t,x)
\right|
\leq
C
\mathcal{M}
\mathbb{M}^2
\bigl(\widetilde\Delta_j f_{n} \bigr)(t,x).
\end{align}
Using the Littlewood-Paley inequality, \eqref{eq:dyadic-space-time-maximal}, it follows that
\begin{align*}
\|\phi(\Delta)u_{n}\|_{L_{p,p}(w_{2},w_{1},T)}
& \leq \|\phi(\Delta)u_{n}\|_{L_p(\mathbb R^{d+1};w_{2}(t)w_1(x)\,\mathrm{d}t\mathrm{d}x)}
\nonumber
\\
&\leq
C
\left\|
\left(
\sum_{j\in\mathbb Z}
\left|\Delta_{j}\phi(\Delta)u_{n} \right|^{2}
\right)^{1/2}
\right\|_{L_p(\mathbb R^{d+1};w_{2}(t)w_1(x)\,\mathrm{d}t\mathrm{d}x)}
\nonumber\\
&\leq
C
\left\|
\left(
\sum_{j\in\mathbb Z}
\left|
\mathcal{M}\mathbb M_x^2
\widetilde\Delta_jf_{n}
\right|^2
\right)^{1/2}
\right\|_{L_p(\mathbb R^{d+1};w_{2}(t)w_1(x)\,\mathrm{d}t\mathrm{d}x)}
\nonumber\\
&\leq
C
\left\|
\left(
\sum_{j\in\mathbb Z}
|\widetilde\Delta_jf_{n}|^2
\right)^{1/2}
\right\|_{L_p(\mathbb R^{d+1};w_{2}(t)w_1(x)\,\mathrm{d}t\mathrm{d}x)}
\nonumber\\
&\leq C \|f_{n}\|_{L_p(\mathbb R^{d+1};w_{2}(t)w_1(x)\,\mathrm{d}t\mathrm{d}x)}
=
C
\|f_{n}\|_{L_{p,p}(w_{2},w_1,T)}
\end{align*}
for any $1<p<\infty$, and $w_{1} \in A_{p}(\mathbb{R}^{d})$ and $w_{2} \in A_{p}(\mathbb{R})$. Also note that, for the fourth inequality above, we used the fact that the Hardy-Littlewood maximal functions are bounded on $L_{r}(\ell_{2} ; w\mathrm{d}x)$ for any $1<r<\infty$ and $w \in A_{r}(\mathbb{R}^{d})$. In other words, 
\begin{align*}
\int_{0}^{T} \|\phi(\Delta)u_{n}(t,\cdot)\|^{q}_{L_{p}(w_{1})} w_{2}(t)\mathrm{d}t \leq C \int_{0}^{T} \|f_{n}(t,\cdot)\|^{q}_{L_{p}(w_{1})}w_{2}(t)\mathrm{d}t
\end{align*}
holds for any $1<p<\infty$, $w_{1} \in A_{p}(\mathbb{R}^{d})$, and $q=q_{0}=p$ and $w_{2}\in A_{q_{0}}(\mathbb{R})$.  Hence, applying the extrapolation theorem \cite[Theorem 1.4]{extrapolation} (also see  \cite[Theorem 2.5]{DK18} for detail) we have \eqref{eqn 08.13.17:56}. 

Moreover, by mimicking the argument in \textbf{Step 1} for $w_{n}$ (recall \eqref{eqn 08.20.17:25}), we also obtain
\begin{align}\label{eqn 08.19.16:31}
\|u_{n}\|_{L_{q,p}(w_{2},w_{1},T)} \leq C(T) \|f_{n}\|_{L_{q,p}(w_{2},w_{1},T)}.
\end{align}
Combining Theorem \ref{thm 01.15.15:07}, \eqref{eqn 08.13.17:56} and \eqref{eqn 08.19.16:31}, it follows that
\begin{align*}
\|\mathcal{L}^{a}_{V}u_{n}\|_{L_{q,p}(w_{2},w_{1},T)} + \|\phi(\Delta)u_{n}\|_{L_{q,p}(w_{2},w_{1},T)} \leq C \|\phi(\Delta)u_{n}\|_{L_{q,p}(w_{2},w_{1},T)} \leq C\|f_{n}\|_{L_{q,p}(w_{2},w_{1},T)},
\end{align*}
and
\begin{align*}
\|u_{n}\|_{\mathbb{H}^{\phi,2}_{q,p,0}(w_{2},w_{1},T)} \leq C_{0} \|f_{n}\|_{L_{q,p}(w_{2},w_{1},T)}.
\end{align*}
Hence, $u_{n}$ is a Cauchy sequence in $\mathbb{H}^{\phi,2}_{q,p,0}(w_{2},w_{1},T)$. Therefore, by taking $u$ as the limit of $u_{n}$ in $\mathbb{H}^{\phi,2}_{q,p,0}(w_{2},w_{1},T)$, we deduce that $u$ is a solution to \eqref{eqn 11.10.16:36} satisfying estimates \eqref{eqn 08.19.16:37} and \eqref{eqn 08.19.16:37-2}. The lemma is proved.
\end{proof}

\begin{lem}\label{lem 11.24.16:13}
Let $\phi$ be a Bernstein function satisfying Assumption \ref{asm 07.21.13:06}($\delta_0,\delta_1$) with $0<\delta_0\leq\delta_1<1$, $1<p,q<\infty$, $w_{1} \in A_{p}(\mathbb{R}^{d})$, $w_{2} \in A_{q}(\mathbb{R})$, and $0<T<\infty$.
Suppose that $V(y)=j_d^{\phi}(|y|)$, and the function $a=a(t,x,y)$ satisfies Assumption \ref{asm 11.10.16:43}, and $f\in L_{q,p}(w_{2},w_{1},T)$.  Then there exists $\varepsilon_{0}$ such that
for any $\varepsilon \leq \varepsilon_{0}$, $z\in \mathbb{R}^{d}$, and for any $u\in \mathbb{H}^{\phi,2}_{q,p,0}(w_{2},w_{1},T)$, whose support in spatial variable is contained in $B_{\varepsilon}(z)$ satisfying
\begin{align*}
\partial_{t}u = \mathcal{L}_{V}^{a}u +f,  \quad u(0,\cdot)=0 \quad (x\in \mathbb{R}^{d}),
\end{align*}
the following estimation
\begin{align*}
\|\phi(\Delta)u\|_{L_{q,p}(w_{2},w_{1},t)} \leq C \left( \|f\|_{L_{q,p}(w_{2},w_{1},t)} + \phi(\varepsilon^{-2}) \|u\|_{L_{q,p}(w_{2},w_{1},t)} \right) \quad \forall\, 0<t\leq T
\end{align*}
holds, where the constant $C=C(c_{0},c_{1},\delta_{0},\delta_{1},a,d,p,q,[w_{1}]_{A_{p}(\mathbb{R}^{d})},[w_{2}]_{A_{q}(\mathbb{R})})$ is independent of $t,z,\varepsilon$.
\end{lem}
\begin{proof}
Let $\varepsilon>0$ and let $z\in \mathbb{R}^{d}$. Take a smooth function $\varphi \in C^{\infty}_{c}(B_{4\varepsilon}(z))$ such that $0\leq \varphi \leq 1$, and $\varphi=1$ on $B_{2\varepsilon}(z)$. We can further assume that $u(t,x)$ is sufficiently smooth in $x$. Denote $u_{z}(t,x)=u(t,x+z)$, $\varphi_{z}(x)=\varphi(x+z)$, $f_{z}(t,x) = f(t,x+z)$, and $a_{z}(t,x,y) = a(t,x+z,y)$, and define the operators
\begin{align*}
\mathcal{L}^{z}_{V}u(x) = \int_{\mathbb{R}^{d}} \left(  u(x+y) - u(x) -y \cdot \nabla u(x) \chi_{\delta_{0},\delta_{1}}  \right) a_{z}(t,x,y) j^{\phi}_{d}(|y|)\mathrm{d}y,
\\
\mathcal{L}^{0}_{V}u(x) = \int_{\mathbb{R}^{d}} \left(  u(x+y) - u(x) -y \cdot \nabla u(x) \chi_{\delta_{0},\delta_{1}}  \right) a(t,z,y) j^{\phi}_{d}(|y|)\mathrm{d}y.
\end{align*}
Then $u_{z} \in \mathbb{H}^{\phi,2}_{q,p,0}(w_{2},w_{1},T)$ satisfies
\begin{align*}
\partial_{t}u_{z}(t,x) &= \mathcal{L}^{z}_{V}u_{z}(t,x) + f_{z}(t,x)
\\
& = \mathcal{L}^{0}_{V}u_{z}(t,x) + \varphi_{z}(\mathcal{L}^{z}_{V}-\mathcal{L}^{0}_{V})u_{z}(t,x) + (1-\varphi_{z})(\mathcal{L}^{z}_{V}-\mathcal{L}^{0}_{V})u_{z}(t,x) +f_{z}(t,x).
\end{align*}
Hence, for any $0<t\leq T$, Lemma \ref{lem 11.21.17:27} with $w_{1,z}(x)=w_{1}(x+z) \in A_{p}(\mathbb{R}^{d})$ yields the following estimation
\begin{align*}
\|\phi(\Delta)u_{z}\|_{L_{q,p}(w_{2},w_{1,z},t)} &\leq C \|\varphi_{z}(\mathcal{L}^{z}_{V}-\mathcal{L}^{0}_{V})u\|_{L_{q,p}(w_{2},w_{1,z},t)} 
\\
&\quad + C \|(1-\varphi_{z})(\mathcal{L}^{z}_{V}-\mathcal{L}^{0}_{V})u\|_{L_{q,p}(w_{2},w_{1,z},t)} +C \|f_{z}\|_{L_{q,p}(w_{2},w_{1,z},t)},
\end{align*}
where the constant $C$ is independent of $z, \varepsilon,t$. Applying Lemma \ref{lem 07.17.19:49} to $\varphi_{z}(\mathcal{L}^{z}_{V}-\mathcal{L}^{0}_{V})u$, we have
\begin{align}\label{eqn 11.21.18:50}
\int_{\mathbb{R}^{d}}\left| \varphi_{z}(\mathcal{L}^{z}_{V}-\mathcal{L}^{0}_{V})u_{z}(t,x) \right|^{p} w_{1,z}(x) \mathrm{d}x
\leq C \mathcal{C}(\varepsilon,\varphi,\nu) \int_{\mathbb{R}^{d}} \left| \phi(\Delta) u_{z}(t,x) \right|^{p} w_{1,z}(x) \mathrm{d}x,
\end{align}
where the constant $C$ is independent of $z,\varepsilon,t$ and
\begin{align*}
\mathcal{C}(\varepsilon,\varphi,\nu)
= \sup_{z\in\mathbb{R}^{d}}\left( |\varphi(z)| + \varepsilon |D\varphi(z)| \right) \left( (1+\varepsilon)\nu(8\varepsilon) + \varepsilon^{\beta} \int_{|y|\leq 8\varepsilon}\frac{\nu(|y|)}{|y|^{d+\beta}} \right)
\end{align*}
Then one can observe that $\mathcal{C}(\varepsilon,\varphi,\nu) \downarrow 0$ as $\varepsilon \downarrow 0$. Hence, we can take $\varepsilon_{0}>0$ small enough so that
\begin{align}\label{eqn 11.21.18:50-2}
\|\varphi_{z}(\mathcal{L}^{z}_{V} - \mathcal{L}^{0}_{V})u\|_{L_{q,p}(w_{2},w_{1,z},t)} \leq \frac{1}{2} \| \phi(\Delta)u_{z}(t,x) \|_{L_{q,p}(w_{2},w_{1,z},t)} \quad \forall\, \varepsilon \leq \varepsilon_{0}.
\end{align}
Since $(1-\varphi_{z}(x))=0$ for $|x| < 2\varepsilon$,   $(1-\varphi_{z})(\mathcal{L}^{z}_{V} - \mathcal{L}^{0}_{V})u$ satisfies
\begin{align*}
& \int_{\mathbb{R}^{d}} \left|  (1-\varphi_{z})(x) \left( \mathcal{L}^{z}_{V} - \mathcal{L}^{0}_{V})u_{z}(t,x) \right) \right|^{p} w_{1}(x) \mathrm{d}x
\\
& = \int_{\mathbb{R}^{d}} \left| 1-\varphi_{z}(x) \right|^{p} \left|  \int_{|y| \geq \varepsilon} u_{z}(x+y) (a_{z}(t,x,y) - a(t,z,y)) j^{\phi}_{d}(|y|) \mathrm{d}y \right|^{p} w_{1,z}(x) \mathrm{d}x.
\end{align*}
One can follow the argument in \eqref{eqn 11.21.18:44} (with $\varepsilon 2^{j}$ in place of $2^{j}$) to obtain
\begin{align*}
\int_{|y| \geq \varepsilon} u_{z}(x+y) (a_{z}(t,x,y) - a(t,z,y)) j^{\phi}_{d}(|y|) \mathrm{d}y \leq C \left(\sum_{j=0}^{\infty}2^{-2\delta_{0}(j+1)} \right) \phi(\varepsilon^{-2}) \mathbb{M}u_{z}(x).
\end{align*}
Hence, it follows that
\begin{align}\label{eqn 11.21.18:51}
\|(1-\varphi_{z})(\mathcal{L}^{z}_{V}-\mathcal{L}^{0}_{V})u_{z}\|_{L_{q,p}(w_{2},w_{1,z},t)} \leq C \phi(\varepsilon^{-2}) \|u_{z}\|_{L_{q,p}(w_{2},w_{1,z},t)}.
\end{align}
Combining \eqref{eqn 11.21.18:50}, \eqref{eqn 11.21.18:50-2}, and \eqref{eqn 11.21.18:51}, for any $0<t \leq T$,
\begin{align*}
\|\phi(\Delta)u_{z}\|_{L_{q,p}(w_{2},w_{1,z},t)} \leq C \|f_{z}\|_{L_{q,p}(w_{2},w_{1,z},t)} + C \phi(\varepsilon^{-2}) \|u_{z}\|_{L_{q,p}(w_{2},w_{1,z},t)} \quad \forall\, \varepsilon \leq \varepsilon_{0}.
\end{align*}
Finally, using the change of variable $x+z \mapsto x$, we get the desired estimate. The lemma is proved.
\end{proof}

\begin{thm}\label{thm 01.28.15:20}
Let $1<p_{0},p_{1}<\infty$, $w_{i} \in A_{p_{i}}(\mathbb{R}^{d})$, $\gamma_{i} \in \mathbb{R}$ ($i=0,1$) and let $\theta \in (0,1)$. Also, let $\phi$ be a Bernstein function satisfying Assumption \ref{asm 07.21.13:06}($\delta_0,\delta_1$) with $0<\delta_0\leq\delta_1\leq 1$. Suppose that
\begin{gather*}
\gamma = (1-\theta)\gamma_{0} +\theta \gamma_{1}, \quad w = (w_{0})^{p(1-\theta)/p_{0}}(w_{1})^{p\theta/p_{1}},
\\
\frac{1}{p} = \frac{1-\theta}{p_{0}} + \frac{\theta}{p_{1}}, 
\end{gather*}
Then it follows that
\begin{align}\label{eqn 01.28.14:53}
[ H^{\phi,\gamma_{0}}_{p_{0}}(w_{0}), H^{\phi,\gamma_{1}}_{p_{1}}(w_{1})  ]_{[\theta]} = H^{\phi,\gamma}_{p}(w),
\end{align}
where the space $[ H^{\phi,\gamma_{0}}_{p_{0}}(w_{0}), H^{\phi,\gamma_{1}}_{p_{1}}(w_{1})  ]_{[\theta]}$ is the complex interpolation space. Moreover, it also holds that
\begin{align}\label{eqn 01.28.15:17}
\|u\|_{H^{\phi,\gamma}_{p}(w)} &\leq  \|u\|^{1-\theta}_{H^{\phi,\gamma_{0}}_{p_{0}}(w_{0})} \|u\|^{\theta}_{H^{\phi,\gamma_{1}}_{p_{1}}(w_{1})}   \nonumber
\\
&\leq \varepsilon \|u\|_{H^{\phi,\gamma_{0}}_{p_{0}}(w_{0})} + C(\varepsilon)\|u\|_{H^{\phi,\gamma_{1}}_{p_{1}}(w_{1})} \quad \forall\, \varepsilon>0.
\end{align}
\end{thm}
\begin{proof}
Note that \eqref{eqn 01.28.15:17} is a direct consequence of the definition of complex interpolation space, and Young's inequality. Therefore, it suffices to show \eqref{eqn 01.28.14:53}. Bring the standard Littlewood-Paley projection operators in \textbf{Step 2} of the proof of Lemma \ref{lem 11.21.17:27}, and define
\begin{align*}
\mathcal{F}[\mathcal{S}_{0}] = 1 - \sum_{j=1}^{\infty}\mathcal{F}[\Psi_{j}].
\end{align*}
Since $\phi$ is a Bernstein function satisfying the scaling condition, for any $1<q<\infty$, $s \in \mathbb{R}$, and any $\bar{w} \in A_{q}(\mathbb{R}^{d})$, we certainly derive the following equivalence 
\begin{align}\label{eqn 01.28.13:41}
\|u\|_{H^{\phi,s}_{q}(\bar{w})} \simeq \|S_{0} \ast u\|_{L_{q}(\bar{w})} + \left\| \left( \sum_{j=1}^{\infty} |\phi(2^{2j})^{s/2} \Psi_{j}\ast u |^{2} \right)^{1/2} \right\|_{L_{q}(\bar{w})}
\end{align}
by exploiting Assumption \ref{asm 07.21.13:06}, and  \eqref{eqn 01.28.13:23}. Indeed, Assumption \ref{asm 07.21.13:06} gives
\begin{align*}
\phi(2^{2j}) \simeq \phi(|\xi|^{2}) \quad \forall\, 2^{j-1} \leq |\xi| \leq 2^{j+1}.
\end{align*}
Using this, one can derive \eqref{eqn 01.28.13:41} by following the argument in \cite[Theorem 3.2]{P23} (we also recall Remark \ref{rmk 01.28.14:36} when one encounters some Fourier multipliers). 
Using the relation \eqref{eqn 01.28.13:41}, one can prove the desired result by following the proof of \cite[Theorem 3.3]{P23} with the space $\ell^{\phi,s}_{2}$ of sequences $\{a_{0},a_{1},\cdots\}$ with
\begin{align*}
|a|_{\ell^{\phi,s}_{2}} := \left(  \sum_{j=0}^{\infty} |\phi(2^{2j})^{s/2} a_{j}|^{2} \right)^{1/2}
\end{align*}
in place of $\ell^{s}_{2}$ therein.
\end{proof}

\begin{proof}[Proof of Theorem \ref{thm 11.24.15:58}]
Due to the method of continuity and Lemma \ref{lem 11.21.17:27}, we only need to prove the a priori estimate of solutions. Suppose that we prove  \eqref{eqn 12.02.15:54} for the case when $p=q$. Then what we have is that the estimate
\begin{align*}
\int_{0}^{T} \|u(t,\cdot)\|^{q}_{H^{\phi,2}_{p}(w_{1})} w_{2}(t)\mathrm{d}t \leq C \int_{0}^{T} \|f(t,\cdot)\|^{q}_{L_{p}(w_{1})} w_{2}(t) \mathrm{d}t
\end{align*}
holds for $q=q_{0}=p$ and for any $w_{2} \in A_{q_{0}}(\mathbb{R})$. Hence, applying the extrapolation theorem \cite[Theorem 1.4]{extrapolation} (also see  \cite[Theorem 2.5]{DK18} for detail) we get \eqref{eqn 12.02.15:54} for all $q \in (1,\infty)$ and for all $w_{2} \in A_{q}(\mathbb{R})$. Therefore, it suffices to consider only the case $p=q$. Suppose that $u\in \mathbb{H}^{\phi,2}_{p,p,0}(w_{2},w_{1},T)$ is a solution to \eqref{eqn 11.24.16:10}.
Take $\varepsilon_{0}$ from Lemma \ref{lem 11.24.16:13} and take a nonnegative smooth function $\varphi \in C^{\infty}_{c}(B_{\varepsilon})$ ($\varepsilon \leq \varepsilon_{0}$) with
\begin{align*}
\int_{\mathbb{R}^{d}} \left|  \varphi(x) \right|^{p} \mathrm{d}x = 1.
\end{align*}
For each $z\in \mathbb{R}^{d}$, denote $\varphi_{z}(x) = \varphi(x-z)$. Also, for a function $b$, define
\begin{align}\label{eqn 12.02.14:46}
\mathcal{I}^{b}(u,\varphi_{z}) = \int_{\mathbb{R}^{d}} \left(  u(x+y) - u(x)  \right) \left( \varphi(x+y-z) - \varphi(x-z) \right) b(t,x,y) j^{\phi}_{d}(|y|) \mathrm{d}y.
\end{align}
First, observe that
\begin{align*}
&\int_{\mathbb{R}^{d}} \left| \phi(\Delta)u(t,x) \right|^{p} w_{1}(x) \mathrm{d}x
\\
&= \int_{\mathbb{R}^{d}} \int_{\mathbb{R}^{d}} \left| \varphi(x-z) \phi(\Delta)u(t,x) \right|^{p} \mathrm{d}z w_{1}(x) \mathrm{d}x
\\
& \leq C \int_{\mathbb{R}^{d}} \int_{\mathbb{R}^{d}} \left| \phi(\Delta)(u\varphi_{z})(t,x) \right|^{p} \mathrm{d}z w_{1}(x) \mathrm{d}x
\\
& \quad + C \int_{\mathbb{R}^{d}} \int_{\mathbb{R}^{d}} \left| u(t,x) \phi(\Delta)[\varphi(\cdot-z)](x) \right|^{p} \mathrm{d}z w_{1}(x) \mathrm{d}x
\\
& \quad + C \int_{\mathbb{R}^{d}} \int_{\mathbb{R}^{d}} \left|\mathcal{I}^{1}(u(t,\cdot),\varphi_{z})(x) \right|^{p} \mathrm{d}z w_{1}(x) \mathrm{d}x.
\end{align*}
By Lemma \ref{lem 07.17.19:30}, 
\begin{align*}
\int_{\mathbb{R}^{d}} \int_{\mathbb{R}^{d}} \left| \mathcal{I}^{1}(u(t,\cdot),\varphi_{z})(x) \right|^{p} \mathrm{d}z w_{1}(x) \mathrm{d}x \leq C \|\varphi\|^{p}_{H^{1}_{p}} \|u(t,\cdot)\|^{p}_{H^{\phi,2\gamma}_{p}(w_{1})},
\end{align*}
where $1-\frac{1}{2\delta_{1}} < \gamma <1/2$.  
Since $\phi(\Delta)[\varphi(\cdot-z)](x) = \phi(\Delta)\varphi(x-z)$, it is easy to see that  (recall \eqref{eqn 07.03.16:05})
\begin{align}\label{eqn 11.28.11:35}
\int_{\mathbb{R}^{d}} \int_{\mathbb{R}^{d}} \left| u(t,x) \phi(\Delta)[\varphi(\cdot-z)] (x) \right|^{p} \mathrm{d}z w_{1}(x)\mathrm{d}x
 \leq C \|\varphi\|^{p}_{H^{2}_{p}} \|u(t,\cdot)\|^{p}_{L_{p}(w_{1})}.
\end{align}
Hence, we deduce that
\begin{align}\label{eqn 11.28.11:48}
&\|\phi(\Delta)u\|_{L_{p,p}(w_{2},w_{1},t)}   \nonumber
\\
& \leq C \|\varphi\|_{H^{2}_{p}} \|u\|_{L_{p,p}(w_{2},w_{1},t)} + C \|\varphi\|_{H^{1}_{p}} \|u\|_{H^{\phi,2\gamma}_{p,p}(w_{2},w_{1},t)}     \nonumber
\\
& \quad \quad + C \left( \int_{0}^{t}  \int_{\mathbb{R}^{d}} \int_{\mathbb{R}^{d}} \left| \phi(\Delta)(u(s,\cdot)\varphi_{z})(x) \right|^{p} \mathrm{d}z w_{1}(x)\mathrm{d}x w_{2}(s)\mathrm{d}s \right)^{1/p}.
\end{align}
For $u\varphi_{z}$, one can check that $u\varphi_{z}$ satisfies (recall \eqref{eqn 12.02.14:46})
\begin{align*}
\partial_{t}(u\varphi_{z}) = \mathcal{L}^{a}_{V}(u\varphi_{z}) + f\varphi_{z} - u\mathcal{L}^{a}_{V}\varphi_{z} - \mathcal{I}^{a}(u,\varphi_{z}).
\end{align*}
Applying Lemma \ref{lem 11.24.16:13} to $u\varphi_{z} \in \mathbb{H}^{\phi,2}_{q,p,0}(w_{2},w_{1},T)$, it follows that
\begin{align}\label{eqn 11.28.11:08}
& \|\phi(\Delta)(u\varphi_{z})\|^{p}_{L_{p,p}(w_{2},w_{1},t)}   \nonumber
\\
&\quad \leq C \left( \|f\varphi_{z}\|^{p}_{L_{p,p}(w_{2},w_{1},t)} + \|u\mathcal{L}^{a}_{V}\varphi_{z}\|^{p}_{L_{p,p}(w_{2},w_{1},t)} + \|\mathcal{I}^{a}(u,\varphi_{z})\|^{p}_{L_{p,p}(w_{2},w_{1},t)} \right)   \nonumber
\\
& \quad \quad + C \phi(\varepsilon^{-2}) \|(u\varphi_{z})\|^{p}_{L_{p,p}(w_{2},w_{1},t)}.
\end{align}
Using Lemma \ref{lem 07.17.19:30} again,
\begin{align}\label{eqn 11.28.11:41}
\int_{\mathbb{R}^{d}} \int_{\mathbb{R}^{d}} \left| \mathcal{I}^{a}(u(t,\cdot),\varphi_{z})(x) \right|^{p} \mathrm{d}z w_{1}(x) \mathrm{d}x \leq C \|\varphi\|^{p}_{H^{1}_{p}} \|u(t,\cdot)\|^{p}_{H^{\phi,2\gamma}_{p}(w_{1})} \quad \forall\, 0<t \leq T.
\end{align}
The same argument for \eqref{eqn 11.28.11:35} deduces
\begin{align}\label{eqn 11.28.11:41-2}
\int_{\mathbb{R}^{d}} \int_{\mathbb{R}^{d}} \left|u(t,x)\mathcal{L}^{a}_{V}\varphi_{z}(x) \right|^{p} \mathrm{d}zw_{1}(x)\mathrm{d}x \leq C \|\varphi\|^{p}_{H^{2}_{p}} \|u(t,\cdot)\|^{p}_{L_{p}(w_{1})}.
\end{align}
Since we assume $p=q$, by integrating both sides of \eqref{eqn 11.28.11:08} in $z$, using Fubini's theorem, and applying \eqref{eqn 11.28.11:41}, \eqref{eqn 11.28.11:41-2}, we have
\begin{align*}
&\int_{0}^{t} \int_{\mathbb{R}^{d}} \int_{\mathbb{R}^{d}} \left|  \phi(\Delta)(u\varphi_{z})(s,x) \right|^{p} \mathrm{d}z w_{1}(x) \mathrm{d}x w_{2}(s) \mathrm{d}s
\\
& \quad \leq C \|f\|^{p}_{L_{p,p}(w_{2},w_{1},t)} + C \|\varphi\|^{p}_{H^{1}_{p}(\mathrm{d}x)} \|u\|^{p}_{H^{\phi,2\gamma}_{p,p}(w_{2},w_{1},t)}
\\
& \quad \quad +   C \|\varphi\|^{p}_{H^{2}_{p}}\|u\|^{p}_{L_{p,p}(w_{2},w_{1},t)} + C \phi(\varepsilon^{-2}) \|(u\varphi_{z})\|^{p}_{L_{p,p}(w_{2},w_{1},t)}.
\end{align*}
for all $0<t \leq T$. Plugging the above inequality into \eqref{eqn 11.28.11:48}, it follows that
\begin{align*}
\|\phi(\Delta)u\|^{p}_{L_{p,p}(w_{2},w_{1},t)} & \leq C\|f\|^{p}_{L_{p,p}(w_{2},w_{1},t)} + C \|\varphi\|^{p}_{H^{1}_{p}} \|u\|^{p}_{H^{\phi,2\gamma}_{p,p}(w_{2},w_{1},t)}
\\
& \quad + C\|\varphi\|^{p}_{H^{2}_{p}}\|u\|^{p}_{L_{p,p}(w_{2},w_{1},t)} +   C\phi(\varepsilon^{-2}) \|u\|^{p}_{L_{p,p}(w_{2},w_{1},t)} \quad \forall\, 0<t \leq T.
\end{align*}
Due to the condition $1-1/2\delta_{1} < \gamma <1/2$, one can use Theorem \ref{thm 01.28.15:20} so that
\begin{align}\label{eqn 11.28.17:59}
\|\phi(\Delta)u\|^{p}_{L_{p,p}(w_{2},w_{1},t)} &\leq C \|f\|^{p}_{L_{q,p}(w_{2},w_{1},t)} + \varepsilon \|u\|^{p}_{H^{\phi,2}_{p,p}(w_{2},w_{1},t)} + C(\varepsilon)\|u\|^{p}_{L_{p,p}(w_{2},w_{1},t)}
\end{align}
for all $0<t \leq T$. Since $u\in \mathbb{H}^{\phi,2}_{p,p,0}(w_{2},w_{1},t)$ is a solution to \eqref{eqn 11.24.16:10}, one can check that
\begin{align*}
u(t,x) = \int_{0}^{t} \mathcal{L}^{a}_{V}u(s,x) + f(s,x) \mathrm{d}s \quad  \forall\, 0<s\leq t \leq T.
\end{align*}
Using H\"older's inequality, 
\begin{align*}
|u(s,x)|^{p} &= \left| \int_{0}^{s} \mathcal{L}^{a}_{V}u(r,x) + f(r,x) \mathrm{d}r \right|^{p}
\\
&= \left| \int_{0}^{s} \left( \mathcal{L}^{a}_{V}u(r,x) + f(r,x) \right) \frac{w^{1/p}_{2}(r)}{w^{1/p}_{2}(r)}\mathrm{d}r \right|^{p}
\\
& \leq C \left(\int_{0}^{t} w^{-1/(p-1)}_{2}(r) \mathrm{d}r \right)^{p-1} \int_{0}^{s} \left( |\mathcal{L}^{a}_{V}u(r,x)|^{p} + |f(r,x)|^{p} \right) w_{2}(r) \mathrm{d}r
\\
& \leq C \left(\int_{0}^{T} w^{-1/(p-1)}_{2}(r) \mathrm{d}r \right)^{p-1} \int_{0}^{s} \left( |\mathcal{L}^{a}_{V}u(r,x)|^{p} + |f(r,x)|^{p} \right) w_{2}(r) \mathrm{d}r,
\end{align*}
which yields that (let $\bar{w}_{2}(T) = \int_{0}^{T} w^{-1/(p-1)}_{2}(r) \mathrm{d}r$)
\begin{align*}
 \int_{\mathbb{R}^{d}} |u(s,x)|^{p} w_{1}(x) \mathrm{d}x &\leq C \bar{w}_{2}(T)  \int_{0}^{s} \int_{\mathbb{R}^{d}} \left( | \mathcal{L}^{a}_{V}u(r,x)|^{p} + |f(r,x)|^{p} \right) w_{1}(x)\mathrm{d}x w_{2}(r)\mathrm{d}r
\\ 
& \leq C \bar{w}_{2}(T) \|u\|^{p}_{H^{\phi,2}_{p,p}(w_{2},w_{1},s)} + C \bar{w}_{2}(T) \|f\|^{p}_{L_{p,p}(w_{2},w_{1},s)} 
\\
& \leq C \bar{w}_{2}(T)\|u\|^{p}_{H^{\phi,2}_{p,p}(w_{2},w_{1},s)} + C \bar{w}_{2}(T) \|f\|^{p}_{L_{q,p}(w_{2},w_{1},t)}.
\end{align*}
Thus 
\begin{align*}
\|u\|^{p}_{L_{p,p}(w_{2},w_{1},t)} & \leq C\bar{w}_{2}(T) \int_{0}^{t} w_{2}(s) \|u\|^{p}_{H^{\phi,2}_{p,p}(w_{2},w_{1},s)} \mathrm{d}s 
\\
& \quad + C\bar{w}_{2}(T)\left(\int_{0}^{T}w_{2}(t)\mathrm{d}t \right)\|f\|^{p}_{L_{p,p}(w_{2},w_{1},t)}
\end{align*}
follows. Combining the above inequality and \eqref{eqn 11.28.17:59},
\begin{align*}
&\|u\|^{p}_{H^{\phi,2}_{p,p}(w_{2},w_{1},t)} 
\\
& \leq C \bar{w}_{2}(T) \int_{0}^{t} w_{2}(s) \|u\|^{p}_{H^{\phi,2}_{p,p}(w_{2},w_{1},s)} \mathrm{d}s + C(\bar{w}_{2},T) \|f\|^{p}_{L_{p,p}(w_{2},w_{1},t)}
\\
& \quad + C \varepsilon \|u\|^{p}_{H^{\phi,2}_{p,p}(w_{2},w_{1},t)} + C(\varepsilon) \|u\|^{p}_{L_{p,p}(w_{2},w_{1},t)}  
\\
& \leq C(\bar{w}_{2},T,\varepsilon) \int_{0}^{t} w_{2}(s) \|u\|^{p}_{H^{\phi,2}_{p,p}(w_{2},w_{1},t)} \mathrm{d}s + C_{1}\varepsilon \|u\|^{p}_{H^{\phi,2}_{p,p}(w_{2},w_{1},t)} + C(\bar{w}_{2},T) \|f\|^{p}_{L_{p,p}(w_{2},w_{1},t)}
\end{align*}
holds for all $0<t\leq T$, where the constants $C_{1}$ is independent of $\varepsilon$. Hence, taking $\varepsilon$ small enough, and then applying Gronwall's inequality, we get
\begin{align*}
\|u\|^{p}_{H^{\phi,2}_{p,p}(w_{2},w_{1},T)} \leq C \|f\|^{p}_{L_{p,p}(w_{2},w_{1},T)}.
\end{align*}
The theorem is proved.
\end{proof}

\begin{thm}\label{thm 01.26.16:45}
Let $\phi$ be a Bernstein function satisfying Assumption \ref{asm 07.21.13:06}($\delta_0,\delta_1$) with $0<\delta_0\leq\delta_1<1$, $1<p<\infty$, and $w_{1} \in A_{p}(\mathbb{R}^{d})$.
Suppose that $V(y)=j_d^{\phi}(|y|)$, and the function $a=a(x,y)$ satisfies Assumption \ref{asm 11.10.16:43}.  
Then there exists a number $\lambda_{0}=\lambda_{0}(c_{0},c_{1},\delta_{0},\delta_{1},a,d,p,[w_{1}]_{A_{p}(\mathbb{R}^{d})})$ such that  for any $\lambda \geq \lambda_{0}$, and  $f\in L_{p}(w_{1})$, there exists a unique solution $u \in H^{\phi,2}_{p}(w_{1})$ to the equation
\begin{align*}
\mathcal{L}^{a}_{V}u(x) -\lambda u(x) = f(x) \quad \forall \, x\in \mathbb{R}^{d}
\end{align*}
satisfying
\begin{align*}
\|\phi(\Delta)u\|_{L_{p}(w_{1})} + \lambda \|u\|_{L_{p}(w_{1})} \leq C \|f\|_{L_{p}(w_{1})},
\end{align*}
where the constant $C=C_{0}(c_{0},c_{1},\delta_{0},\delta_{1},a,d,p,[w_{1}]_{A_{p}(\mathbb{R}^{d})})$ is independent of $\lambda$. Moreover, if $a=a(y)$ does not depend on $x$, then the statement holds for all $\lambda >0$. 
\end{thm}

\begin{proof}
\textbf{Step 1.} Suppose that $a=a(y)$ depends only on $y$.
\\
Without loss of generality, we assume that $f \in C^{\infty}_{c}(\mathbb{R}^{d})$. Define
\begin{align*}
G_{\lambda}(x) = \int_{0}^{\infty} \mathrm{e}^{-\lambda t} p_{d}(t,x) \mathrm{d}t,
\end{align*}
where $p_{d}(t,x)$ is the transition density of a L\'evy process with the triplet $(0,0,\nu)$, where
\begin{align*}
\nu(B) := \int_{B} a(y) j^{\phi}_{d}(\mathrm{d}y).
\end{align*}
Then it is well-known that
\begin{align*}
u(x) : = [G_{\lambda} \ast (-f)] (x) = \int_{\mathbb{R}^{d}} G_{\lambda}(x-y) f(y) \mathrm{d}y
\end{align*}  
Is the classical solution to $=\mathcal{L}^{a}_{V}u - \lambda u = f$.
If we rewrite $u$ as
\begin{align*}
u(x) &= \int_{\mathbb{R}^{d}} \int_{0}^{\infty} \mathrm{e}^{-\lambda t} p_{d}(t,y-x) \mathrm{d}y f(y) \mathrm{d}y= \int_{0}^{\infty} \mathrm{e}^{-\lambda t} \int_{\mathbb{R}^{d}} p_{d}(t,y-x) f(y) \mathrm{d}y \mathrm{d}t,
\end{align*}
then one can follow the argument in \eqref{eqn 12.18.14:18} to get
\begin{align*}
u(x) \leq C \int_{0}^{\infty} \mathrm{e}^{-\lambda t}  \mathbb{M}f(x) \mathrm{d}t = \frac{C}{\lambda} \mathbb{M}f(x).
\end{align*}
This certainly gives $\lambda \|u\|_{L_{p}(w_{1})} \leq C \|f\|_{L_{p}(w_{1})}$, where the constant $C$ does not depend on $\lambda$. 
Define $v(t,x) = \mathrm{e}^{\lambda t}u(x) - u(x)$. Then for each $T>0$, $v \in \mathbb{H}^{\phi,2}_{p,p}(1,w_{1},T)$ is a solution to the $\partial_{t}v = \mathcal{L}^{a}_{V}v + \left( \mathcal{L}^{a}_{V}u - \mathrm{e}^{\lambda t} f \right)$.
We apply Lemma \ref{lem 11.21.17:27} to deduce
\begin{align}\label{eqn 01.16.17:03}
\| \phi(\Delta)v \|_{L_{p,p}(1,w_{1},T)} \leq C  \|\mathcal{L}^{a}_{V}u - \mathrm{e}^{\lambda \cdot}f\|_{L_{p,p}(1,w_{1},T)},
\end{align}
where the constant $C$ is independent of $f$ and $T$. Using the relation between $u$ and $v$, \eqref{eqn 01.16.17:03} turns into
\begin{align*}
\int_{0}^{T} (\mathrm{e}^{\lambda t}-1)^p \| \phi(\Delta) u\|^{p}_{L_{p}(w_{1})} \mathrm{d}t \leq C \left(  T  \|\mathcal{L}^{a}_{V}u\|^{p}_{L_{p}(w_{1})} +  \int_{0}^{T} \mathrm{e}^{\lambda p t} \|f\|^{p}_{L_{p}(w_{1})}\mathrm{d}t \right).
\end{align*}
Using Theorem \ref{thm 01.15.15:07}, we have
\begin{align*}
\int_{0}^{T} (\mathrm{e}^{\lambda t}-1)^p \| \phi(\Delta) u\|^{p}_{L_{p}(w_{1})} \mathrm{d}t \leq C \left(   T \|\phi(\Delta) u\|^{p}_{L_{p}(w_{1})} + \int_{0}^{T} \mathrm{e}^{\lambda p t} \|f\|^{p}_{L_{p}(w_{1})}\mathrm{d}t \right). 
\end{align*}
Since $\lambda >0$,  $\int_{0}^{T} (\mathrm{e}^{\lambda t}-1)^p \mathrm{d}t > 2 C T$ for all large $T$. 
By taking such $T$, we have 
\begin{align*}
\|\phi(\Delta)u \|_{L_{p}(w_{1})}   \leq C \| f \|_{L_{p}(w_{1})},
\end{align*}
where the constant $C$ does not depend on $\lambda$.

\textbf{Step 2.} Suppose that $a=a(x,y)$.
\\
Due to the result in \textbf{Step 1}, and the method of continuity, in this step, we only consider the a priori estimate of solutions. Suppose that $u \in H^{\phi,2}_{p}(w_{1})$ is a solution to 
\begin{align}\label{eqn 01.21.15:23-2}
\mathcal{L}^{a}_{V}u - \lambda u = f,
\end{align}
Then by similarly defining $v(t,x) := \mathrm{e}^{\lambda t}u(x) - u(x) \in \mathbb{H}^{\phi,2}_{p,p}(1,w_{1},1)$, one can deduce 
\begin{align*}
\int_{0}^{1} (\mathrm{e}^{\lambda t}-1)^p \| u \|^{p}_{H^{\phi,2}_{p}(w_{1})} \mathrm{d}t \leq C \|u \|^{p}_{H^{\phi,2}_{p}(w_{1})} + C \int_{0}^{1} \mathrm{e}^{\lambda p t} \|f\|^{p}_{L_{p}(w_{1})} \mathrm{d}t 
\end{align*}
by following \textbf{Step 1} (in this step, use Theorem \ref{thm 11.24.15:58} instead of Lemma \ref{lem 11.21.17:27}). 
Note that the constant $C$ does not depend on $\lambda$. 
Hence, by taking $\lambda_{0}$ large enough so that $\int_{0}^{1} (\mathrm{e}^{\lambda_{0} t}-1)^p \mathrm{d}t > 2C$, it follows that
\begin{align}\label{eqn 01.21.15:23}
 \|\phi(\Delta)u\|^{p}_{L_{p}(w_{1})} \leq  C\| u \|^{p}_{H^{\phi,2}_{p}(w_{1})} \leq 2C \|f\|^{p}_{L_{p}(w_{1})} \quad \forall\, \lambda \geq \lambda_{0}, 
\end{align}
where the constant $C$ does not depend on $\lambda$. 
Using \eqref{eqn 01.21.15:23}  and \eqref{eqn 01.21.15:23-2}, we also obtain
\begin{align*}
\lambda \|u\|_{L_{p}(w_{1})} = \| \mathcal{L}^{a}_{V}u - f \|_{L_{p}(w_{1})} \leq C\|u\|_{H^{\phi,2}_{p}(w_{1})} + \|f\|_{L_{p}(w_{1})} \leq C \|f\|_{L_{p}(w_{1})},
\end{align*}
where the constant $C$ does not depend on $\lambda$.
\end{proof}

\section{Proof of weighted Sobolev--H\"older embeddings: Theorem \ref{thm 05.19.17:17}}

We first record the spatial version of the embedding theorem. 
For $p_0\in(R_w,p)$, set
$$
    \delta_{\phi,\gamma,p_0,p}(r)
    :=
    \phi(r^{-2})^{-\gamma} r^{-dp_0/p},
    \qquad 0<r\le1.
$$
For an open set $D\subset \mathbb R^d$, define
$$
    [u]_{C^{\phi,\gamma}_{p_0,p}(D)}
    :=
    \sup_{\substack{x,y\in D\\0<|x-y|\le1}}
    \frac{|u(x)-u(y)|}
    {\delta_{\phi,\gamma,p_0,p}(|x-y|)}
$$
and
$$
    \|u\|_{C^{\phi,\gamma}_{p_0,p}(D)}
    :=
    \|u\|_{L_\infty(D)}
    +
    [u]_{C^{\phi,\gamma}_{p_0,p}(D)} .
$$
We write $C^{\phi,\gamma}_{p_0,p,\mathrm{loc}}(D)$ for the corresponding local space.

\begin{corollary}[Elliptic weighted Sobolev--H\"older embedding]
\label{cor 02.04.14:21}
Let $\phi$ be a Bernstein function satisfying Assumption \ref{asm 07.21.13:06}$(\delta_0,\delta_1)$ with $0<\delta_0\le\delta_1<1$, and let $\gamma\in(0,1]$ satisfy $2\delta_1\gamma<1$.
Suppose that $w\in A_\infty(\mathbb R^d)$, and take $p\in(1,\infty)$ and $p_0\in(R_w,p)$ such that $w\in A_p(\mathbb R^d)$ and
\begin{equation*}
    2\delta_0\gamma-\frac{dp_0}{p}>0 .
\end{equation*}
Then every $u\in H_p^{\phi,2\gamma}(w)$ admits a representative belonging to $C^{\phi,\gamma}_{p_0,p,\mathrm{loc}}(\mathbb R^d)$.
In particular,
$$
    H_p^{\phi,2\gamma}(\mathbb R^d,w(x)\,dx)
    \hookrightarrow
    C^{\phi,\gamma}_{p_0,p,\mathrm{loc}}(\mathbb R^d).
$$
\end{corollary}

It is natural to ask whether the parameter $p_0$ in Corollary \ref{cor 02.04.14:21} can be replaced by the endpoint $R_w$. 
Equivalently, one may ask whether the embedding
\begin{align}\label{eqn 07.09.14:07}
    H_p^{\phi,2\gamma}(w)
    \hookrightarrow
    C^{\phi,\gamma}_{R_w,p,\mathrm{loc}}(\mathbb R^d)
\end{align}
holds in general, where the modulus would be
\begin{align}\label{eqn 07.03.16:49}
    \delta_{\phi,\gamma,R_w,p}(r)
    =
    \phi(r^{-2})^{-\gamma}r^{-dR_w/p}.
\end{align}
The following example shows that \eqref{eqn 07.03.16:49} is false in general.
In particular, even in the stable case $\phi(\lambda)=\lambda^{\sigma/2}$, one cannot, in general, replace $p_0$ by $R_w$ in the corresponding local H\"older embedding.
\begin{example}
\label{rem:endpoint-Rw-counterexample}
We give an one-dimensional example. 
Let
$$
\frac12<\sigma<1,
\qquad
\phi(\lambda)=\lambda^{\sigma/2},
\qquad
\gamma=1,
\qquad
p=2.
$$
Then $\delta_0=\delta_1=\sigma/2$, and hence
$$
2\delta_1\gamma=\sigma<1.
$$
Moreover,
$
H_2^{\phi,2}(w)
=
H_2^\sigma(w)
$
up to the equivalence of norms. Let
$$
\alpha:=\sigma-\frac12>0, \quad \kappa>0, \quad \theta>2\kappa +1.
$$
Define, near the origin,
$$
\ell(x):=\log\frac{e}{|x|},
\qquad
w(x):=\ell(x)^{-\theta}.
$$
More precisely, we take $w(x)=\ell(x)^{-\theta}$ for $0<|x|<1/2$, and extend $w$ to a positive smooth function on $\mathbb R\setminus(-1/2,1/2)$ such that $w$ is bounded above and below away from zero there.
This modification outside a neighborhood of the origin does not affect the local $A_q$-behavior of the weight.

We claim that
\begin{equation}\label{eqn 07.09.13:12}
R_{w} = \inf \{q\in [1,\infty): w\in A_{q}(\mathbb{R}) \} = 1, \quad w \notin A_{1}(\mathbb{R}).
\end{equation}
To show this we show that $w \in A_{q}(\mathbb{R})$ for any $q >1$, but $w \notin A_{1}(\mathbb{R})$.
Indeed, the only intervals that require attention are those meeting the origin. 
For $a\in\mathbb R$, we have
\begin{equation}
\label{eqn:log-average}
\int_0^r
\left(\log\frac{\mathrm{e}}{t}\right)^a\,\mathrm dt
\simeq
r\left(\log\frac{\mathrm{e}}{r}\right)^a,
\quad 0<r<1/2.
\end{equation}
Indeed, putting $L=\log(\mathrm{e}/r)$ and changing variables $s=\log(\mathrm{e}/t)$, we get
$$
\int_0^r
\left(\log\frac{\mathrm{e}}{t}\right)^a\,\mathrm dt
=
\mathrm{e}\int_L^\infty s^a \mathrm{e}^{-s}\,\mathrm ds
=
rL^a
\int_0^\infty
\left(1+\frac{\tau}{L}\right)^a \mathrm{e}^{-\tau}\,\mathrm d\tau .
$$
Since $L\geq \log(2\mathrm{e})>1$ for $0<r<1/2$, the last integral is bounded above and below by positive constants depending only on $a$. 
Hence
$$
\int_0^r
\left(\log\frac{e}{t}\right)^a\,\mathrm dt
\simeq
rL^a
=
r\left(\log\frac{e}{r}\right)^a .
$$
Applying \eqref{eqn:log-average} with $a=-\theta$ gives
$$
\frac1{2r}\int_{-r}^r w(x)\,\mathrm dx
\simeq
\ell(r)^{-\theta},
\quad 0<r<1/2.
$$
Similarly, for $q>1$,
$$
w(x)^{-1/(q-1)}
=
\ell(x)^{\theta/(q-1)}
$$
near the origin, and hence
$$
\frac1{2r}
\int_{-r}^r
w(x)^{-1/(q-1)}\,\mathrm dx
\simeq
\ell(r)^{\theta/(q-1)}.
$$
Therefore,
$$
\left(
\frac1{2r}\int_{-r}^r w(x)\,\mathrm dx
\right)
\left(
\frac1{2r}
\int_{-r}^r w(x)^{-1/(q-1)}\,\mathrm dx
\right)^{q-1}
\simeq
\ell(r)^{-\theta}
\left(\ell(r)^{\theta/(q-1)}\right)^{q-1}
\simeq 1.
$$
Intervals not meeting the origin are harmless, since on such intervals the weight is locally comparable to a positive constant, or else they are controlled by an interval centered at the origin.
Hence $w\in A_q(\mathbb R)$ for every $q>1$.

On the other hand, $w\notin A_1(\mathbb R)$.
To see this, let $I_0=(-1/2,1/2)$.
If $0<|x|<1/4$, then
$$
\mathbb M w(x)
\geq
\frac1{|I_0|}
\int_{I_0}w(y)\,\mathrm dy
=:c_0>0,
$$
whereas
$$
w(x)=\left(\log\frac{e}{|x|}\right)^{-\theta}\to0
\qquad\text{as }x\to0.
$$
Thus
$$
\frac{\mathbb M w(x)}{w(x)}\to\infty
\qquad\text{as }x\to0,
$$
which contradicts the $A_1$ condition, and hence, \eqref{eqn 07.09.13:12} follows.

Let $\eta\in C_c^\infty((-1,1))$ satisfy $\eta=1$ on $(-1/2,1/2)$, and set $u(x):=\eta(x)|x|^\alpha\ell(x)^\kappa$.
Then $u(0)=0$, but
$$
\frac{|u(x)-u(0)|}{|x|^\alpha}
=
\ell(x)^\kappa
\to\infty
\qquad
\text{as } x\to0.
$$
Thus $u$ does not belong to $C^{\alpha}_{\mathrm{loc}}(\mathbb{R})$. 
Now we prove the estimate
$$
|(-\Delta)^{\sigma/2}u(x)|
\leq
C  |x|^{\alpha-\sigma-1} \left( \mathbf{1}_{|x| \leq 1} |x| \ell(x/2)^\kappa + \mathbf{1}_{|x| >1} 1 \right) \quad \forall\, x\in \mathbb{R}.
$$
Since $0<\sigma<1$, for $x\neq0$ we use the singular integral representation
$$
(-\Delta)^{\sigma/2}u(x)
=
c_\sigma
\int_{\mathbb R}
\frac{u(x+h)-u(x)}{|h|^{1+\sigma}}\,\mathrm dh.
$$
We split the integral into four regions.

First, if $|h|\leq |x|/2$, then $|x|/2 \leq |x+h| \leq 2|x|$.
Hence, by mean value theorem,
$$
|u(x+h)-u(x)|
\leq
C |h| \mathbf{1}_{|x| \leq 1} |x|^{\alpha-1}\ell(x/2)^\kappa \quad \forall\, |h| \leq |x|/2.
$$
Thus
\begin{equation*}
\begin{aligned}
\int_{|h|\leq |x|/2}
\frac{|u(x+h)-u(x)|}{|h|^{1+\sigma}}\,\mathrm dh
& \leq
C \mathbf{1}_{|x| \leq 1} |x|^{\alpha-1}\ell(x/2)^\kappa
\int_0^{|x|/2} h^{-\sigma}\,\mathrm dh
\\
&\leq
C \mathbf{1}_{|x| \leq 1} |x|^{\alpha-\sigma}\ell(x/2)^\kappa.
\end{aligned}
\end{equation*}

Second, if $|x|/2<|h|\leq2|x|$, then $|x+h|\leq3|x|$. Thus we obtain  
$$
|u(x+h)|+|u(x)|
\leq
C |x|^\alpha \left( \eta(x) \ell(x)^\kappa + \eta(x+h) \ell(x+h)^{\kappa} \right).
$$
Therefore, using \eqref{eqn:log-average}, and the fact that $\ell(h)^{\kappa}$ is a decreasing function in $h$,
\begin{align*}
 &\int_{|x|/2<|h|\leq 2|x|}
\frac{|u(x+h)-u(x)|}{|h|^{1+\sigma}}\,\mathrm dh
\\
&\leq
C |x|^\alpha \eta(x) \ell(x)^\kappa
\int_{|x|/2}^{2|x|}h^{-1-\sigma}\,\mathrm dh
 + C |x|^{\alpha}\int_{|x|/2}^{2|x|} h^{-1-\sigma} \eta(x+h) \ell(|x+h|)^{\kappa} \, \mathrm{d}h
 \\
&\leq 
C |x|^{\alpha} \eta(x) \ell(x)^{\kappa}\int_{|x|/2}^{2|x|}h^{-1-\sigma}\,\mathrm dh
 + C |x|^{\alpha-1-\sigma}\int_{0}^{(3|x| \wedge 1)} \eta(h) \ell(h)^{\kappa} \, \mathrm{d}h
 \\
& \leq 
C |x|^{\alpha-\sigma-1} \left( \mathbf{1}_{|x| \leq 1} |x| \ell(x/2)^{\kappa} + C \mathbf{1}_{|x| >1} 1  \right)
\end{align*}

Third, if $2|x|<|h|$, then $|h|/2 \leq |x+h| \leq 2|h|$.
Hence,
$$
|u(x+h)-u(x)|
\leq
C \mathbf{1}_{|x| \leq 1}\left(|h|^\alpha\ell(h/2)^\kappa+|x|^\alpha \ell(x)^\kappa\right) \quad \forall\, 2|x| \leq |h|.
$$
Consequently, using the fact that both $h^{\alpha-1-\sigma}$ and $\ell(h)^{\kappa}$ are decreasing functions in $h$,
\begin{align*}
&\int_{2|x|<|h|}
\frac{|u(x+h)-u(x)|}{|h|^{1+\sigma}}\,\mathrm dh
\\
\leq&
C
\mathbf{1}_{|x| \leq 1}\left( \int_{2|x|}^{\infty}
h^{\alpha-1-\sigma}\ell(h/2)^\kappa\,\mathrm dh
+
C |x|^\alpha \ell(x)^\kappa
\int_{2|x|}^{\infty}h^{-1-\sigma}\,\mathrm dh \right)
\\
\leq & C \mathbf{1}_{|x| \leq 1} |x|^{\alpha-\sigma} \ell(x/2)^{\kappa}.
\end{align*}

Combining the four estimates gives
$$
|(-\Delta)^{\sigma/2}u(x)|
\leq
C  |x|^{\alpha-\sigma-1} \left( \mathbf{1}_{|x| \leq 1} |x| \ell(x/2)^\kappa + \mathbf{1}_{|x| >1} 1 \right).
$$
Since $\alpha=\sigma-\frac12$, this becomes
$$
|(-\Delta)^{\sigma/2}u(x)|
\leq
C |x|^{-1/2-1}\left( \mathbf{1}_{|x| \leq 1} |x| \ell(x/2)^\kappa + \mathbf{1}_{|x| >1} 1 \right).
$$
Therefore,
$$
\int_{\mathbb{R}}
|(-\Delta)^{\sigma/2}u(x)|^2w(x)\,\mathrm{d}x
\leq
C
\int_0^{1}
x^{-1}
\left(\log\frac{\mathrm{e}}{x}\right)^{2\kappa-\theta}
\,\mathrm{d}x
+ 
C
\int_{1}^{\infty} x^{-3} 
\,\mathrm{d}x
$$
With the change of variables $t=\log(\mathrm{e}/x)$, the right hand side becomes
$$
C
\int_{1}^{\infty}
t^{2\kappa-\theta}\,\mathrm{d}t
+ 
C
\int_{1}^{\infty}
x^{-3}\,\mathrm{d}x
$$
which is finite because $\theta>2\kappa+1$.
Also,
$$
\int_{\mathbb R}|u(x)|^2w(x)\,\mathrm{d}x<\infty.
$$
Consequently,
$$
u\in H_2^\sigma(\mathbb R;w(x)\,\mathrm{d}x)
=
H_2^{\phi,2}(\mathbb R;w(x)\,\mathrm{d}x).
$$
So far we have constructed a function
$u\in H_2^{\phi,2}(w)$
which does not belong to
$C_{\mathrm{loc}}^{\alpha}(\mathbb R)$.
Therefore, we finally observe that the embedding \eqref{eqn 07.09.14:07} cannot be obtained in general.
\end{example}
We shall use the following notation. 
Let $E\subset\mathbb R^d$ and $I\subset\mathbb R$ be measurable sets with $0<w_1(E)<\infty$ and $0<w_2(I)<\infty$. 
We write
$$
    w_1(E):=\int_E w_1(x)\,\mathrm{d}x,
    \quad
    w_2(I):=\int_I w_2(t)\,\mathrm{d}t,
$$
and, for a set $I\times E \subset \mathbb{R}\times \mathbb{R}^{d}$,
$$
    w_{2,1}(I\times E)
    :=
    \int_I\int_E w_1(x)w_2(t)\,\mathrm{d}x\mathrm{d}t
    =
    w_2(I)w_1(E).
$$
For locally integrable functions, we denote the corresponding weighted
averages by
$$
    f|_{E,w_1}
    :=
    \frac{1}{w_1(E)}
    \int_E f(x)w_1(x)\,\mathrm{d}x,\quad
    g|_{I,w_2}
    :=
    \frac{1}{w_2(I)}
    \int_I g(t)w_2(t)\,\mathrm{d}t,
$$
and
$$
    h|_{I\times E,w_{2,1}}
    :=
    \frac{1}{w_2(I)w_1(E)}
    \int_I\int_E h(t,x)w_1(x)w_2(t)\,\mathrm{d}x\,\mathrm{d}t .
$$
For $\rho>0$, we define the intrinsic time interval $\mathbb{I}_{\phi,\gamma}(t,\rho):=[t,t+\phi(\rho^{-2})^{-\gamma})$, and the corresponding intrinsic cylinder $Q_{\rho}(t,x):=\mathbb{I}_{\phi,\gamma}(t,\rho)\times B_\rho(x)$.
With this notation, Theorem \ref{thm 05.19.17:17} follows from the
following proposition by a standard density argument.
\begin{prop}
\label{prop 05.19.17:17}
Suppose that the assumptions in Theorem \ref{thm 05.19.17:17} hold.
Then, for every $u\in C_c^\infty([0,\infty)\times\mathbb R^d)$, $(t,x) \in [0,\infty)\times \mathbb{R}^{d}$, and $\rho>0$, we have
\begin{equation}
\label{eqn 02.03.15:02}
\begin{aligned}
|u(t,x)|
&\leq
C
\frac{
\phi(\rho^{-2})^{-\gamma}
(\|\partial_tu\|_{L_{q,p}(w_2,w_1)}+\|\phi^\gamma(\Delta)u\|_{L_{q,p}(w_2,w_1)})
}{
[w_2(\mathbb{I}_{\phi,\gamma}(t,\rho))]^{1/q}
[w_1(B_\rho(x))]^{1/p}
}
\\
&\quad
+
C
\frac{\|u\|_{L_{q,p}(w_2,w_1)}}{[w_2(\mathbb{I}_{\phi,\gamma}(t,\rho))]^{1/q}[w_1(B_\rho(x))]^{1/p}}.
\end{aligned}
\end{equation}
Moreover, for any $r \geq 0$, $z\in\mathbb R^d$, $t,s\in(r-1,r+1) \cap [0,\infty)$, and $x,y\in B_1(z)$ with $|t-s|\leq 1$ and $|x-y|\leq1$, we have
\begin{equation}\label{eqn 02.03.11:16}
\begin{aligned}
& |u(t,x)-u(s,y)|
\\
&\leq
C
\frac{
|t-s|^{1-\frac{q_0}{q}}\phi^{-1}(|t-s|^{-1/\gamma})^{\frac{dp_0}{2p}}
}{
[w_{2}(\mathbb{I}_{\phi,\gamma}(r,4))]^{1/q}
[w_1(B_4(z))]^{1/p}
} (\|\partial_{t}u\|_{L_{q,p}(w_{2},w_{1})} +\|\phi^{\gamma}(\Delta)u\|_{L_{q,p}(w_{2},w_{1})})
\\
&\quad
+
C
\frac{
\phi(|x-y|^{-2})^{-\gamma\left(1-\frac{q_0}{q}\right)}
|x-y|^{-dp_0/p}
}{
[w_{2}(\mathbb{I}_{\phi,\gamma}(r,4))]^{1/q}
[w_1(B_4(z))]^{1/p}
} (\|\partial_{t}u\|_{L_{q,p}(w_{2},w_{1})} +\|\phi^{\gamma}(\Delta)u\|_{L_{q,p}(w_{2},w_{1})}).
\end{aligned}
\end{equation}
Here $C$ is independent of $u,t,s,x,y,z,r$, and $\rho$, and depends only on
$$
c_0,c_1,\delta_0,\delta_1,\gamma,d,p,q,p_0,q_0,
[w_1]_{A_p(\mathbb R^d)},[w_1]_{A_{p_0}(\mathbb R^d)},
[w_2]_{A_q(\mathbb R)},[w_2]_{A_{q_0}(\mathbb R)}.
$$
\end{prop}

\begin{proof}
First, we do some auxiliary computations. 

Let $(r,z) \in [0,\infty)\times \mathbb{R}^{d}$, $(t,x) \in Q_{R}(r,z)$, and $0< \rho \leq R$. We also denote
\begin{align*}
I_{1}(t,s,x) = \int_{s}^{t} |\partial_{t}u(\lambda,x)| \mathrm{d}\lambda, \quad I_{2}(s,x) = \mathbb{M}\phi^{\gamma}(\Delta)u(s,x).
\end{align*}
Then by the fundamental theorem of calculus and Proposition \ref{lem 10.30.10:35}, it follows that
\begin{equation}
\label{26.06.29.11.35}
\begin{aligned}
&\left| u(t,x) - u(s,\cdot)|_{B_{\rho}(z),w_{1}}  \right| 
\\
&= \left| u(t,x) -  \frac{1}{w(B_{\rho}(z))}\int_{B_{\rho}(z)} u(s,y) \, w_{1}(y) \mathrm{d}y  \right|
\\
& \leq \left| I_{1}(t,s,x) \right|+ \frac{C}{w(B_{\rho}(z))} \int_{B_{\rho}(z)} \phi(|x-y|^{-2})^{-\gamma} \left( I_{2}(s,x) + I_{2}(s,y) \right) w_{1}(y) \mathrm{d}y.
\end{aligned}
\end{equation}
By H\"older's inequality,
\begin{equation}
\label{26.06.29.11.36}
\begin{aligned}
    &\frac{1}{w(B_{\rho}(z))}\int_{B_{\rho}(z)} \phi(|x-y|^{-2})^{-\gamma} \left( I_{2}(s,x) + I_{2}(s,y) \right) w_{1}(y) \mathrm{d}y\\
    &\leq  \frac{\phi((2R)^{-2})^{-\gamma}}{w(B_{\rho}(z))}\int_{B_{\rho}(z)} \left( I_{2}(s,x) + I_{2}(s,y) \right) w_{1}(y)\mathrm{d}y 
\\
&\leq  \phi((2R)^{-2})^{-\gamma} \left(  \frac{1}{w(B_{\rho}(z))} \int_{B_{\rho}(z)} \left( \left| I_{2}(s,x) \right|^{p} + \left| I_{2}(s,y) \right|^{p}  \right) w_{1}(y)\mathrm{d}y \right)^{1/p}
\\
&\leq C\phi(R^{-2})^{-\gamma}  |I_{2}(s,x)|  + C \phi(R^{-2})^{-\gamma}  \left( \frac{1}{w(B_{r}(z))} \int_{\mathbb{R}^{d}} \left| I_{2}(s,y) \right|^{p} w_{1}(y) \mathrm{d}y \right)^{1/p}
\\
& = C \phi(R^{-2})^{-\gamma} \left( \left| I_{2}(s,x) \right| + \frac{1}{[w_{1}(B_{r}(z))]^{1/p}} \|I_{2}(s,\cdot) \|_{L_{p}(w_{1})}  \right).
\end{aligned}
\end{equation}
Combining \eqref{26.06.29.11.35} and \eqref{26.06.29.11.36}, we have
\begin{align*}
\left| u(t,x) - u|_{Q_{\rho}(r,z),w_{2,1}} \right|
&\leq \frac{C}{w_2(\mathbb{I}_{\phi,\gamma}(r,\rho))}\int_{\mathbb{I}_{\phi,\gamma}(r,\rho)}|u(t,x) - u(s,\cdot)|_{B_{\rho}(z),w_{1}} |w_{2}(s)\mathrm{d}s\\
&\leq C   \frac{1}{w_{2}(\mathbb{I}_{\phi,\gamma}(r,\rho))}\int_{\mathbb{I}_{\phi,\gamma}(r,\rho)} \left| I_{1}(t,s,x) \right| w_{2}(s) \mathrm{d}s + C \phi(R^{-2})^{-\gamma} I_{4}(s,x),
\end{align*}
where
\begin{align*}
I_{4}(s,x) &= \frac{1}{w_{2}(\mathbb{I}_{\phi,\gamma}(r,\rho))} \left(\int_{\mathbb{I}_{\phi,\gamma}(r,\rho)} \left| I_{2}(s,x) \right| w_{2}(s) \mathrm{d}s \right)
\\
&\quad + \frac{1}{w_{2}(\mathbb{I}_{\phi,\gamma}(r,\rho))}\left( \frac{1}{[w_{1}(B_{\rho}(z))]^{1/p}}\int_{\mathbb{I}_{\phi,\gamma}(r,\rho)} \|I_{2}(s,\cdot)\|_{L_{p}(w_{1})} w_{2}(s) \mathrm{d}s \right)
\\
& \leq \frac{1}{w_{2}(\mathbb{I}_{\phi,\gamma}(r,\rho))} \left(\int_{\mathbb{I}_{\phi,\gamma}(r,\rho)} \left|  I_{2}(s,x) \right| w_{2}(s) \mathrm{d}s \right)
\\
& \quad + \frac{1}{[w_{2}(\mathbb{I}_{\phi,\gamma}(r,\rho))]^{1/q}} \frac{1}{[w_{1}(B_{\rho}(z))]^{1/p}} \| \phi^{\gamma}(\Delta)u\|_{L_{q,p}(w_{2},w_{1})}.
\end{align*}
Integrating both sides over $Q_{R}(r,z)$, we get
\begin{align*}
&\int_{Q_{R}(r,z)} \left| u(t,x) - u|_{Q_{\rho}(r,z),w_{2}w_{1}} \right|w_1(x)\mathrm{d}xw_2(t)\mathrm{d}t
\\
& \leq   \frac{C}{w_{2}(\mathbb{I}_{\phi,\gamma}(r,\rho))} \int_{Q_{R}(r,z)} \int_{\mathbb{I}_{\phi,\gamma}(r,\rho)} \left| I_{1}(t,s,x) \right| w_{2}(s)\mathrm{d}s w_{1}(x)\mathrm{d}x w_{2}(t)\mathrm{d}t
\\
& \quad +  C \phi(R^{-2})^{-\gamma} \int_{Q_{R}(r,z)} I_{4}(s,x) w_{1}(x) \mathrm{d}x w_{2}(t)\mathrm{d}t
\\
& \leq   \frac{C}{w_{2}(\mathbb{I}_{\phi,\gamma}(r,\rho))} \int_{Q_{R}(r,z)} \int_{\mathbb{I}_{\phi,\gamma}(r,\rho)} \left| I_{1}(t,s,x) \right| w_{2}(s)\mathrm{d}s w_{1}(x)\mathrm{d}x w_{2}(t)\mathrm{d}t
\\
& \quad +  C \phi(R^{-2})^{-\gamma} \int_{\mathbb{I}_{\phi,\gamma}(r,R)} \left(  \frac{1}{w_{2}(\mathbb{I}_{\phi,\gamma}(r,\rho))} \int_{\mathbb{I}_{\phi,\gamma}(r,\rho)} \int_{B_{R}(z)} \left| I_{2}(s,x) \right| w_{1}(x) \mathrm{d}x w_{2}(s) \mathrm{d}s \right) w_{2}(t)\mathrm{d}t
\\
& \quad  + C \phi(R^{-2})^{-\gamma} \frac{w_{1}(B_{R}(z))w_{2}(\mathbb{I}_{\phi,\gamma}(r,R))}{[w_{2}(\mathbb{I}_{\phi,\gamma}(r,\rho))]^{1/q}[w_{1}(B_{\rho}(z))]^{1/p}}  \| \phi^{\gamma}(\Delta)u\|_{L_{q,p}(w_{2},w_{1})}\\
&=:\mathrm{A}_1(r,z,\rho,R)+\mathrm{A}_2(r,z,\rho,R)+\mathrm{A}_3(r,z,\rho,R).
\end{align*}
Since
\begin{align*}
\mathrm{A}_2(r,z,\rho,R)&\leq C \phi(R^{-2})^{-\gamma}\frac{w_{1}(B_{R}(z))w_{2}(\mathbb{I}_{\phi,\gamma}(r,R))}{[w_{1}(B_{R}(z))]^{1/p}[w_{2}(\mathbb{I}_{\phi,\gamma}(r,\rho))]^{1/q}} \| \phi^{\gamma}(\Delta)u\|_{L_{q,p}(w_{2},w_{1})}
\\
& \leq C\phi(R^{-2})^{-\gamma} \frac{w_{1}(B_{R}(z))w_{2}(\mathbb{I}_{\phi,\gamma}(r,R))}{[w_{1}(B_{\rho}(z))]^{1/p}[w_{2}(\mathbb{I}_{\phi,\gamma}(r,\rho))]^{1/q}} \| \phi^{\gamma}(\Delta)u\|_{L_{q,p}(w_{2},w_{1})}\\
&\leq C\mathrm{A}_3(r,z,\rho,R),
\end{align*}
and
\begin{align*}
& \mathrm{A}_1(r,z,\rho,R)
\\
& \leq \frac{1}{w_{2}(\mathbb{I}_{\phi,\gamma}(r,\rho))} \int_{\mathbb{I}_{\phi,\gamma}(r,R)} \int_{\mathbb{I}_{\phi,\gamma}(r,\rho)} \left| \int_{s}^{t}  \frac{w_{1}(B_{R}(z))}{[w_{1}(B_{R}(z))]^{1/p}} \|\partial_{t}u(\lambda,\cdot)\|_{L_{p}(w_{1})} \mathrm{d}\lambda \right| w_{2}(s)\mathrm{d}s w_{2}(t)\mathrm{d}t
\\
& \leq \int_{\mathbb{I}_{\phi,\gamma}(r,R)}   \phi(R^{-2})^{-\gamma}  \frac{w_{1}(B_{R}(z))}{[w_{1}(B_{R}(z))]^{1/p}} \mathcal{M}[ \|\partial_{t}u\|_{L_{p}(w_{1})}](t)   w_{2}(t)\mathrm{d}t
\\
&\leq C \phi(R^{-2})^{-\gamma} \frac{w_{2}(\mathbb{I}_{\phi,\gamma}(r,R))}{[w_{2}(\mathbb{I}_{\phi,\gamma}(r,R))]^{1/q}}  \frac{w_{1}(B_{R}(z))}{[w_{1}(B_{R}(z))]^{1/p}} \|\mathcal{M}[\|\partial_{t}u\|_{L_{p}(w_{1})}] \|_{L_{q}(w_{2})}
\\
&\leq C \phi(R^{-2})^{-\gamma} \frac{w_{2}(\mathbb{I}_{\phi,\gamma}(r,R))}{[w_{2}(\mathbb{I}_{\phi,\gamma}(r,R))]^{1/q}} \frac{w_{1}(B_{R}(z))}{[w_{1}(B_{R}(z))]^{1/p}} \|\partial_{t}u\|_{L_{q,p}(w_{2},w_{1})}
\\
&\leq C \phi(R^{-2})^{-\gamma}\frac{w_{2}(\mathbb{I}_{\phi,\gamma}(r,R))}{[w_{2}(\mathbb{I}_{\phi,\gamma}(r,\rho))]^{1/q}} \frac{w_{1}(B_{R}(z))}{[w_{1}(B_{\rho}(z))]^{1/p}} \|\partial_{t}u\|_{L_{q,p}(w_{2},w_{1})},
\end{align*}
we obtain
\begin{equation}
\label{eqn 06.01.14:40}
\begin{aligned}
& \left| u|_{Q_{R}(r,z),w_{2,1}} - u|_{Q_{\rho}(r,z),w_{2,1}}  \right|  
\\
&\leq \frac{1}{w_1(B_R(z))w_2(\mathbb{I}_{\phi,\gamma}(r,R))}\int_{Q_{R}(r,z)} \left| u(t,x) - u|_{Q_{\rho}(r,z),w_{2,1}} \right|w_1(x)\mathrm{d}x w_2(t)\mathrm{d}t\\
& \leq C  \phi(R^{-2})^{-\gamma}\frac{\|\partial_{t}u\|_{L_{q,p}(w_{2},w_{1})}+ \|\phi^{\gamma}(\Delta)u\|_{L_{q,p}(w_{2},w_{1})}}{[w_{2}(\mathbb{I}_{\phi,\gamma}(r,\rho))]^{1/q}[w_{1}(B_{\rho}(z))]^{1/p}}  .
\end{aligned}
\end{equation}
Now, we consider \eqref{eqn 02.03.15:02}.
Due to the Lebesgue differentiation theorem,
\begin{align*}
|u(t,x)| &= |u(t,x) - u|_{Q_{\rho}(t,x),w_{2,1}}|+| \, u|_{Q_{\rho}(t,x),w_{2,1}} \, |
\\
 &\leq \sum_{k=0}^{\infty} |u|_{Q_{2^{-k-1}\rho}(t,x),w_{2,1}} - u|_{Q_{2^{-k}\rho}(t,x),w_{2,1}}| + | \, u|_{Q_{\rho}(t,x),w_{2,1}} \, |
\\
& =: \mathcal{A}_{\rho}(t,x) + | \, u|_{Q_{\rho}(t,x),w_{2,1}}  |.
\end{align*}
A direct computation yields
\begin{align*}
| \, u|_{Q_{\rho}(t,x),w_{2,1}}| \leq \frac{1}{[w_{2}(\mathbb{I}_{\phi,\gamma}(t,\rho))]^{1/q}[w_{1}(B_{\rho}(x))]^{1/p}} \|u\|_{L_{q,p}(w_{2},w_{1})}.
\end{align*}
Applyig \eqref{eqn 06.01.14:40} (with $R=2^{-k}\rho$, and $\rho=2^{-k-1}\rho$),
\begin{equation}
\label{eqn 02.03.15:10}
\begin{aligned}
\mathcal{A}_{\tau,\rho}(t,x) 
&\leq C \sum_{k=0}^{\infty} \frac{\phi(2^{2k}\rho^{-2})^{-\gamma} (\|\partial_{t}u\|_{L_{q,p}(w_{2},w_{1})} +  \|\phi^{\gamma}(\Delta)u\|_{L_{q,p}(w_{2},w_{1})} )}{[w_{2}(\mathbb{I}_{\phi,\gamma}(t,2^{-k-1}\rho))]^{1/q}[w_{1}(B_{2^{-k-1}\rho}(x))]^{1/p}}.
\end{aligned}
\end{equation}
Note that, using \eqref{eqn 02.03.21:00},  we have
\begin{align}\label{eqn 02.03.16:58}
\frac{1}{w_{1}(B_{2^{-k-1}\rho}(x))} \leq \frac{C}{w_{1}(B_{\rho}(x))} 2^{dp_{0}} (2^{dk})^{p_{0}},\quad \forall p_{0}\in(R_{w_{1}},p],
\end{align}
and
\begin{equation}\label{eqn 05.29.14:54}
\begin{aligned}
&\frac{\phi(2^{2k}\tau^{-2})^{-\gamma}[w_{2}(\mathbb{I}_{\phi,\gamma}(t,\tau))]^{1/q}}{[w_{2}(\mathbb{I}_{\phi,\gamma}(t,2^{-k-1}\tau))]^{1/q}}
\\
&\leq C  \phi(2^{2k}\tau^{-2})^{-\gamma}\left(\frac{\phi(\tau^{-2})^{-\gamma}}{\phi(2^{2(k+1)}\tau^{-2})^{-\gamma}}\right)^{\frac{q_{0}}{q}}
\\
&=C\phi(\tau^{-2})^{-\gamma}\left(\frac{\phi(2^{2(k+1)}\tau^{-2})^{\gamma}}{\phi(2^{2k}\tau^{-2})^{\gamma}}\right)^{\frac{q_0}{q}}\left(\frac{\phi(\tau^{-2})^{\gamma}}{\phi(2^{2k}\tau^{-2})^{\gamma}}\right)^{1-\frac{q_0}{q}}
\\
&\leq C\phi(\tau^{-2})^{-\gamma}  2^{-2k\delta_0\gamma\left(1-\frac{q_0}{q}\right)},
\end{aligned}
\end{equation}
for all $q_0\in(R_{w_2},q]$.
Hence, by selecting $p_{0},q_{0}$ satisfying \eqref{26.02.03.21.14}, and using \eqref{eqn 02.03.16:58} and \eqref{eqn 05.29.14:54}, it follows that
\begin{equation}
\label{eqn 02.03.15:25}
\begin{aligned}
&\mathcal{A}_{\rho}(t,x)
\\
& \leq C  \frac{\phi(\rho^{-2})^{-\gamma}(\|\partial_{t}u\|_{L_{q,p}(w_{2},w_{1})} +\|\phi^{\gamma}(\Delta)u\|_{L_{p}(w_1,w_2)} )}{[w_{2}(\mathbb{I}_{\phi,\gamma}(t,\rho))]^{1/q}[w_{1}(B_{\rho}(x))]^{1/p}} \sum_{k=0}^{\infty} 2^{-k(2\delta_0\gamma- \frac{2\delta_0\gamma q_0}{q}-\frac{dp_0}{p})}
\\
& \leq C \frac{\phi(\rho^{-2})^{-\gamma}(\|\partial_{t}u\|_{L_{q,p}(w_{2},w_{1})} +\|\phi^{\gamma}(\Delta)u\|_{L_{p}(w_1,w_2)} )}{[w_{2}(\mathbb{I}_{\phi,\gamma}(t,\rho))]^{1/q}[w_{1}(B_{\rho}(x))]^{1/p}}.
\end{aligned}
\end{equation}
Therefore, we have
\begin{align*}
|u(t,x)| &\leq \mathcal{A}_{\rho}(t,x) + |\, u|_{Q_{\rho}(t,x),w_{2,1}}| 
\\
& \leq C \frac{\phi(\rho^{-2})^{-\gamma}}{[w_{2}(\mathbb{I}_{\phi,\gamma}(t,\rho))]^{1/q}[w_{1}(B_{\rho}(x))]^{1/p}} \left( \|\partial_{t}u\|_{L_{q,p}(w_{2},w_{1})} + \|\phi^{\gamma}(\Delta)u\|_{L_{q,p}(w_{2},w_{1})} \right)
\\
& \quad  + \frac{C}{[w_{2}(\mathbb{I}_{\phi,\gamma}(t,\rho))]^{1/q}[w_{1}(B_{\rho}(x))]^{1/p}} \|u\|_{L_{q,p}(w_{2},w_{1})}.
\end{align*}
Thus \eqref{eqn 02.03.15:02} is proved. 

It only remains to consider \eqref{eqn 02.03.11:16}. 
For this, let $(r,z) \in [0,\infty) \times \mathbb{R}^{d}$, $t,s\in (r-1,r+1) \cap [0,\infty)$ and $x,y \in B_{1}(z)$. Take $0<\rho<1$ small enough so that
\begin{align}\label{eqn 09.02.14:08}
\phi(\rho^{-2})^{-\gamma} \leq |t-s| ,\quad  \rho \leq |x-y|.
\end{align}
Then due to the Lebesgue differentiation theorem, it follows that
\begin{align*}
& \left| u(t,x) - u(s,y) \right| 
\\
& \leq  \left| u(t,x) - u|_{Q_{2\rho}(t,x),w_{2,1}}  \right| + \left| u|_{Q_{2\rho}(t,x),w_{2,1}} - u|_{Q_{2\rho}(s,y),w_{2,1}} \right| + \left| u|_{Q_{2\rho}(s,y),w_{2,1}} - u(s,y) \right|
\\
& \leq \sum_{k=0}^{\infty} \left| u|_{Q_{2^{-k+1}\rho}(t,x),w_{2,1}} - u|_{Q_{2^{-k}\rho}(t,x),w_{2,1}} \right| 
+ \left| u|_{Q_{2\rho}(t,x),w_{2,1}} - u|_{Q_{2\rho}(s,y),w_{2,1}}  \right|
\\
& \quad + \sum_{k=0}^{\infty} \left| u|_{Q_{2^{-k+1}\rho}(s,y),w_{2,1}} - u|_{Q_{2^{-k}\rho}(s,y),w_{2,1}} \right|
\\
& := \mathcal{A}_{2\rho}(t,x) + \left| u|_{Q_{2\rho}(t,x),w_{2,1}} - u|_{Q_{2\rho}(s,y),w_{2,1}}  \right| + \mathcal{A}_{2\rho}(s,y).
\end{align*}
By following the argument from \eqref{eqn 02.03.15:10} to \eqref{eqn 02.03.15:25}, and using \eqref{eqn 02.03.21:00},
\begin{align*}
&\mathcal{A}_{2\rho}(t,x) + \mathcal{A}_{2\rho}(s,y) 
\\
&\leq  C \frac{\phi(\rho^{-2})^{-\gamma}}{[w_{2}(\mathbb{I}_{\phi,\gamma}(t,2\rho))]^{1/q}[w_{1}(B_{2\rho}(x))]^{1/p}}  (\|\partial_{t}u\|_{L_{q,p}(w_{2},w_{1})} +\|\phi^{\gamma}(\Delta)u\|_{L_{q,p}(w_{2},w_{1})}) 
\\
& \quad +  C \frac{\phi(\rho^{-2})^{-\gamma}}{[w_{2}(\mathbb{I}_{\phi,\gamma}(s,2\rho))]^{1/q}[w_{1}(B_{2\rho}(y))]^{1/p}}  (\|\partial_{t}u\|_{L_{q,p}(w_{2},w_{1})} +\|\phi^{\gamma}(\Delta)u\|_{L_{q,p}(w_{2},w_{1})}).
\end{align*}
follows. Similarly, 
\begin{align*}
& \left|  u|_{Q_{2\rho}(t,x),w_{2,1}} - u|_{Q_{2\rho}(s,y),w_{2,1}} \right|
\\
& \leq C \frac{\phi(\rho^{-2})^{-\gamma}}{[w_{2}(\mathbb{I}_{\phi,\gamma}(t,2\rho))]^{1/q}[w_{1}(B_{2\rho}(x))]^{1/p}}  (\|\partial_{t}u\|_{L_{q,p}(w_{2},w_{1})} +\|\phi^{\gamma}(\Delta)u\|_{L_{q,p}(w_{2},w_{1})})
\\
& \quad +  C \frac{\phi(\rho^{-2})^{-\gamma}}{[w_{2}(\mathbb{I}_{\phi,\gamma}(s,2\rho))]^{1/q}[w_{1}(B_{2\rho}(y))]^{1/p}}  (\|\partial_{t}u\|_{L_{q,p}(w_{2},w_{1})} +\|\phi^{\gamma}(\Delta)u\|_{L_{q,p}(w_{2},w_{1})}).
\end{align*}
Therefore, by exploiting \eqref{eqn 02.03.21:00} again, we obtain
\begin{align*}
& \left| u(t,x) - u(s,y)   \right|
\\
&\leq  C \frac{\phi(\rho^{-2})^{-\gamma}}{[w_{2}(\mathbb{I}_{\phi,\gamma}(t,2\rho))]^{1/q}[w_{1}(B_{2\rho}(x))]^{1/p}}  (\|\partial_{t}u\|_{L_{q,p}(w_{2},w_{1})} +\|\phi^{\gamma}(\Delta)u\|_{L_{q,p}(w_{2},w_{1})}) 
\\
& \quad +  C \frac{\phi(\rho^{-2})^{-\gamma}}{[w_{2}(\mathbb{I}_{\phi,\gamma}(s,2\rho))]^{1/q}[w_{1}(B_{2\rho}(y))]^{1/p}}  (\|\partial_{t}u\|_{L_{q,p}(w_{2},w_{1})} +\|\phi^{\gamma}(\Delta)u\|_{L_{q,p}(w_{2},w_{1})})
\\
&\leq C \frac{\phi(\rho^{-2})^{-\gamma\left(1-\frac{q_0}{q}\right)}\rho^{-\frac{dp_0}{p}}}{[w_{2}(\mathbb{I}_{\phi,\gamma}(t,2))]^{1/q}[w_{1}(B_{2}(x))]^{1/p}}  (\|\partial_{t}u\|_{L_{q,p}(w_{2},w_{1})} +\|\phi^{\gamma}(\Delta)u\|_{L_{q,p}(w_{2},w_{1})})
\\
& \quad +  C \frac{\phi(\rho^{-2})^{-\gamma\left(1-\frac{q_0}{q}\right)}\rho^{-\frac{dp_0}{p}}}{[w_{2}(\mathbb{I}_{\phi,\gamma}(s,2))]^{1/q}[w_{1}(B_{2}(y))]^{1/p}}  (\|\partial_{t}u\|_{L_{q,p}(w_{2},w_{1})} +\|\phi^{\gamma}(\Delta)u\|_{L_{q,p}(w_{2},w_{1})})
\\
&\leq C \frac{\phi(\rho^{-2})^{-\gamma\left(1-\frac{q_0}{q}\right)}\rho^{-\frac{dp_0}{p}}}{[w_{2}(\mathbb{I}_{\phi,\gamma}(r,4))]^{1/q}[w_{1}(B_{4}(z))]^{1/p}}  (\|\partial_{t}u\|_{L_{q,p}(w_{2},w_{1})} +\|\phi^{\gamma}(\Delta)u\|_{L_{q,p}(w_{2},w_{1})}).
\end{align*}
Note that from \eqref{26.02.03.21.14} 
\begin{align*}
1\leq \left(\frac{R}{r} \right)^{2\delta_{0}\gamma - 2\delta_{0}\gamma\frac{q_{0}}{q}-\frac{dp_{0}}{p}} \leq C \frac{\phi(R^{-2})^{-\gamma\left(1-\frac{q_{0}}{q}\right)}R^{-\frac{dp_{0}}{p}}}{\phi(r^{-2})^{-\gamma\left(1-\frac{q_{0}}{q}\right)}r^{-\frac{dp_{0}}{p}}} \quad \forall\, 0<r\leq R<\infty.
\end{align*}
Therefore, due to the choice of $\rho$ (recall \eqref{eqn 09.02.14:08})
\begin{align*}
& \left| u(t,x) - u(s,y)   \right|
\\
&\leq C \frac{\phi(\rho^{-2})^{-\gamma\left(1-\frac{q_0}{q}\right)}\rho^{-\frac{dp_0}{p}}}{[w_{2}(\mathbb{I}_{\phi,\gamma}(r,4))]^{1/q}[w_{1}(B_{4}(z))]^{1/p}}  (\|\partial_{t}u\|_{L_{q,p}(w_{2},w_{1})} +\|\phi^{\gamma}(\Delta)u\|_{L_{q,p}(w_{2},w_{1})})
\\
& \leq C \frac{|t-s|^{1-\frac{q_0}{q}}\phi^{-1}(|t-s|^{-1/\gamma})^{\frac{dp_0}{2p}}}{[w_{2}(\mathbb{I}_{\phi,\gamma}(r,4))]^{1/q}[w_{1}(B_{4}(z))]^{1/p}}  (\|\partial_{t}u\|_{L_{q,p}(w_{2},w_{1})} +\|\phi^{\gamma}(\Delta)u\|_{L_{q,p}(w_{2},w_{1})}) 
\\
& \quad + C \frac{ \phi(|x-y|^{-2})^{-\gamma\left(1-\frac{q_0}{q}\right)} |x-y|^{-\frac{dp_0}{p}}}{[w_{2}(\mathbb{I}_{\phi,\gamma}(r,4))]^{1/q}[w_{1}(B_{4}(z))]^{1/p}}   \ (\|\partial_{t}u\|_{L_{q,p}(w_{2},w_{1})} +\|\phi^{\gamma}(\Delta)u\|_{L_{q,p}(w_{2},w_{1})}).
\end{align*}
Thus, we get \eqref{eqn 02.03.11:16}. The proposition is proved.
\end{proof}

\section{Proof of the Well-Posedness of the Martingale problem and Krylov-type estimate: Theorem \ref{thm 02.26.19:05}}
\label{sec:proof-martingale-krylov}

The proof is divided into three parts: uniqueness of martingale solutions, the strong Markov property, and the weighted Krylov-type estimate.

We first bring the following equivalence between time-dependent and time-independent martingale formulations, which comes from \cite[Theorem 1.1]{S75}.
\begin{lem}
\label{lem:equivalent-time-dependent-martingale-problem}
Let $0\leq s<T<\infty$, $x\in\mathbb{R}^d$, and let $\hat{\mathbb{P}}$ be a probability measure on $\mathbb{D}([s,T];\mathbb{R}^d)$ satisfying $\mathbb{P}(\mathcal{X}_s=x)=1$.
Then the following  are equivalent;

(i) For every $\varphi\in C_c^\infty(\mathbb{R}^d)$, the process
\begin{equation*}
M_t^\varphi:=\varphi(\mathcal{X}_t)-\varphi(x)-\int_s^t\mathcal{L}_{\phi,r}^a\varphi(\mathcal{X}_r)\,\mathrm{d}r,
\end{equation*}
is a $(\mathcal{D}_t)_{t\in[s,T]}$-martingale under $\hat{\mathbb{P}}$.

(ii) For every $F \in C^{1,2}_{b}([s,\infty)\times \mathbb{R}^{d})$, the process
\begin{equation*}
N_t^F:=F(t,\mathcal{X}_t)-F(s,x)-
\int_s^t
\left[
\partial_rF(r,\mathcal{X}_r)
+
\mathcal{L}_{\phi,r}^a
\bigl(F(r,\cdot)\bigr)(\mathcal{X}_r)
\right]
\,\mathrm{d}r,
\end{equation*}
is a $(\mathcal{D}_t)_{t\in[s,T]}$-martingale under $\hat{\mathbb{P}}$.

\end{lem}

\begin{lem}\label{lem:stability-approximate-martingale-problems}
Let $0\leq s<T<\infty$ and $x\in\mathbb{R}^d$, and let $\{\mathbb{P}_n\}_{n=1}^{\infty}$ be a tight family of probability measures on $\mathbb{D}([s,T];\mathbb{R}^{d})$ satisfying $\mathbb{P}_n(\mathcal{X}_s=x)=1$ for every $n\in\mathbb{N}$. We suppose the following;

\textbf{A1.} For each $\varphi\in C_c^\infty(\mathbb{R}^d)$, suppose that there exists a jointly Borel measurable function
$$
\Gamma_n^\varphi:
[s,T]\times\mathbb D([s,T];\mathbb R^d)
\longrightarrow
\mathbb R
$$
such that $\omega\mapsto\Gamma_n^\varphi(r,\omega)$ is $\mathcal{D}_r^0$-measurable for every $r\in[s,T]$, and the process
\begin{equation*}\label{eq:approximate-martingale-problem}
M_t^{n,\varphi}:=\varphi(\mathcal{X}_t)-\varphi(x)-\int_s^t\Gamma_n^\varphi(r,\mathcal{X}_{r})\,\mathrm{d}r,
\qquad
t\in[s,T],
\end{equation*}
is a $(\mathcal{D}_t)_{t\in[s,T]}$-martingale under $\mathbb{P}_n$.

\textbf{A2.} Suppose, moreover, that
\begin{equation}
\label{eq:approximate-generator-error}
\lim_{n\to\infty} \mathbb{E}^{\mathbb{P}_n}
\int_s^T
\left|
\Gamma_n^\varphi(r,\mathcal{X})-\mathcal{L}_{\phi,r}^a\varphi(r,\mathcal{X}_r)
\right|
\,\mathrm{d}r = 0
\end{equation}
for every $\varphi\in C_c^\infty(\mathbb{R}^d)$.

Then every weak accumulation point of $\{\mathbb{P}_n\}_{n=1}^{\infty}$ solves the time-dependent martingale problem for $(\mathcal{L}_{\phi,t}^a,s,x)$ on $[s,T]$.
\end{lem}

\begin{proof}
Let $\hat{\mathbb{P}}$ be a weak accumulation point.
Passing to a subsequence, still denoted by $\{\mathbb{P}_n\}$, we assume that $\mathbb{P}_n$ weakly converges to $\hat{\mathbb{P}}$ on the measurable space $\mathbb{D}([s,T];\mathbb{R}^d)$. Since the  mapping $e_{s}: \omega \mapsto \omega(s)$ is continuous with respect to the Skorokhod $J_1$-topology, weak convergence implies
\begin{align*}
\hat{\mathbb{P}}(\mathcal{X}_{s}=x) = \hat{\mathbb{P}}(e_{s}^{-1}(\{x\})) \geq \liminf_{n\to\infty} \mathbb{P}_{n}(e_{s}^{-1}(\{x\})) = \liminf_{n\to\infty}\mathbb{P}_{n}(\mathcal{X}_{s}=x)=1.
\end{align*}
Fix $\varphi\in C_c^\infty(\mathbb{R}^d)$.
For $s\leq u<t\leq T$, define
$$
I_{u,t}^\varphi(\omega)
:=
\int_u^t
\mathcal{L}^{a}_{\phi,r}\varphi(\omega(r)) \,\mathrm{d}r.
$$
Since $\varphi \in C^{\infty}_{c}(\mathbb{R}^{d})$, one can check that $(r,z)\mapsto\mathcal{L}_{\phi,r}^a\varphi(z)$ is jointly Borel measurable, and
\begin{equation}
\label{eq:stability-generator-bound}
\sup_{r\in[s,T]}\sup_{z\in\mathbb{R}^d}|\mathcal{L}_{\phi,r}^a\varphi(z)|<\infty.
\end{equation}
One can also check that there exists a nondecreasing function $\nu_\varphi:[0,\infty)\longrightarrow[0,\infty)$ satisfying $\lim_{\rho\downarrow0}\nu_\varphi(\rho)=0$ such that
\begin{equation}
\label{eq:stability-generator-modulus}
|\mathcal{L}_{\phi,r}^a\varphi(z)-\mathcal{L}_{\phi,r}^a\varphi(z')|
\leq \nu_\varphi(|z-z'|) \quad \forall\, r\in[s,T], \quad \forall\, z,z'\in \mathbb{R}^{d}.
\end{equation}
Suppose that $\omega_m \to \omega$ in the Skorokhod $J_1$-topology.
Then for any $r$, where $\omega\in \mathbb{D}([s,T];\mathbb{R}^{d})$ is continuous, $\omega_m(r) \to \omega(r)$.
Since c\`adl\`ag functions have at most countably many discontinuities, $\omega_{m}(r) \to \omega(r)$ for almost every $r \in [u,t]$.
By \eqref{eq:stability-generator-modulus}, $\mathcal{L}_{\phi,r}^a\varphi(r,\omega_m(r))$ converges to $\mathcal{L}_{\phi,r}^a\varphi(r,\omega(r))$ for almost every $r\in[u,t]$.
Therefore, $I^{\varphi}_{u,t}(\omega_{m}) \to I^{\varphi}_{u,t}(\omega)$ by \eqref{eq:stability-generator-bound} and the dominated convergence theorem.
Hence, $I_{u,t}^\varphi$ is bounded and continuous on $\mathbb{D}([s,T];\mathbb{R}^d)$ equipped with the Skorokhod $J_1$-topology.

Let us denote
\begin{align*}
\Delta\mathcal{X}_{r} = \mathcal{X}_{r} - \lim_{\ell \uparrow r}\mathcal{X}_{\ell}, \quad \mathbb{T}_{\hat{\mathbb{P}}}
:=
\left\{
r\in(s,T):
\hat{\mathbb{P}}(\Delta \mathcal{X}_{r}=0)=1
\right\}
\cup\{s,T\}.
\end{align*}
Then $[s,T]\setminus \mathbb{T}_{\hat{\mathbb{P}}}$ is at most countable (see \textit{e.g.} \cite[Lemma 7.7, Chapter 3]{EK09}). Hence, $\mathbb{T}_{\hat{\mathbb{P}}}$ is dense in $[s,T]$.
Fix $s\leq t_1\leq\cdots\leq t_k\leq u<t\leq T$, where $t_1,\ldots,t_k,u,t\in\mathbb T_{\hat{\mathbb{P}}}$, and let $h_1,\ldots,h_k$ be bounded continuous functions on $\mathbb{R}^{d}$.
Set
$$
H(\omega):=\prod_{j=1}^k h_j(\mathcal \omega(t_{j})).
$$
Since $H$ is bounded and $\mathcal{D}_u$-measurable, \textbf{A1.} gives
\begin{align}
\label{eq:stability-martingale-before-limit}
0=\mathbb{E}^{\mathbb{P}_n}
\Bigg[
\bigg(
\varphi(\mathcal{X}_t)
-
\varphi(\mathcal{X}_u)
-
\int_u^t
\Gamma_n^\varphi(r,\mathcal{X})
\,\mathrm{d}r
\bigg)H\Bigg].
\end{align}
By \eqref{eq:approximate-generator-error},
\begin{equation*}\label{eqn 09.18.16:56}
\begin{aligned}
\left|
\mathbb{E}^{\mathbb{P}_n}
\Bigg[
\bigg(
\int_u^t
\Gamma_n^\varphi(r,\mathcal{X})
\,\mathrm{d}r
-
I_{u,t}^\varphi(\mathcal{X})
\bigg)
H\Bigg]
\right|&\leq
\prod_{j=1}^k\|h_j\|_{L_\infty}
\,
\mathbb{E}^{\mathbb{P}_n}
\int_u^t
\left|
\Gamma_n^\varphi(r,\mathcal{X})
-\mathcal{L}_{\phi,r}^a\varphi(r,\mathcal{X}_r)
\right|\,\mathrm{d}r\\
&\to 0 \quad \text{as} \quad n\to\infty.
\end{aligned}
\end{equation*}
Define
\begin{align*}
\mathcal{H}(\omega)
:=
\left[
\varphi(\omega(t))
-
\varphi(\omega(u))
-
I_{u,t}^\varphi(\omega)
\right]H(\omega)
\end{align*}
Since $I^{\varphi}_{u,t}$ is bounded, $\mathcal{H}$ is also bounded.
Moreover, $\mathcal{H}$ is continuous at every $\omega \in \mathbb{D}([s,T]; \mathbb{R}^{d})$  which is continuous at $t_1,\dots,t_k,u,t$.
Since we take $t_1,\dots,t_k,u,t \in \mathbb{T}_{\hat{\mathbb{P}}}$, the functional $\mathcal{H}$ is continuous ($\hat{\mathbb{P}}$-a.s.).
Therefore, the weak convergence yields $
\mathbb{E}^{\mathbb{P}_n}[\mathcal{H}] \to \mathbb{E}^{\hat{\mathbb{P}}}[\mathcal{H}]$.
Therefore, using \eqref{eq:approximate-generator-error}, and \eqref{eq:stability-martingale-before-limit},  it follows that
\begin{equation}\label{eq:stability-martingale-after-limit}
\begin{aligned}
& \mathbb{E}^{\hat{\mathbb{P}}}\left[  \left( \varphi(\mathcal{X}_{t}) - \varphi(\mathcal{X}_{u}) - \int_{u}^{t} \mathcal{L}^{a}_{\phi,r}\varphi(\mathcal{X}_{r}) \mathrm{d}r  \right) H \right]
\\
& = \mathbb{E}^{\hat{\mathbb{P}}}[\mathcal{H}]
 = \lim_{n\to\infty}\mathbb{E}^{\mathbb{P}_{n}}[\mathcal{H}]
\\
& = \lim_{n\to\infty}\mathbb{E}^{\mathbb{P}_{n}}\left[ \left( \varphi(\mathcal{X}_{t}) - \varphi(\mathcal{X}_{u}) - \int_{u}^{t} \Gamma^{\varphi}_{n}(r,\mathcal{X}) \mathrm{d}r  \right) H  \right] 
\\
& \quad + \lim_{n\to\infty}\mathbb{E}^{\mathbb{P}_{n}}\left[ \left( \int_{u}^{t} \left( \Gamma^{\varphi}_{n}(r,\mathcal{X}) - \mathcal{L}^{a}_{\phi,r}\varphi(r,\mathcal{X}_{r}) \right) \mathrm{d}r \right) H \right]
\\
&=0.
\end{aligned}
\end{equation}
For fixed $u,t\in\mathbb T_{\hat{\mathbb{P}}}$ with $u<t$, monotone class theorem extends (see also \textit{e.g.} \cite[Proposition 7.1, Chapter 3]{EK09} \eqref{eq:stability-martingale-after-limit} so that
\begin{align}\label{eqn 09.22.18:14}
\mathbb{E}^{\hat{\mathbb{P}}}\left[  \left( \varphi(\mathcal{X}_{t}) - \varphi(\mathcal{X}_{u}) - \int_{u}^{t} \mathcal{L}^{a}_{\phi,r}\varphi(\mathcal{X}_{r}) \mathrm{d}r  \right) H \right] = 0 \quad \forall \, H:\text{ bounded $\mathcal{D}^{0}_{u}$-measurable}
\end{align}
follows. Now, for convenience, denote
$$
M_t^\varphi
:=
\varphi(\mathcal X_t)-\varphi(x)
-
\int_s^t
\mathcal L_{\phi,r}^a
\varphi(\mathcal X_r)
\,\mathrm{d}r.
$$
The process $M^\varphi$ is c\`adl\`ag and satisfies
\begin{equation}
\label{eq:stability-limit-martingale-bound}
\sup_{s\leq r\leq T}
|M_r^\varphi|
\leq
2\|\varphi\|_{L_\infty}
+
(T-s)
\sup_{r,z}
|\mathcal{L}_{\phi,r}^a\varphi(z)|<\infty.
\end{equation}
Let $s\leq u<t\leq T$ be arbitrary.
Choose $u_m,t_m\in\mathbb{T}_{\hat{\mathbb{P}}}$ such that $u_m\downarrow u$, $t_m\downarrow t$ and $u_m<t_m$.
If $t=T$, we take $t_m=T$.
For every bounded $\mathcal{D}_u^0$-measurable random variable $H$, we have $H\in\mathcal{D}_{u_m}^0$, and hence 
$$
\mathbb{E}^{\hat{\mathbb{P}}}
\left[
\left(
M_{t_m}^\varphi-M_{u_m}^\varphi
\right)H
\right]
=0.
$$
due to \eqref{eqn 09.22.18:14}. By the right-continuity of $M^\varphi$,
\eqref{eq:stability-limit-martingale-bound}, and the dominated
convergence theorem,
$$
\mathbb E^{\hat{\mathbb{P}}}
\left[
\left(
M_t^\varphi-M_u^\varphi
\right)H
\right]
=0.
$$
Thus, $M^\varphi$ is a $(\mathcal{D}^{0}_{t})_{t\in[s,T]}$-martingale. Finally, let $H$ be bounded and $\mathcal D_u$-measurable.
Choose $u<v_m<t$, and $v_m\downarrow u$.
Since $\mathcal{D}_u\subset\mathcal{D}_{v_m}^0$, it follows that
$$
\mathbb{E}^{\hat{\mathbb{P}}}
\left[
\left(
M_t^\varphi-M_{v_m}^\varphi
\right)H
\right]
=0.
$$
Letting $m\to\infty$, we conclude that
$$
\mathbb{E}^{\hat{\mathbb{P}}}
\left[
\left(
M_t^\varphi-M_u^\varphi
\right)H
\right]
=0.
$$
Therefore, $M^\varphi$ is a $(\mathcal{D}_{t})_{t\in[s,T]}$-martingale. Therefore, using Lemma \ref{lem:equivalent-time-dependent-martingale-problem}, we deduce that $\hat{\mathbb{P}}$ solves the martingale problem for $(\mathcal{L}_{\phi,t}^a,s,x)$ on $[s,T]$.
\end{proof}

\begin{rem}\label{rmk 09.29.17:29}
Under the assumptions of Lemma \ref{lem:stability-approximate-martingale-problems}, suppose in addition that the martingale problem for $(\mathcal{L}_{\phi,t}^a,s,x)$ on $[s,T]$ has at most one solution.

Since $\{\mathbb{P}_n\}_{n=1}^{\infty}$ is tight, Prokhorov's theorem implies that it has a weakly convergent subsequence.
By Lemma \ref{lem:stability-approximate-martingale-problems}, every weak accumulation point is a solution to the martingale problem. 
Hence, the martingale problem has a unique solution, which we denote by $\mathbb{P}_{s,T,x}$.

Let $\{\mathbb{P}_{n_k}\}_{k=1}^{\infty}$ be an arbitrary subsequence. 
By tightness, it admits a further subsequence $\{\mathbb{P}_{n_{k_\ell}}\}_{\ell=1}^{\infty}$ converging weakly to some probability measure $\mathbb{Q}$. 
Lemma \ref{lem:stability-approximate-martingale-problems} implies  that $\mathbb{Q}$ solves the martingale problem for $(\mathcal{L}_{\phi,t}^a,s,x)$. 
By uniqueness, $\mathbb{Q}=\mathbb{P}_{s,T,x}$.
Thus, every subsequence of $\{\mathbb{P}_n\}_{n=1}^{\infty}$ has a further subsequence converging to $\mathbb{P}_{s,T,x}$. Therefore, $\mathbb{P}_{n}$ weakly converges to $\mathbb{P}_{s,T,x}$.
\end{rem}

\begin{prop}[Existence of solutions to time-dependent martingale problems]
\label{lem:existence-martingale-solution}
For every $s\geq0$ and $x\in\mathbb{R}^d$, there exists a probability measure $\mathbb{P}_{s,x}$ such that, for every $T>s$, its restriction to $\mathbb{D}([s,T];\mathbb{R}^d)$ is a solution to the time-dependent martingale problem for $(\mathcal{L}_{\phi,t}^a,s,x)$ on $[s,T]$.
\end{prop}
\begin{proof}
\textbf{Step 1.} We construct a sequence of c\`adl\`ag adapted processes.
\\
Fix $0\leq s<T<\infty$ and $x\in\mathbb{R}^d$.
Let $\pi_{n}$ be a sequence of partitions of $[s,T]$ such that 
$$
\pi_n=\{s=t_0^n<t_1^n<\cdots<t_{N_n}^n=T\}, \quad \lim_{n\to\infty}|\pi_{n}|:=\lim_{n\to\infty}\sup_{k}|t^{n}_{k+1}-t^{n}_{k}|=0
$$
Define $\kappa_n(r):=t_k^n$ for $r\in(t_k^n,t_{k+1}^n]$. Let $N(\mathrm{d}r,\mathrm{d}y,\mathrm{d}v)$ be a Poisson random measure on $[s,T]\times\mathbb{R}^d\times[0,\infty)$ with intensity $\mathrm{d}r\,j^{\phi}_{d}(y)\mathrm{d}y\,\mathrm{d}v$, and let
$$
\widetilde N(\mathrm{d}r,\mathrm{d}y,\mathrm{d}v)
:=
N(\mathrm{d}r,\mathrm{d}y,\mathrm{d}v)
-
\mathrm{d}r\,j^{\phi}_{d}(y)\mathrm{d}y\,\mathrm{d}v.
$$
Construct recursively a c\`adl\`ag adapted process $X^n$ satisfying
\begin{align*}
X_t^n
&=
x+
\int_s^t
\int_{|y|\leq1}
\int_0^\infty
yH_r^n(y,v)\,
\widetilde N(\mathrm{d}r,\mathrm{d}y,\mathrm{d}v)
+
\int_s^t
\int_{|y|>1}
\int_0^\infty
yH_r^n(y,v)\,
N(\mathrm{d}r,\mathrm{d}y,\mathrm{d}v)
\\
& \quad + \int_{s}^{t} b(r,Z^{n}_{r}) \mathrm{d}r,
\end{align*}
where $H_r^n(y,v):=\mathbf{1}_{[0,a(r,Z_r^n,y)]}(v)$, and  $Z_r^n:=X_{\kappa_n(r)}^n$, and 
\begin{align*}
b(r,z) = \int_{\mathbb{R}^{d}} y \left( \mathbf{1}_{|y|\leq 1} - \chi_{\delta_{0},\delta_{1}} \right) a(r,z,y) j^{\phi}_{d}(y)\mathrm{d}y.
\end{align*}
For $z,\xi\in\mathbb{R}^d$, and $\varphi\in C^{\infty}_{c}(\mathbb{R}^{d})$, set
\begin{equation*}
\mathcal{A}_{r,z}\varphi(\xi)
:=
\int_{\mathbb{R}^d}
\left( \varphi(\xi+y) - \varphi(\xi) -\nabla \varphi(\xi) \cdot y  \chi_{\delta_{0},\delta_{1}} \right)
a(r,z,y)\,
j^{\phi}_{d}(y)\mathrm{d}y.
\end{equation*}
Applying It\^o's formula to $\varphi(X_t^n)$ gives
\begin{align*}
\varphi(X_t^n)-\varphi(x)
=
\int_s^t
\mathcal{A}_{r,Z_r^n}
\varphi(X_{r-}^n)\,
\mathrm{d}r+\mathcal{M}_t^{n,\varphi},
\end{align*}
where
\begin{align*}
\mathcal{M}_t^{n,\varphi}
:=
\int_s^t
\int_{\mathbb{R}^d}
\int_0^\infty
\bigl(
\varphi(X_{r-}^n+y)-\varphi(X_{r-}^n)
\bigr)
H_r^n(y,v)\,
\widetilde N(\mathrm{d}r,\mathrm{d}y,\mathrm{d}v).
\end{align*}
The integrand is predictable, and 
\begin{align*}
\mathbb E
\left[
|\mathcal{M}_T^{n,\varphi}|^2
\right]
\leq
a_1(T-s)
\Bigg[
\|\nabla\varphi\|_{L_\infty}^2
\int_{|y|\leq1}
|y|^2\,j^{\phi}_{d}(y)\mathrm{d}y
+
4\|\varphi\|_{L_\infty}^2
\int_{|y|>1}
j^{\phi}_{d}(y)\mathrm{d}y
\Bigg]
<\infty.
\end{align*}
Thus, $\mathcal{M}^{n,\varphi}$ is a square-integrable martingale. Therefore, the process
\begin{equation*}
M_t^{n,\varphi}
:=
\varphi(X_t^n)-\varphi(x)
-
\int_s^t
\mathcal A_{r,Z_r^n}
\varphi(X_r^n)\,
\mathrm{d}r
\end{equation*}
is a square-integrable martingale since $X_{r-}^n=X_r^n$ for almost every $r$.
\\
\textbf{Step 2.} We prove that the family of probability measures $\{\mathbb{P}_{s,T,x}\}^{\infty}_{n=1}=\{\mathbb{P}\circ (X^{n})^{-1}\}^{\infty}_{n=1}$ is tight on $\mathbb{D}([s,T];\mathbb{R}^{d})$.

Let $s \leq \tau_n \leq T$ be a stopping time, let $0\leq\theta \leq T-s$, and set $\sigma_n:=(\tau_n+\theta)\wedge T$.
For simplicity, we denote that
\begin{align*}
&X^{n}_{t} 
\\
&= x + \int_{s}^{t}\int_{|y|\leq 1}\int_{0}^{\infty} yH^{n}_{r}(y,v)\widetilde{N}(\mathrm{d}r,\mathrm{d}y,\mathrm{d}v) + \int_{s}^{t}\int_{|y|>1}\int_{0}^{\infty} yH^{n}_{r}(y,v) N(\mathrm{d}r,\mathrm{d}y,\mathrm{d}v) 
\\
& \quad + \int_{s}^{t} b(r,Z^{n}_{r}) \mathrm{d}r
\\
&:=x+S^{n}_{t}+J^{n}_{t} + D^{n}_{t}.
\end{align*}
For $S^{n}$, observe that
\begin{align}\label{eqn 10.01.15:08}
\mathbb{E}
\left[
|S_{\sigma_n}^n-S_{\tau_n}^n|^2
\right]=
\mathbb{E}
\int_{\tau_n}^{\sigma_n}
\int_{|y|\leq1}
\int_0^\infty
|y|^2H_r^n(y,v)
\,\mathrm{d}v\,j^{\phi}_{d}(y)\mathrm{d}y\,\mathrm{d}r\leq
C \theta.
\end{align}
For every Borel set $A\subset\{y\in\mathbb R^d:|y|>1\}$ define
$$
\mathcal{N}_t^{n,A}
:=
\int_s^t
\int_A
\int_0^\infty
H_r^n(y,v)\,
N(\mathrm{d}r,\mathrm{d}y,\mathrm{d}v), \quad \left( \mathcal{N}^{n,\{|y|>1\}}_{t}=\mathcal{N}^{n,>1}_{t} \right).
$$
Its predictable compensator is
$$
\Lambda_t^{n,A}
:=
\int_s^t
\int_A
a(r,Z_r^n,y)\,
j^{\phi}_{d}(y)\mathrm{d}y\,\mathrm{d}r, \quad \left( \Lambda^{n,\{|y|>1\}}_{t}=\Lambda^{n,>1}_{t} \right).
$$
In particular, $\mathcal{N}_t^{n,A}-\Lambda_t^{n,A}$ is a martingale. Observe that
\begin{align}\label{eqn 10.01.15:09}
\mathbb{P}
\left(
|J_{\sigma_n}^n-J_{\tau_n}^n|>\frac{\varepsilon}{2}
\right) &\leq
\mathbb{P}
\left(
\mathcal{N}_{\sigma_n}^{n,>1}
-
\mathcal{N}_{\tau_n}^{n,>1}
\geq1
\right)    \nonumber
\\
& \leq
\mathbb{E}
\left[
\mathcal{N}_{\sigma_n}^{n,>1}
-
\mathcal{N}_{\tau_n}^{n,>1}
\right]=
\mathbb{E}
\left[
\Lambda_{\sigma_n}^{n,>1}
-
\Lambda_{\tau_n}^{n,>1}
\right]   \nonumber
\\
&\leq
C\theta.
\end{align}
Indeed, the first inequality holds since $J^{n}_{\sigma_{n}}=J^{n}_{\tau_{n}}$ if there is no large jump in the interval $(\tau_n,\sigma_n]$. 
Observe that
\begin{align*}
\sup_{r,z} |b(r,z)| \leq C \left(\int_{0}^{1} \phi(\rho^{-2})\mathrm{d}\rho \mathbf{1}_{\delta_{1}<1/2} + 0 \mathbf{1}_{\delta_{0}\leq 1/2\leq \delta_{1}} + \int_{1}^{\infty} \phi(\rho^{-2}) \mathrm{d}\rho \mathbf{1}_{\delta_{0}>1/2} \right) \leq C
\end{align*}
Therefore,
\begin{align}\label{eqn 11.21.12:30}
\mathbb{E}[|D^{n}_{\sigma_{n}} - D^{n}_{\tau_{n}} |] \leq C \theta.
\end{align}
Combining \eqref{eqn 10.01.15:08}, \eqref{eqn 10.01.15:09}, and \eqref{eqn 11.21.12:30},
\begin{align*}
\mathbb{P}(|X^{n}_{\sigma_{n}}-X^{n}_{\tau_{n}}| > \varepsilon) &\leq \mathbb{P}\left(|S^{n}_{\sigma_{n}}-S^{n}_{\tau_{n}}| > \frac{\varepsilon}{3}\right) + \mathbb{P}\left(|J^{n}_{\sigma_{n}}-J^{n}_{\tau_{n}}| > \frac{\varepsilon}{3}\right) +\mathbb{P}\left( |D^{n}_{\sigma_{n}} - D^{n}_{\tau_{n}} | \geq \frac{\varepsilon}{3} \right)
\\
& \leq C \frac{\mathbb{E}[|S^{n}_{\sigma_{n}}-S^{n}_{\tau_{n}}|^{2}]}{\varepsilon^{2}} + C\theta + C\theta\frac{1}{\varepsilon} \leq C \left(\frac{1}{\varepsilon^{2}} + \frac{1}{\varepsilon}+ 1 \right) \theta,
\end{align*}
where the constant $C$ is independent of $|\pi_{n}|$, $\tau_n$, and $\theta$. Therefore, it follows that
\begin{equation}\label{eqn 10.01.17:06}
\limsup_{\theta \downarrow 0}
\sup_{n}
\sup_{\tau_{n}}
\mathbb P
\left(
|X_{(\tau_n+\theta)\wedge T}^n-X_{\tau_n}^n|
>\varepsilon
\right)
=0.
\end{equation}

By Doob's inequality,
\begin{equation}
\label{eq:appendix-small-jump-maximal}
\sup_n
\mathbb E
\left[
\sup_{s\leq t\leq T}|S_t^n|^2
\right]
\leq
4 a_{1}(T-s) \int_{|y|\leq 1}|y|^{2}j^{\phi}_{d}(y)\mathrm{d}y = C(T-s).
\end{equation}
Fix $K>1$ and decompose $J_t^n=J_t^{n,1,K}+J_t^{n,>K}$, where
\begin{align*}
J_t^{n,1,K}
&:=
\int_s^t
\int_{1<|y|\leq K}
\int_0^\infty
yH_r^n(y,v)\,
N(\mathrm{d}r,\mathrm{d}y,\mathrm{d}v),
\\
J_t^{n,>K}
&:=
\int_s^t
\int_{|y|>K}
\int_0^\infty
yH_r^n(y,v)\,
N(\mathrm{d}r,\mathrm{d}y,\mathrm{d}v).
\end{align*}
Set $\mathcal{N}_t^{n,1,K}:=\mathcal{N}_t^{n,\{1<|y|\leq K\}}$ and $\mathcal{N}_t^{n,>K}:=\mathcal{N}_t^{n,\{|y|>K\}}$.
Since the number of jumps up to $t$, whose sizes are in $(1,K]$ is $\mathcal{N}^{n,1,K}_{t}$, 
$$
\sup_{s\leq t\leq T}|J_t^{n,1,K}|
\leq
K\mathcal N_T^{n,1,K}.
$$
Thus,
\begin{align}
\mathbb{P}
\left(
\sup_{s\leq t\leq T}|J_t^{n,1,K}|>M
\right)
\leq \frac{\mathbb{E}[\mathcal{N}^{n,1,K}_{T}]}{M}
\leq
\frac{a_1K(T-s)}{M}
\int_{1<|y|\leq K}
j^{\phi}_{d}(y)\mathrm{d}y.
\label{eq:appendix-medium-jump-bound}
\end{align}
Similarly,
\begin{align}
\label{eq:appendix-tail-jump-bound}
\mathbb{P}
\left(
\sup_{s\leq t\leq T}|J_t^{n,>K}|>0
\right)
\leq \mathbb{P}(\mathcal{N}^{n,>K}_{T} \geq 1)
\leq
a_1(T-s)
\int_{|y|>K}
j^{\phi}_{d}(y)\mathrm{d}y.
\end{align}
\begin{align}\label{eqn 11.21.01:17}
\mathbb{P}\left( \sup_{s\leq t\leq T} |D^{n}_{t}| > M \right) \leq C(T-s)/M.
\end{align}
Therefore, for any $M\geq 5|x|$, 
using \eqref{eq:appendix-small-jump-maximal},
\eqref{eq:appendix-medium-jump-bound}, 
\eqref{eq:appendix-tail-jump-bound}, and \eqref{eqn 11.21.01:17},
\begin{align*}
\mathbb{P}
\left(
\sup_{s\leq t\leq T}|X_t^n|>M
\right)\leq&
\mathbb{P}
\left(
\sup_{s\leq t\leq T}|S_t^n|>\frac M5
\right)+
\mathbb{P}
\left(
\sup_{s\leq t\leq T}|J_t^{n,1,K}|>\frac M5
\right)
\\
&+
\mathbb{P}
\left(
\sup_{s\leq t\leq T}|J_t^{n,>K}|>0
\right) + \mathbb{P}\left( \sup_{s\leq t\leq T} |D^{n}_{t}| > \frac{M}{5} \right)
\\
\leq&
\frac{C(T-s)}{M^2}+
\frac{5a_1K(T-s)}{M}
\int_{1<|y|\leq K}j^{\phi}_{d}(y)\mathrm{d}y
\\
&+
a_1(T-s)
\int_{|y|>K}j^{\phi}_{d}(y)\mathrm{d}y+ \frac{C(T-s)}{M}.
\end{align*}
Since
$$
\lim_{K\to\infty} \int_{|y|>K}j^{\phi}_{d}(y)\mathrm{d}y=0,
$$
we conclude that
\begin{equation*}
\lim_{M\to\infty}
\sup_n
\mathbb{P}
\left(
\sup_{s\leq t\leq T}|X_t^n|>M
\right)
=0.
\end{equation*}
Hence, the family $\{X_t^n\}_{n=1}^{\infty}$ is tight for every fixed $t\in[s,T]$. Therefore, $\{\mathbb{P}^{n}_{s,T,x}\}^{\infty}_{n=1}$ satisfies Aldous' tightness criterion \cite[Theorem 1]{Al78} (also recall \eqref{eqn 10.01.17:06}). Thus, $\{\mathbb{P}_{s,T,x}^n\}_{n=1}^{\infty}$ is tight on $\mathbb{D}([s,T];\mathbb{R}^d)$.
\\
\textbf{Step 3.} Now we show that for every $\varphi\in C_c^\infty(\mathbb R^d)$,
\begin{equation}
\label{eq:frozen-generator-consistency}
\lim_{n\to\infty}\mathbb{E}^{\mathbb{P}_{s,T,x}^n}
\int_s^T
\left|
\mathcal{A}_{r,\mathcal{X}_{\kappa_n(r)}}
\varphi(\mathcal{X}_r)
-
\mathcal{L}_{\phi,r}^a
\varphi(\mathcal{X}_r)
\right|
\,\mathrm{d}r
= 0.
\end{equation}
Since $\mathbb{P}^{n}_{s,T,x}=\mathbb{P}\circ (X^{n})^{-1}$, and $\mathcal{X}_{r}(\omega)=\omega(r)$ for any $\omega \in \mathbb{D}([s,T];\mathbb{R}^{d})$, to prove \eqref{eq:frozen-generator-consistency}, it is enough to show
\begin{align}\label{eqn 10.01.18:10}
\lim_{n\to\infty} \mathbb{E}\int_{s}^{T} \left| \mathcal{A}_{r,Z^{n}_{r}}\varphi(X^{n}_{r}) - \mathcal{L}^{a}_{\phi,r}\varphi(X^{n}_{r})  \right| \mathrm{d}r = 0.
\end{align}
Set $G_n(r):=|\mathcal{A}_{r,Z_r^n}\varphi(X_r^n)-\mathcal{L}_{\phi,r}^a\varphi(X_r^n)|$.
Since $\mathcal{L}_{\phi,r}^a\varphi(X_r^n)=\mathcal{A}_{r,X_r^n}\varphi(X_r^n)$, and $\varphi\in C^{\infty}_{c}(\mathbb{R}^{d})$, one can check that
\begin{align*}
G_n(r)\leq C \nu(|Z_r^n-X_r^n|), \quad \sup_{n,r}|G_{n}(r)| \leq C,
\end{align*}
where the constant $C=C(\varphi)$ is independent of $n,r$. Let $\varepsilon>0$. Splitting the range of integral, we get
\begin{align*}
\mathbb{E}
\int_s^T
G_n(r)\,\mathrm{d}r &= \mathbb{E}\int_{s}^{T}\mathbf{1}_{|Z^{n}_{r}-X^{n}_{r}| \leq  \varepsilon} G_{n}(r)\mathrm{d}r  + \mathbb{E}\int_{s}^{T}\mathbf{1}_{|Z^{n}_{r}-X^{n}_{r}| >  \varepsilon} G_{n}(r)\mathrm{d}r
\\
&\leq
C\nu(\varepsilon)+
C
\sup_{r\in[s,T]}
\mathbb P
\left(
|Z^{n}_{r}-X^{n}_{r}|>\varepsilon
\right).
\end{align*}
Since $Z^{n}_{r}=X^{n}_{\kappa_{n}(r)}$, and $0\leq r-\kappa_n(r)\leq|\pi_n|$, the estimate \eqref{eqn 10.01.17:06}, applied with $\tau_n=\kappa_n(r)$ and $\theta=r-\kappa_n(r)$, shows that
$$
\lim_{n\to\infty} \sup_{r\in[s,T]}
\mathbb P
\left(
|X_r^n-X_{\kappa_n(r)}^n|>\varepsilon
\right)
=0.
$$
Hence, it follows that
\begin{align*}
\lim_{n\to\infty} \mathbb E
\int_s^T
G_n(r)\,\mathrm{d}r
=0.
\end{align*}
This implies \eqref{eqn 10.01.18:10}, and the desired assertion is proved.
\\
\textbf{Step 4.}
Choose a sequence of partitions $\pi_n$ of $[s,\infty)$ such that
\begin{align*}
\lim_{n\to\infty} |\pi^{T}_{n}| := \lim_{n\to\infty} |(\pi_{n}\cap [s,T]) \cup \{T\}| = 0 \quad \forall \, s<T<\infty.
\end{align*}
Construct process $X^{n}$ as in \textbf{Step 1}, using $\pi_{n}$, so that the restriction $X^{n}|_{[s,T]}$ fulfills any assertions from \textbf{Step 1} to \textbf{Step 3} on $[s,T]$ with $\pi^{T}_{n}$. Fix $T>s$, and let $\mathbb{P}_{s,T,x}^n:=\mathbb{P}\circ(X^n|_{[s,T]})^{-1}$.
By \textbf{Step 2}, the family $\{\mathbb{P}_{s,T,x}^n\}_{n=1}^{\infty}$ is tight on $\mathbb{D}([s,T];\mathbb{R}^d)$. A diagonal argument yields a subsequence, still denoted by $\{X^n\}_{n=1}^{\infty}$ so that $\mathbb{P}\circ(X^{n}|_{[s,s+k]})$ weakly converges (let us denote the limit as $\mathbb{P}_{s,s+k,x}$) for all $k \in \mathbb{N}$.
Define
\begin{align*}
\Delta\mathcal{X}_{r} = \mathcal{X}_{r} -\lim_{\ell\uparrow r}\mathcal{X}_{\ell}, \quad D_{k}:=\{r \in (s,s+k): \mathbb{P}_{s,s+k,x}(\Delta\mathcal{X}_{r}\neq0)>0\},
\end{align*}
and
\begin{align*}
\mathbb{T}_{\mathbb{P}_{s,x}} := \bigcup_{k=1}^{\infty}\left( (s,s+k)\setminus D_{k} \right).
\end{align*}
Then since $\mathbb{T}_{\mathbb{P}_{s,x}}$ is dense in $[s,\infty)$, we may take a sequence $T_{m} \uparrow \infty$ in $\mathbb{T}_{\mathbb{P}_{s,x}}$. Define $\mathbb{P}_{s,T_{m},x}$ as the restriction of $\mathbb{P}_{s,s+k,x}$ to $[s,T_{m}]$ ($T_{m}<s+k$). 
Note that $\mathbb{P}_{s,T_{m},x}$ is invariant under the choice of $k$, and hence it is well-defined. Indeed, if we define a restriction mapping $R^{T_{2}}_{T_{1}}: \mathcal{X}|_{[s,T_{2}]} \mapsto \mathcal{X}|_{[s,T_{1}]}$ ($T_{1}<T_{2})$, then for any $k<\ell$, and $T_{m}<s+k<s+\ell$,
\begin{align*}
\mathbb{P}_{s,s+k,x}\circ\left( \mathcal{X}|_{[s,T_{m}]}^{-1} \right) &= \mathbb{P}_{s,s+k,x}\circ \left( (R^{s+k}_{T_{m}}(\mathcal{X}))^{-1} \right)
\\
&= \mathbb{P}_{s,s+\ell,x} \circ \left( R^{s+\ell}_{s+k}(R^{s+k}_{T_{m}}(\mathcal{X})) \right)^{-1} = \mathbb{P}_{s,s+\ell,x}\left( R^{s+\ell}_{T_{m}}(\mathcal{X}) \right)^{-1}
\end{align*}
Define $\mathbb{P}_{s,x}$ so that the restriction of $\mathbb{P}_{s,x}$ to $[s,T_{m}]$ is $\mathbb{P}_{s,T_{m},x}$. 
Then one can check that $\mathbb{P}_{s,x}$ is well-defined, and the restriction of $\mathbb{P}_{s,x}$ to $[s,T]$ (as a restriction of $\mathbb{P}_{s,T_{m},x}$ to $[s,T]$ for $T<T_{m}$) is also well-defined for any $s<T<\infty$.  
For every
$\varphi\in C_c^\infty(\mathbb R^d)$, define
\begin{equation*}
\Gamma_n^\varphi(r,\omega)
:=
\mathcal A_{r,\omega(\kappa_n(r))}
\varphi(\omega(r)).
\end{equation*}
The process
\begin{align*}
\varphi(\mathcal X_t)-\varphi(x)
-
\int_s^t
\Gamma_n^\varphi(r,\mathcal{X}_{r})
\,\mathrm{d}r
\end{align*}
is a martingale under
$\mathbb P_{s,T_{m},x}^n$. Furthermore,
\textbf{Step 3} yields
\begin{align*}
\lim_{n\to\infty}\mathbb E^{\mathbb P_{s,T_{m},x}^n}
\int_{s}^{T_{m}}
\left|
\Gamma_n^\varphi(r,\mathcal X)
-
\mathcal L_{\phi,r}^a
\varphi(\mathcal X_r)
\right|
\,\mathrm{d}r
=0.
\end{align*}
Since $\varphi \in C^{\infty}_{c}(\mathbb{R}^{d})$, we apply Lemma \ref{lem:stability-approximate-martingale-problems} to obtain that $\mathbb{P}_{s,T_{m},x}$ solves the time-dependent martingale problem on $[s,T_{m}]$. Finally, due to the construction of $P_{s,T,x}$ for general $T$, we also deduce that $P_{s,T,x}$ solves the time-dependent martingale problem on $[s,T]$ by taking $T<T_{m}$, and the corresponding martingale solution $\mathbb{P}_{s,T_{m},x}$ and using the restriction.
\end{proof}

\subsection{Uniqueness of solutions to time-dependent martingale problems}

We begin with an occupation identity that is independent of the
choice of a martingale solution.

\begin{lem}
\label{lem:occupation-identity-uniqueness}
Take  $\beta>0$  from Assumption \ref{24.08.18.13.17}, and choose
\begin{align*}
\gamma_*:=\min\left\{\frac{1}{2},\frac{1}{4\delta_1}\right\}, \quad q_*>4, \quad p_*>\max\left\{2,\frac{d}{\beta},\frac{2d}{\delta_0\gamma_*}\right\}
\end{align*}
so that
\begin{align*}
2\delta_1\gamma_*<1, \quad \beta-\frac{d}{p_*}>0, \quad 2\delta_0\gamma_*\left(1-\frac{2}{q_*}\right)>\frac{2d}{p_*}.
\end{align*}
Let $0\leq s<T<\infty$ and $f\in C_c^\infty((0,T)\times\mathbb{R}^d)$.
Suppose that $u \in \mathbb{H}^{\phi,2}_{q_{\ast},p_{\ast}}(T)$ is the solution to
\begin{equation}
\label{eq:uniqueness-backward-equation}
\partial_ru
+
\mathcal L_{\phi,r}^au
=
-f \quad \text{on} \quad (0,T)\times\mathbb{R}^{d}, 
\quad
u(T,\cdot)=0,
\end{equation}
given by Theorem \ref{thm 11.24.15:58} after time reversal.
Then every martingale solution $\hat{\mathbb{P}}$ starting from $(s,x)$ satisfies
\begin{equation}
\label{eq:uniqueness-occupation-identity}
u(s,x)
=
\mathbb E^{\hat{\mathbb P}}
\int_s^T
f(r,\mathcal X_r)
\,\mathrm{d}r.
\end{equation}
\end{lem}

\begin{proof}
Theorem \ref{thm 11.24.15:58} gives
\begin{equation}
\label{eq:uniqueness-maximal-regularity}
\|\partial_ru\|_{L_{q_*,p_*}(T)}+\|u\|_{H_{q_*,p_*}^{\phi,2}(T)}\leq C\|f\|_{L_{q_*,p_*}(T)}.
\end{equation}
Due to Theorem \ref{thm 05.19.17:17}, we can assume that $u$ is uniformly continuous in $(s,T)\times \mathbb{R}^{d}$. Let $\eta\in C^{\infty}_{c}(B_{1})$ be a nonnegative function with $\|\eta\|_{L_{p^{\ast}}}=1$, and set
\begin{align*}
\eta_{\varepsilon_{1}}(z) = \varepsilon_{1}^{-d/p_{\ast}}\eta(z/\varepsilon_{1}), \quad g^{\varepsilon_{1}}(r,x) = (g(r,\cdot) \ast \eta^{p_{\ast}}_{\varepsilon_{1}})(x) = \int_{\mathbb{R}^{d}} g(r,x-z) \eta^{p_{\ast}}_{\varepsilon_{1}}(z) \mathrm{d}z.
\end{align*}
Define
\begin{equation*}
R_{\varepsilon_{1}}(r,z)
:=
\mathcal L_{\phi,r}^au^{\varepsilon_{1}}(r,z)
-
(
\mathcal L_{\phi,r}^a u )^{\varepsilon_{1}}(r,z).
\end{equation*}
From \eqref{eq:uniqueness-backward-equation}, one can check that
\begin{equation}
\label{eq:uniqueness-mollified-equation}
\partial_ru^{\varepsilon_{1}}
+
\mathcal L_{\phi,r}^au^{\varepsilon_{1}}
=
-f^{\varepsilon_{1}}+R_{\varepsilon_{1}}
\end{equation}
in the sense of distributions.

For fixed $r$ and $z$, define $\overline{a}_{r,z}(x,y):=a(r,z,y)-a(r,x,y)$.
Applying Lemma \ref{lem 07.17.19:49} with $w\equiv1$ and $p=p_*$, coefficient $\overline a_{r,z}$, and $\eta_{\varepsilon_{1}}(z-\cdot)$, and then applying H\"older's inequality, it follows that
\begin{equation}\label{eq:uniqueness-commutator}
\begin{aligned}
\sup_{z\in\mathbb R^d}
|R_{\varepsilon_{1}}(r,z)|
&\leq
\int_{\mathbb{R}^{d}} \left| \mathcal{L}^{\bar{a}_{r,z}}_{\phi,r}u(r,z-x) \right| \eta^{p_{\ast}}_{\varepsilon_{1}}(x) \mathrm{d}x
\\
&=\int_{\mathbb{R}^{d}} \left| \mathcal{L}^{\bar{a}_{r,z}}_{\phi,r}u(r,z-x)\right| \eta_{\varepsilon_{1}}(x)  \eta^{p_{\ast}-1}_{\varepsilon_{1}}(x) \mathrm{d}z
\\
&\leq \|\mathcal{L}^{\bar{a}_{r,z}}_{\phi,r}u(r,z-\cdot)\eta_{\varepsilon_{1}}\|_{L_{p_{\ast}}} \|\eta\|_{L_{p_{\ast}}}
\leq
C \|\eta\|_{L_{p_{\ast}}} \mathcal{C}_{*}(\varepsilon_{1})
\|\phi(\Delta)u(r,\cdot)\|_{L_{p_*}},
\end{aligned}
\end{equation}
where
\begin{equation*}
\mathcal{C}_{*}(\varepsilon_{1})
:=
{\varepsilon_{1}}^{-d/p_*}\|\eta\|_{L_{\infty}}
\left(
\nu(2\varepsilon_{1})
+
{\varepsilon_{1}}^\beta
\int_0^\varepsilon
\frac{\nu(\rho)}
{\rho^{1+\beta}}
\,\mathrm{d}\rho
\right) + {\varepsilon_{1}}^{-d/p_{*}}\|\nabla \eta\|_{L_{\infty}}\nu (2\varepsilon_{1}). 
\end{equation*}
Using Remark \ref{rmk 07.16.16:08} and $\beta-\frac{d}{p_*}>0$, we conclude that $\mathcal{C}_*(\varepsilon_{1})\rightarrow0$ as $\varepsilon_{1}\downarrow0$.
Hence, from \eqref{eq:uniqueness-maximal-regularity} and \eqref{eq:uniqueness-commutator},
\begin{align}
\label{eq:uniqueness-remainder-path}
\mathbb{E}^{\hat{\mathbb{P}}}
\int_s^T
|R_{\varepsilon_{1}}(r,\mathcal{X}_r)|
\,\mathrm{d}r\leq
C \mathcal{C}_{*}(\varepsilon_{1})
(T-s)^{1-1/q_*}
\|\phi(\Delta)u\|_{L_{q_*,p_*}(T)}
\to 0 \quad \text{as} \quad \varepsilon_{1} \downarrow 0
\end{align}
follows. Since $u$ is uniformly continuous, we extend $u^{\varepsilon_{1}}$ so that
\begin{gather*}
u^{\varepsilon_{1}}(t,x) = u^{\varepsilon_{1}}(s,x) \quad \forall\, t\leq s,
\quad\quad
u^{\varepsilon_{1}}(t,x) = u^{\varepsilon_{1}}(T,x) \quad \forall\, t\geq T.
\end{gather*}
Now, take a nonnegative function $\tilde{\eta} \in C^{\infty}_{c}((0,1))$, and define
\begin{align*}
u^{\varepsilon_{1},\varepsilon_{2}}(t,x) = \int_{-\infty}^{\infty} u^{\varepsilon_{1}}(t-r,x) \varepsilon_{2}^{-1}\tilde{\eta}(r/\varepsilon_{2}) \mathrm{d}r \left( = \int_{0}^{\varepsilon_{2}} u^{\varepsilon_{1}}(t-r,x) \varepsilon_{2}^{-1}\tilde{\eta}(r/\varepsilon_{2}) \mathrm{d}r \right).
\end{align*}
Since $\hat{\mathbb{P}}$ is a martingale solution and $u^{\varepsilon_{1},\varepsilon_{2}} \in C^{1,2}_{b}([s,\infty) \times \mathbb{R}^{d})$, we obtain
\begin{equation}\label{eq:uniqueness-Dynkin}
\begin{aligned}
u^{\varepsilon_{1},\varepsilon_{2}}(s,x) &= \mathbb{E}^{\hat{\mathbb{P}}}[u^{\varepsilon_{1},\varepsilon_{2}}(s,\mathcal{X}_{s})]
\\
&=  \mathbb{E}^{\hat{\mathbb{P}}}[u^{\varepsilon_{1},\varepsilon_{2}}(T,\mathcal{X}_{T})] -  \mathbb{E}^{\hat{\mathbb{P}}} \int_{s}^{T} \left( \partial_{r}u^{\varepsilon_{1},\varepsilon_{2}}(r,\mathcal{X}_{r}) + \mathcal{L}^{a}_{\phi,r}u^{\varepsilon_{1},\varepsilon_{2}}(r,\mathcal{X}_{r})  \right) \mathrm{d}r
\\
&= \mathbb{E}^{\hat{\mathbb{P}}}[u^{\varepsilon_{1},\varepsilon_{2}}(T,\mathcal{X}_{T})]+ \mathbb E^{\hat{\mathbb{P}}}
\int_s^T
f^{\varepsilon_{1}}(r,\mathcal X_r)
\,\mathrm{d}r
-
\mathbb E^{\hat{\mathbb{P}}}
\int_s^T
R_{\varepsilon_{1}}(r,\mathcal X_r)
\,\mathrm{d}r
\\
& \quad +\mathbb{E}^{\hat{\mathbb{P}}}\int_{s}^{T} \left( \partial_{r}+\mathcal{L}^{a}_{\phi,r} \right) \left(  u^{\varepsilon_{1}} - u^{\varepsilon_{1},\varepsilon_{2}}\right)(r,\mathcal{X}_{r}) \mathrm{d}r.
\end{aligned}
\end{equation}
Observe that
\begin{align*}
&  \left| \mathcal{L}^{a}_{\phi,r}u^{\varepsilon_{1}}(r,z) - \mathcal{L}^{a}_{\phi,r}u^{\varepsilon_{1},\varepsilon_{2}}(r,z) \right| 
\\
&  \leq C  |u^{\varepsilon_{1}}(r,\cdot) - u^{\varepsilon_{1},\varepsilon_{2}}(r,\cdot)|_{C^{2}(\mathbb{R}^{d})}  
\\
&  \leq C \varepsilon_{1}^{-2} \int_{-\infty}^{\infty} \left( u^{\varepsilon_{1}}(r-\varepsilon_{2}\rho,z) - u^{\varepsilon_{1}}(r,z) \right) \tilde{\eta}(\rho) \mathrm{d}\rho.
\end{align*}
Since $u^{\varepsilon_{1}}$ is uniformly continuous, we have
\begin{align}\label{eqn 11.18.16:11}
\lim_{\varepsilon_{2}\downarrow 0} \sup_{r,z}\left| \mathcal{L}^{a}_{\phi,r}u^{\varepsilon_{1}}(r,z) - \mathcal{L}^{a}_{\phi,r}u^{\varepsilon_{1},\varepsilon_{2}}(r,z) \right| = 0.
\end{align}
Also, one can check that
\begin{align*}
 |\partial_{r}u^{\varepsilon_{1}}(r,z) - \partial_{r}u^{\varepsilon_{1},\varepsilon_{2}}(r,z)|  
& \leq \int_{0}^{1} \left|  \partial_{r}u^{\varepsilon_{1}}(r-\varepsilon_{2}\rho,z) - \partial_{r}u^{\varepsilon_{1}}(r,z)\right|\tilde{\eta}(\rho) \mathrm{d}\rho
\\
& \leq \int_{0}^{1} \int_{\mathbb{R}^{d}} \left| \partial_{r}u(r-\varepsilon_{2}\rho,z-y)-\partial_ru(r,z-y) \right| \eta_{\varepsilon_{1}}^{p_{\ast}} (y) \mathrm{d}y \tilde{\eta}(\rho)\mathrm{d}\rho
\\
&\leq \int_{0}^{1} \| \partial_{r}u(r-\varepsilon_{2}\rho,\cdot) - \partial_{r}u(r,\cdot)\|_{L_{p_{\ast}}} \tilde{\eta}(\rho)\mathrm{d}\rho.
\end{align*}
Therefore, we have
\begin{equation}\label{eqn 11.18.16:11-2}
\begin{aligned}
&\limsup_{\varepsilon_{2} \downarrow 0}  \int_{s}^{T} \sup_{z\in \mathbb{R}^{d}} \left| \partial_{r}u^{\varepsilon_{1}}(r,z) - \partial_{r}u^{\varepsilon_{1},\varepsilon_{2}}(r,z) \right|  \mathrm{d}r
\\
&\leq \limsup_{\varepsilon_{2} \downarrow 0}C(\varepsilon_{1}) \int_{0}^{1} \int_{s}^{T} \| \partial_{r}u(r-\varepsilon_{2}\rho,\cdot) - \partial_{r}u(r,\cdot)\|_{L_{p_{\ast}}}  \mathrm{d}r  \tilde{\eta}(\rho)\mathrm{d}\rho
\\
&\leq  \limsup_{\varepsilon_{2} \downarrow 0} C(T,s,\varepsilon_{1}) \int_{0}^{1} \| \partial_{r}u(\cdot-\varepsilon_{2}\rho,\cdot) - \partial_{r}u(\cdot,\cdot) \|_{L_{q_{\ast},p_{\ast}}} \tilde{\eta}(\rho) \mathrm{d}\rho = 0,
\end{aligned}
\end{equation}
where for the last inequality, we use the continuity of $L_{q_{\ast},p_{\ast}}$ norm and the fact that $\partial_{r}u \in L_{q_{\ast},p_{\ast}}$. Plugging \eqref{eqn 11.18.16:11}, and \eqref{eqn 11.18.16:11-2} to \eqref{eq:uniqueness-Dynkin} (recall that $u^{\varepsilon_{1},\varepsilon_{2}}$ converges to $u^{\varepsilon_{1}}$ uniformly as $\varepsilon_{2}\downarrow 0$ and $u^{\varepsilon_{1}}(T,\cdot)=0$)
\begin{align}\label{eqn 11.18.16:41}
u^{\varepsilon_{1}}(s,x) = \lim_{\varepsilon_{2}\downarrow0} u^{\varepsilon_{1},\varepsilon_{2}}(s,x)= \mathbb{E}^{\hat{\mathbb{P}}}\int_{s}^{T}f^{\varepsilon_{1}}(r,\mathcal{X}_{r}) \mathrm{d}r - \mathbb{E}^{\hat{\mathbb{P}}}\int_{s}^{T} R_{\varepsilon_{1}}(r,\mathcal{X}_{r})\mathrm{d}r.
\end{align}
Recall that $u(s,\cdot)$ is continuous. Therefore, $u^{\varepsilon_{1}}(s,x)$ converges to $u(s,x)$ as $\varepsilon_{1} \downarrow 0$.
Also, since $f\in C^{\infty}_{c}((s,T)\times \mathbb{R}^{d})$, $f^{\varepsilon_{1}} \to f$ uniformly on $[s,T]\times\mathbb{R}^d$. 
Hence,
$$
\mathbb{E}^{\hat{\mathbb{P}}}
\int_s^T
f^{\varepsilon_{1}}(r,\mathcal{X}_r)
\,\mathrm{d}r
\to 
\mathbb{E}^{\hat{\mathbb{P}}}
\int_s^T
f(r,\mathcal{X}_r)
\,\mathrm{d}r \quad \text{as} \quad \varepsilon_{1} \downarrow 0.
$$
Therefore, letting $\varepsilon_{1}\downarrow0$ in \eqref{eqn 11.18.16:41} and using \eqref{eq:uniqueness-remainder-path}, we have \eqref{eq:uniqueness-occupation-identity}. The lemma is proved.
\end{proof}

\begin{prop}
\label{lem:uniqueness-martingale-solution}
For every $s\geq0$ and $x\in\mathbb R^d$, there is at most one global martingale solution starting from $(s,x)$.
\end{prop}

\begin{proof}
Let $\mathbb{P}^1$ and $\mathbb{P}^2$ be two global martingale solutions starting from $(s,x)$.
Fix $s<T<T_{1}<\infty$, and let $h\in C_c^\infty((s,T_{1}))$, $g\in C_c^\infty(\mathbb{R}^d)$.
We apply Lemma \ref{lem:occupation-identity-uniqueness} with $f(r,z):=h(r)g(z)$ to obtain
\begin{align*}
\int_{s}^{T_{1}}h(t)\mathbb{E}^{\mathbb{P}^1}[g(\mathcal{X}_t)]\,\mathrm{d}t=\int_{s}^{T_{1}}h(t)\mathbb{E}^{\mathbb{P}^2}[g(\mathcal X_t)]\,\mathrm{d}t.
\end{align*}
Consequently, for almost every $t \in (s,T_{1})$,
\begin{equation}
\label{26.08.06.13.41}
\mathbb{E}^{\mathbb{P}^1}[g(\mathcal{X}_t)]=\mathbb{E}^{\mathbb{P}^2}[g(\mathcal{X}_t)]
\end{equation}
holds. Note that for any $t_m\downarrow t$, by the right-continuity of the canonical paths $\mathcal{X}_{t_m}\to\mathcal{X}_t$ almost surely. Hence, for $i=1,2$, the map $t\mapsto\mathbb{E}^{\mathbb{P}^i}[g(\mathcal{X}_t)]$ is right-continuous. Thus, \eqref{26.08.06.13.41} holds for every $t\in[s,T]$.
Since $g\in C^{\infty}_{c}(\mathbb{R}^{d})$ is arbitrary, for any Borel set $\Gamma \subset \mathbb{R}^{d}$
\begin{equation}
\label{eq:one-dimensional-uniqueness}
\mathbb{P}^{1}(\mathcal{X}_{t} \in \Gamma) = \mathbb{P}^{2}(\mathcal{X}_{t} \in \Gamma) \quad \forall\, t \in [s,T]
\end{equation}
follows. Then one can apply a theory for regular conditional probability distributions (\textit{e.g.} \cite[Theorem 1.2]{S75}) and the induction argument to derive (here $\Gamma_{i}$ are Borel sets)
\begin{equation}\label{eqn 07.21.14:53}
\begin{aligned}
&\mathbb{P}^{1}(\mathcal{X}_{t_{1}} \in \Gamma_{1}, \cdots, \mathcal{X}_{t_{n+1}} \in \Gamma_{n+1})
\\
&= \int \left( \mathbf{1}_{\Gamma_{1}}(\mathcal{X}_{t_{1}}) \times \cdots \times \mathbf{1}_{\Gamma_{n}}(\mathcal{X}_{t_{n}}) \right) \mathbb{E}^{1}_{t_{n},\mathcal{X}_{t_{n}}}[\mathbf{1}_{\Gamma_{n+1}}(\mathcal{X}_{t_{n+1}})] \, \mathrm{d}\mathbb{P}^{1}
\\
&= \int \left( \mathbf{1}_{\Gamma_{1}}(\mathcal{X}_{t_{1}}) \times \cdots \times \mathbf{1}_{\Gamma_{n}}(\mathcal{X}_{t_{n}}) \right) \mathbb{E}^{2}_{t_{n},\mathcal{X}_{t_{n}}}[\mathbf{1}_{\Gamma_{n+1}}(\mathcal{X}_{t_{n+1}})] \, \mathrm{d}\mathbb{P}^{1}
\\
&= \int \left( \mathbf{1}_{\Gamma_{1}}(\mathcal{X}_{t_{1}}) \times \cdots \times \mathbf{1}_{\Gamma_{n}}(\mathcal{X}_{t_{n}}) \right) \mathbb{E}^{2}_{t_{n},\mathcal{X}_{t_{n}}}[\mathbf{1}_{\Gamma_{n+1}}(\mathcal{X}_{t_{n+1}})] \, \mathrm{d}\mathbb{P}^{2}
\\
& = \mathbb{P}^{2}(\mathcal{X}_{t_{1}} \in \Gamma_{1}, \cdots, \mathcal{X}_{t_{n+1}} \in \Gamma_{n+1}),
\end{aligned}
\end{equation}
where $\mathbb{E}^{i}_{t_{n},\mathcal{X}_{t_{n}}}$ are expectations with respect to regular conditional probabilities of $\mathbb{P}^{i}$ with respect to $\mathcal{D}_{t_{n}}$, which is a martingale solution starting from $(t_{n},\mathcal{X}_{t_{n}}(\omega))$. Since finite dimensional cylinder sets generate $\mathcal{D}_{T}$, \eqref{eqn 07.21.14:53} yields $\mathbb{P}^{1}=\mathbb{P}^{2}$. Since $s<T<\infty$ is arbitrary, we prove the proposition.
\end{proof}

\subsection{Strong Markov property}

\begin{defn}[Time-inhomogeneous strong Markov family]
\label{def:time-inhomogeneous-strong-Markov-family}
Let
$$
\left\{\mathbb{P}_{s,x}:s\geq0,\ x\in\mathbb{R}^d\right\} = \{ \mathbb{P}_{s,x}\}_{s\geq 0,\, x \in \mathbb{R}^{d}}
$$
be a family of probability measures such that $\mathbb{P}_{s,x}(\mathcal{X}_s=x)=1$.
For $0\leq r\leq t$ and every bounded Borel function $g:\mathbb{R}^d\to\mathbb{R}$, define
$$
P_{r,t}g(z):=\mathbb{E}_{r,z}\left[g(\mathcal{X}_t)\right].
$$
Assume that, for every $\vartheta\geq0$ and every bounded Borel function $g$, the map
$$
(r,z)\longmapsto P_{r,r+\vartheta}g(z)
$$
is Borel measurable.
For a stopping time $\tau$, let
$$
\mathcal{D}_\tau:=\left\{A\in\mathcal{D}_\infty:A\cap\{\tau\leq t\}\in\mathcal{D}_t\text{ for every }t\geq0\right\}.
$$
We say that $\{\mathbb{P}_{s,x}\}_{s\geq0,\,x\in\mathbb R^d}$ is  a 
\emph{time-inhomogeneous strong Markov family} if
\begin{equation}
\label{eq:def-time-inhomogeneous-strong-Markov}
\mathbb E_{s,x}
\left[
g(\mathcal X_{\tau+r})
\mid
\mathcal D_\tau
\right]
=
P_{\tau,\tau+r}g(\mathcal X_\tau)
\quad
\mathbb P_{s,x}\text{-a.s.}
\end{equation}
holds for every $s\geq0$, $x\in\mathbb{R}^d$, every finite (a.s.) stopping time $\tau\geq s$, every $r\geq0$, and every bounded Borel function $g:\mathbb{R}^d\to\mathbb{R}$.
\end{defn}

\begin{prop}[Strong Markov property]
\label{lem:strong-Markov-martingale-solution}
Let $\mathbb{P}_{s,x}$ be the unique global martingale solution given by Propositions \ref{lem:existence-martingale-solution} and \ref{lem:uniqueness-martingale-solution}.
Then
$$
\left\{\mathbb{P}_{s,x}:s\geq0,\ x\in\mathbb{R}^d\right\}
$$
is a time-inhomogeneous strong Markov family.
\end{prop}

\begin{proof}
For any $s \geq 0$, and $x \in \mathbb{R}^{d}$, one can apply \textit{e.g.} \cite[Chapter 4]{EK09} to deduce that the process $(\mathcal{X}_{t})_{t\in [s,T]}$ is a strong Markov process (under the filtration $(\mathcal{D}_{t})_{t \in [s,T]}$). Hence, for any finite (a.s.) stopping time $\tau \geq s$, bounded Borel function $g$, and for any $r\geq0$, strong Markov property yields \eqref{eq:def-time-inhomogeneous-strong-Markov} , and thus $\{\mathbb{P}_{s,x} \}_{s\geq0,x\in \mathbb{R}^{d}}$ is a time-inhomogeneous strong Markov family.
\end{proof}

\subsection{Krylov-type estimate}

\begin{prop}[Krylov-type estimate]
\label{lem:Krylov-estimate-without-beta-condition}
Assume the hypotheses of Theorem \ref{thm 02.26.19:05} $(iii)$.
Then, for every $0\leq s<T<\infty$, $x\in\mathbb{R}^d$, and $f\in C_c^\infty((0,T)\times\mathbb{R}^d)$, estimate~\eqref{eqn 01.26.15:17} holds.
\end{prop}

\begin{proof}
Let $ \eta \in C_c^\infty(\mathbb R^d)$ be a nonnegative function with unit integral such that $\eta(x)=0$ for $|x| \geq 1$. Define
$$
\eta_{n}(z):=n^d\eta(nz),\quad a_n(t,x,y):=\int_{\mathbb{R}^d}a(t,x-z,y)\eta_n(z)\,\mathrm{d}z.
$$
The coefficients $a_n$ satisfy Assumption \ref{asm 11.10.16:43} with constants independent of $n$.
In particular,
$$
|a_n(t,x,y)-a_n(t,x',y)|
\leq
\nu(|x-x'|)
$$
and
\begin{equation*}
\sup_{t,x,y}
|a_n(t,x,y)-a(t,x,y)|
\leq
\nu(n^{-1})
\to 0 \quad \text{as} \quad n\to\infty.
\end{equation*}
Note that for each $n$, the coefficient $a_n$ is infinitely differentiable in $x$.

Let $\mathbb{P}_{s,x}^{(n)}$ denote the unique martingale solution associated with $\mathcal{L}_{\phi,t}^{a_n}$.
Let $u_n^{w}\in\mathbb{H}_{q,p}^{\phi,2}(w_2,w_1,T)$ be the solution to
\begin{equation}\label{eqn 09.10.18:00}
\partial_rv+\mathcal{L}_{\phi,r}^{a_n}v=-f,\qquad v(T,\cdot)=0.
\end{equation}
By Theorem \ref{thm 11.24.15:58},
\begin{align}\label{eqn 09.10.18:13}
\|\partial_ru_n^{w}\|_{L_{q,p}(w_2,w_1,T)}+\|u_n^{w}\|_{H_{q,p}^{\phi,2}(w_2,w_1,T)}\leq C\|f\|_{L_{q,p}(w_2,w_1,T)},
\end{align}
where the constant $C$ is independent of $n$.

Let $u_n^*\in\mathbb{H}_{q_*,p_*}^{\phi,2}(T)$ be the solution to \eqref{eqn 09.10.18:00}, where $p_*,q_*,\gamma_*$ are chosen as in
Lemma \ref{lem:occupation-identity-uniqueness}.
By Remark \ref{26.08.06.14.07}, $u_n^{w}=u_n^*$ in the sense of distribution. Moreover, due to Theorem \ref{thm 05.19.17:17}, we can assume that both $u_{n}^{w}$ and $u^{\ast}_{n}$ are continuous near $(s,x)$. These allow us to denote their common value at $(s,x)$ by $u_n(s,x)$. Applying
Lemma \ref{lem:occupation-identity-uniqueness} to the coefficient
$a_n$ and the martingale solution $\mathbb P_{s,x}^{(n)}$, we
obtain
\begin{equation}
\label{eq:Krylov-occupation-an}
u_n(s,x) = u^{\ast}_{n}(s,x) = \mathbb{E}^{\mathbb{P}^{(n)}_{s,x}}
\int_s^T
f(r,\mathcal{X}_r)
\,\mathrm{d}r.
\end{equation}
By Proposition \ref{prop 05.19.17:17} (with $\rho=\rho_{s,T}$), \eqref{eqn 09.10.18:13}, and \eqref{eq:Krylov-occupation-an},
\begin{equation}
\label{eq:Krylov-uniform-an}
\begin{aligned}
\left|
\mathbb{E}^{\mathbb{P}^{(n)}_{s,x}}
\int_s^T
f(r,\mathcal X_r)
\,\mathrm{d}r
\right|&=|u_n(s,x)| = |u^{w}_{n}(s,x)|
\\
&\leq
C
\frac{
\|\partial_ru_n\|_{L_{q,p}(w_2,w_1,T)}
+
\|u_n\|_{H_{q,p}^{\phi,2}(w_2,w_1,T)}
}{
[w_2(\mathbb{I}_{\phi,\gamma}(s,\rho))]^{1/q}
[w_1(B_\rho(x))]^{1/p}
}
\\
&\leq
C
\frac{
\|f\|_{L_{q,p}(w_2,w_1,T)}
}{
[w_2(\mathbb{I}_{\phi,\gamma}(s,\rho))]^{1/q}
[w_1(B_\rho(x))]^{1/p}
}.
\end{aligned}
\end{equation}
The constant $C$ is independent of $n$.

We next pass to the limit.
For each $n\in\mathbb N$, let $\mathbb{P}_{s,T,x}^{(n,m)}$ denote the sequence of probability measures on $\mathbb{D}([s,T];\mathbb{R}^{d})$ constructed as in Proposition \ref{lem:existence-martingale-solution} (from \textbf{Step 1} to \textbf{Step 3}) associated with $a_{n}$.
Since the coefficients $a_n$ satisfy all the assumptions that $a$ satisfies uniformly in $n$, one can check that
$$
\left\{
\mathbb{P}_{s,T,x}^{(n,m)}
:
n,m\in\mathbb{N}
\right\}
$$
is tight by following the proof of Proposition \ref{lem:existence-martingale-solution}.
For each fixed $n$, Lemma \ref{lem:stability-approximate-martingale-problems} (recall Remark \ref{rmk 09.29.17:29}) implies that $\mathbb{P}_{s,T,x}^{(n,m)}
$ weakly converges to $\mathbb{P}^{(n)}_{s,T,x}$  as $m\to\infty$. Therefore, the family $\{\mathbb{P}^{(n)}_{s,T,x}\}^{\infty}_{n = 1}$ is also tight on $\mathbb{D}([s,T];\mathbb{R}^{d})$.
For every $\varphi\in C_c^\infty(\mathbb{R}^d)$, set $\Gamma_n^\varphi(r,\omega):=\mathcal{L}_{\phi,r}^{a_n} \varphi(\omega(r))$.
Then
$$
\varphi(\mathcal X_t)-\varphi(x)
-
\int_s^t
\Gamma_n^\varphi(r,\mathcal X_{r})
\,\mathrm{d}r
$$
is a martingale under $\mathbb P_{s,T,x}^{(n)}$ since it is a martingale solution.
Moreover, due to the dominated convergence theorem (recall \eqref{eq:Krylov-occupation-an}),
\begin{align*}
& \lim_{n\to\infty}\mathbb{E}^{\mathbb{P}_{s,x}^{(n)}}\left[
\int_s^T
\left|\Gamma_n^\varphi(r,\mathcal{X})-\mathcal{L}_{\phi,r}^{a}
\varphi(\mathcal{X}_r)\right|
\,\mathrm{d}r\right] = 0.
\end{align*}
Hence, we apply Lemma \ref{lem:stability-approximate-martingale-problems}, to deduce that every weak accumulation point of $\{\mathbb{P}_{s,T,x}^{(n)}\}^{\infty}_{n = 1}$ solves the martingale problem for $\mathcal{L}_{\phi,t}^a$.
Therefore, due to the uniqueness of martingale solution, $\mathbb{P}^{(n)}_{s,T,x}$ weakly converges to $\mathbb{P}_{s,T,x}$.

Since $f$ is bounded and continuous, the functional
$$
\omega \mapsto \int_s^Tf(r,\omega(r))\,\mathrm{d}r
$$
is bounded and continuous with respect to the Skorokhod $J_1$-topology.
Therefore,
\begin{align*}
\lim_{n\to\infty}\mathbb{E}^{\mathbb{P}^{(n)}_{s,x}}
\int_s^T
f(r,\mathcal{X}_r)
\,\mathrm{d}r
=
\mathbb{E}_{s,x}
\int_s^T
f(r,\mathcal{X}_r)
\,\mathrm{d}r.
\end{align*}
Passing to the limit in \eqref{eq:Krylov-uniform-an} proves \eqref{eqn 01.26.15:17}. The Proposition is proved.
\end{proof}

\begin{proof}[Proof of Theorem \ref{thm 02.26.19:05}]
(i) Global existence follows from Proposition \ref{lem:existence-martingale-solution}, while uniqueness follows from Proposition \ref{lem:uniqueness-martingale-solution}.
This proves the first assertion.

(ii) This follows from Proposition \ref{lem:strong-Markov-martingale-solution}.

(iii) We first observe that Proposition \ref{lem:Krylov-estimate-without-beta-condition} guarantees \eqref{eqn 01.26.15:17} when
$f \in C_c^\infty((0,T)\times\mathbb{R}^d)$. Moreover, it also implies that the linear functional
\begin{align*}
\mathcal{K}_{s,T,x}(f) = \mathbb{E}_{s,x}\int_{s}^{T}f(r,\mathcal{X}_{r})\mathrm{d}r
\end{align*}
is bounded on the dense subset of $L_{q,p}(w_2,w_1,T)$. Therefore, by extending $\mathcal{K}_{s,T,x}$ on $L_{q,p}(w_2,w_1,T)$, we prove the assertion.
\end{proof}

\section*{Acknowledgements}
\pdfbookmark{Acknowledgements}{Acknowledgements}
J.-H. Choi was supported by the National Research Foundation of Korea (NRF) grant funded by the Korea government (MSIT) (No. RS-2026-25477288). J. Kang was supported by the National Research Foundation of Korea (NRF) grant funded by the Korea government (MSIT) (No. RS-2025-24533342). D. Park was supported by the National Research Foundation of Korea (NRF) grant funded by the Korea government (MSIT) (No. RS-2025-25415913). D. Park was supported by the National Research Foundation of Korea (NRF) grant funded by the Korea government (MSIT) (No. RS-2026-25600180).
\vspace{1em}

\section*{Declaration of interest}
The authors declare no competing interests.

\section*{Declaration of generative AI and AI-assisted technologies in the manuscript preparation process}
During the revision of this manuscript, the authors used OpenAI Codex and ChatGPT to assist with language editing and organization, searches for related literature, and checks of mathematical arguments, notation, and internal consistency.
The authors reviewed and revised the AI-assisted material and take full responsibility for the manuscript, including its mathematical results, proofs, and references.


\end{document}